\documentclass{article}
\usepackage{arxiv}

\usepackage{amsfonts}
\usepackage{algorithm, algorithmic} 

\usepackage{amssymb}
\usepackage{amsmath}
\usepackage{latexsym}
\usepackage{mathrsfs}
\usepackage{dsfont}
\usepackage{theorem}

\usepackage{epstopdf} 
\usepackage{psfrag}
\usepackage{graphicx}
\usepackage{multirow}

\usepackage[disable]{todonotes}

\newcommand{\GP}{\mathcal{G}_{k}^b}

\newcommand{\GPvar}{\mathcal{V}_{k}}
\newcommand{\expect}[1]{\mathbb{E}[#1]}

\newcommand{\var}[1]{\mathbb{V}[#1]}
\newcommand{\cov}[2]{\mathbb{C}\text{ov}[#1, #2]}

\newcommand{\Rb}{R_{k}^b}
\newcommand{\calRb}{\mathcal{R}_{k}^b}
\newcommand{\tildeRb}{\tilde{R}_{k}^b}

\newcommand{\vect}[1]{\mathbf{#1}}
\newcommand{\mat}[1]{\mathbf{#1}}
\newcommand{\vectRb}{\vect{R}}
\newcommand{\vectbootRb}{\vect{\hat{R}}}
\newcommand{\vectcalRb}{\vect{\mathcal{R}}}

\newtheorem{example}{Example}[section]

\newcommand{\SNOWPACnumber}{4}

\newcommand{\argmin}{\operatornamewithlimits{arg\ min}}

\usepackage[utf8]{inputenc} 
\usepackage[T1]{fontenc}    
\usepackage{hyperref}       
\usepackage{url}            
\usepackage{booktabs}       
\usepackage{nicefrac}       
\usepackage{microtype}      
\usepackage[square,numbers,sort&compress]{natbib}
\usepackage{doi}

\title{A trust-region method for derivative-free nonlinear constrained stochastic optimization 
}
  
\author{ 
	Friedrich Menhorn \\
	Department of Informatics\\
	Technical University of Munich\\
	85748 Garching, Germany \\
	\texttt{menhorn@in.tum.de} \\
	\And
	Florian Augustin \\
	The MathWorks Inc\\
	Natick, MA 01760, USA \\
	\texttt{fmaugust@mit.edu} \\
	\And
	Hans-Joachim Bungartz \\
	Department of Informatics\\
	Technical University of Munich\\
	85748 Garching, Germany \\
	\texttt{bungartz@in.tum.de} \\
	\And
	Youssef M. Marzouk \\
	Department of Aeronautics and Astronautics\\
	Massachusetts Institute of Technology\\
	Cambridge, MA 02139, USA \\
	\texttt{ymarz@mit.edu} \\
}

\begin{document}

\maketitle

\begin{abstract}
In this work we introduce the stochastic nonlinear constrained derivative-free optimization method (S)NOWPAC (Stochastic Nonlinear Optimization With Path-Augmented Constraints). The method extends the derivative-free optimizer NOWPAC \cite{Augustin2014} to be applicable for optimization under uncertainty. It is based on a trust-region framework, utilizing local fully quadratic surrogate models combined with Gaussian process surrogates to mitigate the noise in the objective function and constraint evaluations. We show the performance of our algorithm on a variety of robust optimization problems from the CUTEst benchmark suite by comparing to other popular optimization methods. 
While we focus on robust optimization benchmark problems to demonstrate (S)NOWPAC's capabilities, the optimizer can be applied to a broad range of applications in nonlinear constrained stochastic optimization.
\end{abstract}


\section{Introduction}
Consider a wind park with multiple wind turbines where we are able to control the overall power production by steering the turbine heads via their yaw angle. At the given location we experience very uncertain conditions with respect to, e.g., the direction or the velocity of the wind. Additionally, the sensors for the yaw angle of the wind turbines introduce some measurement error or are delayed in their measurements. Furthermore, we also have to consider that the wind turbines influence each other by wind shadowing and wake effects and we have to satisfy certain structural constraints on the turbine. Finally, given these uncertain conditions and constraints we are interested in steering the turbines optimally to maximize the average total power production and, at the same time, reduce the variance in the result. Hence, we want to find an optimal design which is robust with respect to these uncertain and possibly extreme conditions. Besides, the simulation code that models such kind of scenario is only given as a black box code to us, e.g. it is only an binary executable where we can set the input via setting files and access results through output files. Thus, we cannot assume that we have access to the code or quantities like higher order derivatives or even just gradients. In this paper, we present a novel approach for optimization under uncertainty, i.e., for nonlinear constrained stochastic optimization, designed for such kind of application problems. The approach is derivative-free such that we are able to work with black box applications where we do not have access to additional information like gradients.

We consider an objective function $f$ (, e.g., the total power production), and we derive an algorithm for finding optimal design parameters $x$ within a set $X = \{x \;:\; c(x) \leq 0\} \subseteq \mathbb{R}^n$ of admissible feasible design parameters (, e.g., the optimal yaw angle of the wind turbine). The functions $c = (c_1, \ldots, c_r)$ are called constraint functions (, e.g., the turbine's structural constraints). The objective function and constraints are derived from stochastic models of the underlying process and are therefore stochastic themselve and in general considered to be nonlinear and black box, cf. \cite{Bertsimas2011, Pflug1996}. For example, these uncertain conditions may reflect limited accuracy in measurements as in our example in the beginning; or, it may model our lack of knowledge about process parameters and result in stochastic/non-deterministic evaluations of $f$ and $c$. Though we will later on parametrize this uncertainty by a parameter $\theta$ where $\theta$ follows a probability density function $\Theta$ for easier illustration we point out that we do not require knowledge of the underlying distribution.
Since gradient estimation in black box stochastic optimization is challenging or even impossible, we focus on methodologies that do not require gradients, only utilizing black box evaluations of $f$ and the constraints $c$---so-called \textit{derivative-free} approaches. 

Note that henceforth in this work we use the terms optimization under uncertainty, stochastic optimization and robust optimization interchangeably. Our approach is designed for nonlinear constrained optimization problems that show a stochastic behaviour and which can be formulated in the form of an expectation as explained in the upcoming sections. This is the case for most problems in these areas of optimization.

We now define our general problem formulation that we target in this work. As mentioned before we assume uncertain conditions on $f$ and/or $c$ and have to take this variability into account during optimization.
We regard a stochastic optimization problem \cite{Ben-Tal1999, Beyer2007} using {\it robustness measures} $\mathcal{R}^f$ and $\mathcal{R}^c$,
\begin{equation}\label{eq:general_robust_opt_problem}
\begin{split}
&\;\;\;\min \mathcal{R}^f(x)\\
&\mbox{s.t.} \quad \mathcal{R}^c(x) \leq 0.
\end{split}
\end{equation}
where we are interested in finding a robust or reliable solution with respect to the underlying uncertainties. The choice of $\mathcal{R}$ depends on the problem at hand and we state a variety of robustness measures in Section \ref{sec:robustness_formulation} but also refer to the rich literature on risk and deviation measures \cite{Acerbi2002, Artzner1999, Krokhomal2011, Rackwitz2001, Rockafellar2000, Rockafellar2002a, Rockafellar2002, Uryasev2000, Zhang}. 
In order to simplify notation we omit the superscripts $f$ and $c$ subsequently whenever the reference is clear form the context.

In this work we develop a new approach to find a solution for~\eqref{eq:general_robust_opt_problem} by generalizing the trust-region optimization algorithm NOWPAC \cite{Augustin2014} to an algorithm for optimization under uncertainty with inherently noisy evaluations of the objective function and the constraints. Here, we enhance the trust-region management of stochastic derivative-free optimization procedures for noisy evaluations procedures
building up on work by Kannan and Wild \cite{Kannan2012}. The main contributions of this work are threefold. First, we introduce Gaussian process models of the objective function and the constraints to progressively reduce the noise in function evaluations. Second, we also leverage knowledge from the Gaussian process to
progressively improve evaluations of $f$ and $c$. Combined, this allows us to control the structural error in the local surrogate models used in NOWPAC. Third, we propose a procedure to recover feasibility, which is important in steps where the algorithm falsely assumed feasibility due to the noise in the constraint evaluations. We refer to Section~\ref{sec:derivative_free_stochastic} for a detailed discussion of our contributions.

\section{Other work}
Before diving into the description of our method we give a brief history and of existing optimization techniques to solve~\eqref{eq:general_robust_opt_problem} in the following paragraphs. We particularly point to a recent review paper \cite{Larson2019} for a comprehensive summary of derivative-free optimization methods. 

One possible optimization approach is Sample Average Approximation (SAA) \cite{Kim2011}, where a set of samples $\{\theta_i\}_{i=1}^N \subset \Theta$ is chosen for approximating the robustness measures $R_N \approx \mathcal{R}$. This set of samples is fixed throughout the optimization process to minimize the sample approximated objective function $\mathcal{R}^f$. This results in approximate solutions of (\ref{eq:general_robust_opt_problem}) that depend on the particular choice of samples used. In order to reduce the associated approximation error, typically several optimization runs are averaged or the sample size $N \rightarrow \infty$ is increased; see \cite{Ahmed2008, Rubinstein1993, Shapiro2009}.  An error analysis of SAA for constrained optimization problems can be found in \cite{Bayraksan2006}. The advantage of SAA is that it eliminates the noise introduced by different sample approximations between optimization steps and thus deterministic black box optimization methods can be used to solve the optimization problem.

Other approaches draw new samples from the uncertain parameter $\theta$ every time the robustness measures are evaluated. Due to the re-sampling, the evaluations of the approximate robustness measures $R_N(x)$ exhibit {\it sampling noise} 
and thus solving (\ref{eq:general_robust_opt_problem}) requires stochastic optimization methods.
If the noise is small enough, for example if the sample size $N$ is sufficiently large, pattern search methods may be used to solve the optimization problem. Avoiding gradient approximations makes these methods less sensitive to noise in the evaluations of the robust objective and constraints. Since the early works by Hookes and Jeeves \cite{Hooke1961} and Nelder and Mead \cite{Nelder1965, Spendley1962}, there has been a significant research effort in various extensions and developments of excellent direct search optimization procedures \cite{Audet2008, Audet2006, Audet2009, DiPillo2013, Liuzzi2006, Liuzzi2010, Powell1994,  Powell1998}.
Surrogate model based optimization \cite{Augustin2014, Bortz1998, Conn1993, Kelley1999, March2012, Regis2011, Regis2014, Sampaio2015} can also be used to solve  (\ref{eq:general_robust_opt_problem}). Having sufficiently accurate approximations of gradients even convergence results for these methods exist; see \cite{Carter1991, Choi2000, Heinkenschloss2002, Larson2016, Chen2018}.
Here, {\it sufficiently accurate}, however, requires the gradient approximation to become increasingly accurate while approaching an optimal solution. This idea is incorporated in the proposed derivative-free stochastic optimization procedures STRONG \cite{Chang2013} and ASTRO-DF \cite{Shashaani2015}, which which reduce the noise in black box evaluations by taking averages over an increasing number of samples while approaching an optimal design.

Thus far we only discussed optimization methods that rely on a diminishing magnitude of the noise in the robustness measure approximations and we now turn our attention to methods without this requirement. In 1951, Robbins and Monroe \cite{Robbins1951} pioneered by proposing the Stochastic Approximation (SA) method. Since then SA has been generalized to various gradient approximation schemes, e.g. by Kiefer and Wolfowitz (KWSA) \cite{Kiefer1952} and Spall \cite{Spall1992, Spall1998, Wang2003} with the Simultaneous Perturbation Stochastic Approximation (SPSA). We refer to \cite{Bhatnagar2013, Kibzun1996, Kushner1997} for a detailed introduction and theoretical analysis of SA methods and only remark that for all SA methods a variety of technical parameters, like the step and stencil sizes, have to be chosen very carefully. Despite a rich literature and theoretical results, this choice remains a challenging task in applying SA approaches: optimal and heuristic choices exist \cite{Spall1998}, however, they are highly problem dependent and have a strong influence on the performance and efficiency of SA methods. Nevertheless, with the field of machine learning growing fast as the main application field, variants of the stochastic approximation method see many new developments---specifically stochastic gradient descent. We refer to \cite{Bottou2018} for a recent review.

Finally, Bayesian Global Optimization (BGO) \cite{Mockus1974, Mockus1989} can be used to solve (\ref{eq:general_robust_opt_problem}). In BGO the objective function is approximated using a Gaussian process in order to devise an exploration and exploitation scheme for global optimization based on expected improvement or knowledge gradients, see for example \cite{ Frazier2009, Jones1998}. Handling nonlinear constraints only recently gained attention in BGO \cite{Gramacy2015}. One particular approach, constrained Bayesian Optimization (cBO), based on expected constrained improvement optimization can be found in \cite{Gardner2014} and recent developments, e.g., in \cite{Letham2019}.

We recognized a gap in the existing approaches with respect to the combination of stochastic constrained optimization problems and derivative-free approaches which we intend to close in this paper.
In Section~\ref{sec:robustness_formulation} we introduce sampling based approximations of robustness measures along with their confidence intervals for statistical estimation of their sampling errors. We also discuss a variety of robustness measures to rigorously define the robust formulation \eqref{eq:general_robust_opt_problem}. Thereafter, in Section~\ref{sec:derivative_free_stochastic}, we briefly recap the trust-region algorithm NOWPAC \cite{Augustin2014} which we then generalize to make it applicable to stochastic (noisy) robust optimization tasks. We close with numerical examples in Section~\ref{sec:numerical_results} and conclude in Section~\ref{sec:conclusions}.

\section{Problem formulations for optimization under uncertainty}\label{sec:robustness_formulation}
In this section, we discuss sampling approximations of robustness measures which are most commonly used in stochastic optimization problems along with their associated confidence intervals. We introduce a collection of robustness measures $\mathcal{R}^f$ and $\mathcal{R}^{c}$ to model robustness and risk for the robust optimization problem (\ref{eq:general_robust_opt_problem}). To simplify the notation we refer to the objective function $f$ and the constraints $c$ as black box $b$, the corresponding robustness measures will be denoted by $\mathcal{R}^b$. We assume that $b$ is square integrable with respect to $\theta$, i.e. its variance is finite, its cumulative distribution function is assumed to be continuous and invertible at every fixed design point $x \in \mathbb{R}^n$.\\

\subsection{Statistical estimation of robustness measures}\label{sec:robustness_sample_approximation}
Most robustness measures in literature can be written in terms of an expectation,
\begin{equation}\label{eq:unified_notation_for_robustness_measures}
\mathcal{R}^b(x) := \mathbb{E}_{\theta}\left[ B(x, \theta) \right],
\end{equation}
where the function $B$ depends on the actual choice of robustness measure. Throughout this paper we assume that $B$ has finite variance. 

For the approximation of (\ref{eq:unified_notation_for_robustness_measures}) at $x$ we use a sample average $E_N$ based on $N$ samples $\{\theta_i\}_{i=1}^N \sim \mu$,
\begin{equation}\label{eq:sample_approximation_of_expected value}
\mathbb{E}_{\theta}\left[B(x,\theta)\right] = E_N\left[B(x,\theta)\right] + \varepsilon_x = \frac{1}{N}\sum\limits_{i=1}^N B(x, \theta_i)+ \varepsilon_x.
\end{equation}
Here $\varepsilon_x$ represents the error of the sample approximation. From the Central Limit Theorem we know that $\sqrt{N}\varepsilon_x$ is asymptotically normally distributed with zero mean and variance $\sigma^2 = \mathbb{V}_{\theta}[B(x, \theta)]$ for $N \rightarrow \infty$. This allows the definition of a confidence interval around the approximated expected value, $E_N\left[B(x,\theta)\right] $, which contains $\mathbb{E}_{\theta}\left[B(x,\theta)\right]$ with high probability. To get a confidence interval
\[
\left[ 
\mathbb{E}_N\left[ B(x, \theta) \right] - \bar{\varepsilon}_x, \;
\mathbb{E}_N\left[ B(x, \theta) \right] + \bar{\varepsilon}_x
\right] 
\]
that contains $\mathbb{E}_{\theta}\left[ B(x, \theta) \right]$ with a probability exceeding $\nu \in \; ]0, 1[$ we compute the sample estimate $s_N(x)$ of the standard deviation of $\{B(x, \theta_i)\}_{i=1}^N$, 
\begin{equation}\label{eq:variance_estimator_for_sample_average}
s_N(x)^2 = \frac{1}{N-1}\sum\limits_{i=1}^N\left(B(x, \theta_i) - E_N[B(x, \theta)]\right)^2,
\end{equation}
and set
\[
\bar{\varepsilon}_x = \frac{t_{\nu}\, s_N(x)}{\sqrt{N}},
\]
with $t_{\nu}$ being a constant defining the confidence interval with respect to $\nu$. For larger samples sizes this constants reflects the Z-score. 
We choose $t_{\nu} = 2$ in our implementation which yields a confidence level exceeding $0.975$ for a sample size $N \geq 60$.
Finally, in our proposed Algorithm~\hyperref[alg:SNOWPAC]{\SNOWPACnumber} we use $\bar{\varepsilon}_x$ as an indicator for the upper bound on the sampling error $\varepsilon_x \leq \bar{\varepsilon}_x$ with probability exceeding $\nu$. 

\subsection{Common robustness measures}
In the following we will give common examples of robustness measures and their sampling estimators. We refer to \cite{Rockafellar2002, Szego2002} for a detailed discussion about risk assessment strategies. 

The classical first example for a robustness measure is the expected value
\[
\mathcal{R}_0^b(x) := \mathbb{E}_{\theta}\left[ b(x; \theta) \right] = \int\limits_{\Theta} b(x; \theta) d\mu.
\]
Although it may be arguable that the expected value measures robustness with respect to variations in $\theta$, since it does not inform about the spread of $b$ around $\mathcal{R}_0^b(x)$, it is a widely applied measure to handle uncertain parameters in optimization problems. For example, the expected objective value, $\mathcal{R}_{0}^f$, yields a design that performs best on average, whereas $\mathcal{R}_0^c$ specifies feasibility in expectation. 
In order to also account for the spread of realizations of $b$ around $\mathcal{R}_0^b$ for different values of $\theta$ in a statistical sense, justifying the term {\it robustness} measure, a standard deviation term,
\[
\mathcal{R}^b_1(x) := \mathbb{V}^{\frac{1}{2}}_{\theta}\left[ b(x; \theta) \right] = \sqrt{\mathbb{E}_{\theta}\left[ b(x; \theta)^2 \right] - \mathcal{R}^b_0(x)^2},
\]
can be included. We remark that the linear combination
\[
\mathcal{R}^b_2(x) := \gamma c_1 \mathcal{R}^b_0(x) + (1 - \gamma) c_2 \mathcal{R}^b_1(x),
\]
$\gamma \in [0, 1]$, $c_1$, $c_2 > 0$, of $\mathcal{R}_0^b$ and $\mathcal{R}_1^b$ has a natural interpretation in decision making. By minimizing the standard deviation term we gain confidence in the optimal value being well represented by $\mathcal{R}_0^b$. Combining the two goals of objective minimization in expectation and 
the reduction of the spread of possible outcomes, the robustness measure $\mathcal{R}^b_2(x)$ provides a trade off between two possibly contradicting goals. The user's priority in one goal over the other is reflected by the weighting factor $\gamma$. The constants $c_1$ and $c_2$ are required to obtain a proper scaling between $\mathcal{R}_0^b$ and $\mathcal{R}_1^b$. 
Finally we remark that it is well known that $-\mathcal{R}_0^b$ is a coherent risk measure, whereas $\mathcal{R}^b_2$ is not (see \cite{Artzner1999}).

Commonly used robustness measures are probabilistic constraints, also known as chance constraints \cite{Prekopa1970}, Here, a probability level $\beta \in \; ]0, 1[$ is specified up to which the optimal design has to be feasible. The corresponding robustness measure is 
\[
\mathcal{R}^{b,\beta}_3(x) := \mu\left[b(x, \theta) \geq 0 \right] - (1 - \beta) = 
\mathbb{E}_{\theta}\left[ \mathds{1}(b(x, \theta) \geq 0)\right] - (1 - \beta).
\]
Probabilistic constraints are used for economic modelling, for example construction costs of a power plant not exceeding a prescribed budget with probability $\beta$. Another application of probabilistic constraints in physics is the adjustment a gas mixture in a combustion chamber and prevent extinction of the flame with (high) probability $\beta$. A penalty for the associated costs or risks for violating the constraints can be included in the objective function. See \cite{Li2011, Li2010, Li2012} for an efficient method for approximating $\mathcal{R}^{b, \beta}_3$. Under the assumption of an invertible cumulative distribution function, $F_{\mu}$, probabilistic constraints can be formulated in terms of quantile functions,
\[
\mathcal{R}^{b,\beta}_4(x) := \min \left\{ \alpha \in \mathbb{R} \, : \, \mu[ b(x, \theta) \leq \alpha] \geq \beta \right\}.
\]
The two formulations are equivalent in the sense that they yield the same set of feasible points:
$
\{ x \in \mathbb{R}^n \; : \; \mathcal{R}^{c,\beta}_3(x) \leq 0 \} =
\{ x \in \mathbb{R}^n \; : \; \mathcal{R}^{c,\beta}_4(x) \leq 0 \}.
$
However, in the appendix~\autoref{appssec: quantilesamplingestimator} we show 
that $\mathcal{R}^{b,\beta}_4$ often exhibits favourable smoothness properties as compared to $\mathcal{R}^{b,\beta}_3$, making it more suitable to model probabilistic constraints in our optimization procedure. We remark that for $b = f$, the robustness measure $\mathcal{R}^{f,\beta}_4$ is also known as Value at Risk (VaR), a widely used non-coherent risk measure in finance applications. Note, if the underlying distribution is (assumed) normal also $\mathcal{R}^b_2(x)$ is often used in practice as a chance constraint where $\gamma c_1 = 1$ and $(1 - \gamma) c_2$ describes the confidence interval.

The Conditional Value at Risk (CVaR) \cite{Acerbi2002, Rockafellar2002} is a coherent extension of $\mathcal{R}^{b,\beta}_4(x)$. It is defined as the conditional expectation of $b$ exceeding the VaR:
\[
\mbox{CVaR}_{\beta}(x) := \frac{1}{1-\beta} \int\limits_{b(x, \theta) \geq \mathcal{R}^{b,\beta}_4(x)} b(x, \theta) d\mu.
\]
Following \cite{Alexander2006, Rockafellar2002} we define the robustness measure 
\[
\mathcal{R}^{b,\beta}_5(x, \gamma) := \gamma + \frac{1}{1-\beta}\mathbb{E}_{\theta}\left[  \max\{ b(x, \theta)-\gamma, 0\} \right],
\]
which allows us to minimize the CVaR without having to compute $\mathcal{R}^{b,\beta}_4$ first as minimizing $\mathcal{R}^b_5$ over the extended feasible domain $X \times \mathbb{R}$ yields
\[
\min\limits_{x \in X} \mbox{CVaR}_{\beta}(x) = \min\limits_{(x, \gamma) \in X \times \mathbb{R}} \mathcal{R}^b_5(x, \gamma).
\]

For completeness we also mention the measure traditionally most closely associated with robust optimization, the worst case formulation:
\[
\mathcal{R}^{b}_6(x) := \max\limits_{\theta \in \Theta}\left\{ b(x, \theta)  \right\}.
\]
It is, however, often computationally challenging to evaluate $\mathcal{R}^b_6$ and only in special cases of simple non-black box functions $b$
it is possible to analytically compute $\mathcal{R}^{b}_6(x)$, which then yields a deterministic optimization problem, see f.e. \cite{Ben-Tal1998, Henrion2012, Kibzun1998, Shapiro2009, Uryasev2000}. $\mathcal{R}^{b}_6(x)$ does no necessarily require knowledge of the probability distribution of $\theta$ which we assume in this work. Therefore, we exclude worst case formulations from our discussion subsequently.

Note that in case of $\mathcal{R}_0^b(x)$, $\mathcal{R}_3^{b,\beta}(x)$ and $\mathcal{R}_5^{b,\beta}(x)$ the assumption of finite variance of $B$ already follows from the assumption that $b$ is square integrable with respect to $\theta$. However, for the variance of $B$ in $\mathcal{R}_1^b(x)$ and $\mathcal{R}_2^b(x)$ to be finite we require the stronger integrability condition of $b^2$ being square integrable. 

\section{Stochastic nonlinear constrained optimization}\label{sec:derivative_free_stochastic}
Using a finite number of samples to approximate the robustness measures from Section ~\ref{sec:robustness_sample_approximation} at every step of the optimization introduces sampling noise $\varepsilon_x$. We propose a stochastic optimization framework based on the black box optimizer NOWPAC \cite{Augustin2014} to solve
\begin{equation}\label{eq:general_approximiated_robust_opt_problem}
\begin{split}
&\;\;\;\min R_N^f(x)\\
&\mbox{s.t.} \quad R_N^c(x) \leq 0,
\end{split}
\end{equation}
for finite sample approximations $R_N^b(x) \approx \mathcal{R}^b(x) + \bar{\varepsilon}_x$ of robustness measures. 
Within the Section~\ref{sec:reviewNOWPAC} we briefly review NOWPAC's key features to set the stage for its generalization to (S)NOWPAC---(Stochastic) Nonlinear Optimization With Path-Augmented Constraints---in Sections~\ref{sec:noise_adapted_trust_region_management}- \ref{sec:algorithm_SNOWPAC}. 

\subsection{Review of the trust-region framework NOWPAC}\label{sec:reviewNOWPAC}
NOWPAC \cite{Augustin2014} is a deterministic\footnote{Therefore, we know that $\mathcal{R} = R$.} derivative-fee trust-region optimization framework that uses black box evaluations to build fully linear (see \cite{Conn2009}) surrogate models $m_k^{\mathcal{R}^f}$ and $m_k^{\mathcal{R}^c}$ of the objective function and the constraints within a neighbourhood of the current design $x_k$ where $k \in \mathbb{N}$ is the current optimization step. This neighbourhood, $\{x \in \mathbb{R}^n \,:\, \|x - x_k\| \leq \rho_k\}$, is called a trust-region with trust-region radius $\rho_k>0$. We use the short-hand notation $\mathcal{R}^c(x) := (\mathcal{R}^{c_1}(x), \ldots, \mathcal{R}^{c_r}(x))$ and define the feasible domain as $X := \{x \in \mathbb{R}^n\,:\, \mathcal{R}^c(x) \leq 0\}$. The optimization is performed as follows: starting from $x_0 \in X$ a sequence of intermediate points $\{x_k\}_k$ is computed by solving the trust-region subproblems
\begin{equation}\label{eq:robust_surrogate_model_NOWPAC}
\begin{split}
&\quad x_k := \argmin m_k^{\mathcal{R}^f}(x)\\
&\mbox{s.t.} \quad x \in X_k, \; \|x-x_k\| \leq \rho_k
\end{split}
\end{equation}
with the approximated feasible domain 
\begin{equation}\label{eq:ibp_augmented_feasible_domain}
X_k := \left\{ x \in \mathbb{R}^n \;:\; m_k^{\mathcal{R}^c}(x) + h_k(x-x_k) \leq 0\right\}.
\end{equation}
The additive offset $h_k$ to the constraints is called the inner boundary path, a convex offset-function to the constraints ensuring convergence of NOWPAC. We refer to \cite{Augustin2014} for more details on the inner boundary path. Having computed $x_k$ NOWPAC only accepts this trial step if it is feasible with respect to the exact constraints $\mathcal{R}^c$, i.e. if $\mathcal{R}^c(x_k) \leq 0$. Otherwise the trust-region radius is reduced and, after having ensured fully linearity of the models $m_k^{\mathcal{R}^f}$ and $m_k^{\mathcal{R}^c}$, a new trial step $x_k$ is computed. 

To assess closeness to a first-order optimal point the criticality measure 
\begin{equation}\label{eq:criticality_measure}
\alpha_k(\rho_k) := \frac{1}{\rho_k} \left| \min\limits_{\substack{x_k+d \in X_k\\ \|d\| \leq \rho_k}} \left\langle g_k^{\mathcal{R}^f}, d\right\rangle\right|
\end{equation}
is used, where $g_k^{\mathcal{R}^f} = \nabla m_k^{\mathcal{R}^f}(x_k)$ is the gradient of the surrogate model of the objective function $\mathcal{R}^f$ at $x_k$. 

To assess acceptance of the trial point and the update of the trust-region the acceptance ratio is computed as
\begin{equation}
\label{eq:acceptance_ratio}
r_k = \frac{\mathcal{R}^f(x_k) - \mathcal{R}^f(x_k + s_k)}{m_k^{\mathcal{R}^f}(x_k) - m_k^{\mathcal{R}^f}(x_k + s_k)}.
\end{equation}
This ratio reflects the truth versus the prediction of the surrogate. Based on the result we accept the point, i.e. $r_k >= \eta_1$, or otherwise reject it. We furthermore adapt the trust-region accordingly
\begin{equation}
\label{eq:update_trustregion}
\rho_{k+1} = \begin{cases}
\rho_k    & \mbox{if} \quad r_k \geq 2\\
\gamma_{inc}\rho_k & \mbox{if} \quad r_k \geq \eta_2 \quad \mbox(and)  \quad r_k < 2\\
\rho_k             & \mbox{if} \quad \eta_1 \leq r_k < \eta_2,\\
\gamma_{dec}\rho_k       & \mbox{if} \quad r_k < \eta_1.\\
\end{cases}
\end{equation}
We recall the simplified algorithm for NOWPAC within Algorithm~\ref{alg:simplified_NOWPAC}.
\begin{algorithm}[!h]
\begin{algorithmic}[1]
\STATE{Construct the initial fully linear models $m_0^{\mathcal{R}^f}(x_0 + s)$, $m_0^{\mathcal{R}^c}(x_0 + s)$, k = 0\;}
\WHILE{$\rho_{k} >= \rho_{min}$}
\STATE{Compute criticality measure $\alpha_k(\rho_k)$ via~\eqref{eq:criticality_measure}\;}
\STATE{\hspace{-\algorithmicindent}// \textit{STEP 0:} Criticality step}
\WHILE{$\alpha_k(\rho_k)$ is too small}
\STATE{Decrease $\rho_k = \omega \rho_k$ and update $m_k^{\mathcal{R}^f}$ and $m_k^{\mathcal{R}^c}$\;}
\ENDWHILE
\STATE{\hspace{-\algorithmicindent}// \textit{STEP 1:} Step calculation}
\STATE{Compute a trial step $s_k$ via~\eqref{eq:ibp_augmented_feasible_domain}\;}
\STATE{\hspace{-\algorithmicindent}// \textit{STEP 2:} Check feasibility of trial point}
\IF{$\mathcal{R}^c(x_k)(x_k + s_k) > 0$}
\STATE{Set $\rho_k = \gamma \rho_k$ and update $m_k^{\mathcal{R}^f}$ and $m_k^{\mathcal{R}^f}$\;}
\STATE{Go to {\tt STEP 0}\;}
\ENDIF
\STATE{\hspace{-\algorithmicindent}// \textit{STEP 3:} Acceptance of trial point and update trust-region}
\STATE{Compute $r_k$ via~\eqref{eq:acceptance_ratio}\;}
\IF{$r_k \geq \eta_0$}
\STATE{Set $x_{k+1} = x_k + s_k$\;}
\STATE{Include $x_{k+1}$ into the node set and update the models to $m_{k+1}^{\mathcal{R}^f}$ and $m_{k+1}^{\mathcal{R}^c}$\;}
\ELSE{
\STATE{Set $x_{k+1} = x_k$, $m_{k+1}^{\mathcal{R}^f} = m_k^{{R}^f}$ and $m_{k+1}^{\mathcal{R}^c} = m_k^{\mathcal{R}^c}$\;}
}\ENDIF
\STATE{Update $\rho_{k+1}$ via~\eqref{eq:update_trustregion}\;}
\STATE{Update $m_{k+1}^{\mathcal{R}^f}$ and $m_{k+1}^{\mathcal{R}^c}$\;}
\STATE{k = k+1\;}
\ENDWHILE
\end{algorithmic}
\caption{Simplified NOWPAC\label{alg:simplified_NOWPAC}}
\end{algorithm}

\subsection{Noise-adapted trust-region managment}\label{sec:noise_adapted_trust_region_management}
The efficiency of Algorithm~\ref{alg:simplified_NOWPAC} depends on the accuracy of the surrogate models $m_k^{\mathcal{R}^b}$ and subsequently our ability to predict a good reduction of the objective function within the subproblem (\ref{eq:robust_surrogate_model_NOWPAC}). It is thus necessary to make Algorithm~\ref{alg:simplified_NOWPAC} robust with respect to the sampling noise in finite sampling approximations of the robustness measures. To achieve this goal, we firstly  introduce a noise-adapted trust-region management to NOWPAC to couple the structural error in the surrogate approximations and the sampling error in the evaluation of $R_N$. Secondly we propose the construction of Gaussian processes to reduce the sampling noise in the finite sample approximation of the robustness measures.  

We know from \cite[Thm. 2.2]{Kannan2012} that fully linear surrogate models being constructed from finite sample approximations using noise corrupted black box evaluations satisfy the error bound
\begin{align}\label{eq:error_bound_on_noisy_surrogate_models}
\begin{split}
\left\| \mathcal{R}^b(x_k + s) - m_k^{{R}^b}(x_k + s) \right\| &\leq \kappa_1\, \rho_k^2,\\
\left\|\nabla \mathcal{R}^b(x_k + s) - \nabla m_k^{{R}^b}(x_k + s) \right\| &\leq \kappa_2\, \rho_k.
\end{split}
\end{align}
The constants $\kappa_1$ and $\kappa_2$ depend on the poisedness constant $\Lambda \geq 1$--- a measure for the spread of the surrogate points in the domain---as well as on the estimates of the statistical upper bounds for the noise term, 
$
\bar{\varepsilon}_{max}^k = \max\limits_{i=1}^{\bar{n}}\{\bar{\varepsilon}_i\}
$ 
from Section~\ref{sec:robustness_sample_approximation}. Hence, if the maximal noise term $\bar{\varepsilon}_{max}^k$ is of order $\rho_k^2$ we know that the bounds in \eqref{eq:error_bound_on_noisy_surrogate_models} apply as shown in \cite[Thm. 2.2]{Kannan2012}. Otherwise, in the presence of noise, i.e. $\bar{\varepsilon}_{max}^k > 0$, the term $\bar{\varepsilon}_{max}^k \rho_k^{-2}$ and thus $\kappa_1$ and $\kappa_2$ grow unboundedly for a shrinking trust-region radius, violating the fully linearity property of $m_k^{R^b}$. 
Thus, in order to ensure fully linearity of the surrogate models, we have to enforce an upper bound on the error term.

\begin{algorithm}[!htb]
\begin{algorithmic}
\STATE{Input:  trust-region factor $a \in \{1, \gamma_{\text{dec}},\, \gamma_{\text{inc}}, \, \omega\}$.\;}
\STATE{Set $\rho_{k+1} =  \max \left\{a\rho_k, \;\lambda_t\sqrt{\bar{\varepsilon}_{\text{max}}^k}\right\}$\;}
\IF{$\rho_{k+1} > \rho_{\text{max}}$}
\STATE{Set $\rho_{k+1} = \rho_{\text{max}}$\;}
\ENDIF
\end{algorithmic}
\caption{Noise adapted updating procedure for trust-region radius. \label{alg:noisy_trust_region_management}}
\end{algorithm}

This noise-adapted trust-region management couples the structural error of the fully linear approximation with the highly probable upper bound on the error in the approximation of the robustness measures. This coupling, however, also prevents the trust-region radii from converging to $0$, therefore limiting the level of accuracy of the surrogate models $m_k^{R^b}$ and thus the accuracy of the optimization result. 

\subsection{Gaussian process supported noise correction}\label{sec:gaussian_process_supported_noise_correction}
In order to increase the accuracy of the optimization result, we need to reduce the magnitude of the noise term $\bar{\varepsilon}_{max}^k$. 
Since the straight-forward solution to increase the number of samples is too costly, we, instead, suggest 
a different strategy: we introduce Gaussian process (GP) surrogates of $\mathcal{R}^b$ by using the finite sample approximations at already evaluated optimization points $\{(x_i, R^b_i)\}_{i=1}^{K}$. Using this second surrogate  we can reduce the error and smoothen the resulting estimator by taking into account more global information by the GP since the GP converges to its target function with a increased number of points. Here, we leverage consistency properties and smooth behavior of the GP, e.g., shown in \cite{Scheuerer2013, Wendland2004, Stuart2016}. This is similar to a contral variate approaches for variance reduction in Monte Carlo \cite{Lemieux2017}. As we will see in the following sections this helps us to smoothen the noisy evaluations and decrease the magnitude of the noise term $\bar{\varepsilon}_{max}^k$. 

\subsubsection{GP construction}
For the construction of the GPs we only take points with a distance smaller than $\zeta_1 \rho_k$ around the current best design point $x_k$ into account, i.e. 
\begin{align}
\label{eq:GP_training_set}
(\mat{X}, \vect{R}^b) = \{ (x_j, R^b_j) : \| x_j - x_k \|_2 \leq \zeta_1 \rho_k, j = 1,..., K \}.
\end{align} 
This focuses our approximation to a localized neighbourhood and we can assume stationarity of the GP surrogates. Note, that we thereby specifically require the point of evaluation $x_k$ to be included in the training set. By default, we use a value of $\zeta_1 = 3$ to incorporate enough global information around the current design \todo{FM: constant}.

The GP estimators for mean and variance employing general training data $(\mat{X}, \vect{y})$ are given as 
\begin{align}
\GP[\vect{y}] := \mathcal{G}^b({x}_k; \mat{X}, \vect{y}) = \vect{k}_{{x}_k \mat{X}} (K^b_{\mat{X} \mat{X}} + \mat{N})^{-1} \vect{y} \\
\GPvar := \mathcal{V}(x_k; \mat{X}) = \vect{k}^b_{{x}_k {x}_k} - \vect{k}^b_{{x}_k \mat{X}} (K^b_{\mat{X} \mat{X}} + \mat{N})^{-1} \vect{k}^b_{\mat{X} {x}_k}
\end{align}
Here, $\vect{k}^b_{{x}_k \mat{X}}$ and  $K^b_{\mat{X} \mat{X}}$ denote the kernel vector and kernel matrix, respectively, evaluated at every pair $(x_p, x_q), x_p, x_q \in X$. Additionally, $\mathbf{N}$ denotes a general noise matrix with noise estimates $\bar{\varepsilon}_{k}^b$ on its diagonal. Note, that an independent GP is built for each of the objective and constraint functions. For more details about GPs we refer to \cite{Rasmussen2006}. In the following we employ the short-hand notation $\GP := \GP[\vect{y}]$ and $\GPvar := \GPvar[\vect{y}]$ if the training set is clear from the context.

\subsubsection{GP smoothing}
\label{sec:GPsmoothing}
In order to reduce the noise in the finite sample approximations $R_k^b$ we balance their contribution with the GP surrogate estimates $\GP[\vectRb]$ as they becomes more and more accurate with an increasing amount of evaluations $\vectRb$ during the optimization procedure:
\begin{align}
\label{eq:adjusted_black_box_evaluations}
\tilde{R}_{k}^{b} &= \gamma_k \GP[\vectRb] +(1-\gamma_k) {R}_k^b.
\end{align}
Here we use a linear combination of the noisy sampling estimates ${R}_k^b$ and the mean estimator of the GP $\GP[\vectRb]$  on the current evaluation $x_k$. The weighting of the two contributions is balanced by $\gamma_k$. 

We are interested in finding the optimal $\gamma_k$ to minimize the error of the new estimator $\tilde{R}_{k}^{b} $.
Under the assumption that $R_{k}^b$ is itself an unbiased estimator, we therefore compute the 
root mean squared error (RMSE)
\begin{equation}
\begin{split}
\label{eq:rmse}
\text{RMSE}(\tildeRb)
&= [\gamma_k \expect{\GP[\Rb] - \Rb}]^2 + \var{ \gamma_k \GP[\Rb] + (1-\gamma) \Rb} \\
&= [\gamma_k ( \GP[\calRb] - \calRb)]^2 + \gamma_k^2 \var{ \GP[\Rb] } + (1-\gamma_k)^2 \var{\Rb} \\
&+ 2 \gamma_k (1-\gamma_k) \cov{ \GP[\Rb]}{\Rb},
\end{split} 
\end{equation} 
of the estimator \eqref{eq:adjusted_black_box_evaluations} depending on $\gamma_k$ and use it as noise estimate
\begin{align}
\tilde{\varepsilon}_{k}^b = t_{\nu} \cdot \min_{\gamma_k} \text{RMSE}(\tildeRb).
\end{align}

Taking the derivative \eqref{eq:rmse} for $\gamma_k$ and setting it equal to 0 gives us its optimal value to minimize the error
\begin{align}
\label{eq:mingamma}
\gamma_k= \frac{\var{\Rb} - \cov{\GP[\Rb]}{\Rb}}{( \GP[\calRb] - \calRb )^2 + \var{\GP[\Rb]} + \var{\Rb} - 2 \cov{ \GP[\Rb]}{\Rb} }
\end{align}
which is subsequently used in \eqref{eq:adjusted_black_box_evaluations}. In the end, the approximations $\tilde{R}_{k}^{b}$ as well as the associated noise level $\tilde{\varepsilon}_k^b$ are used to build the local surrogate models $m_k^{R^f}$ and $m_k^{R^c}$.

Due to the linearity of the Gaussian process mean operator we can compute and approximate the quantities in \eqref{eq:mingamma} in closed form. The variance is given as 
\begin{equation}
\label{eq:vargp}
\begin{split}
\var{\GP[\Rb]} &= \sum_{i=1}^N \var{R^b_i} (\sum_{j=1}^N \vect{k}_{{x}_i {x}_j} ((K_{\mat{X} \mat{X}} + \mathbf{N})^{-1})_{[i, j]})^2  \\
&\approx \sum_{i=1}^N (\frac{\vect{\bar{\varepsilon}_i^b}}{t_{\nu}})^2 (\sum_{j=1}^N \vect{k}_{{x}_i {x}_j} ((K_{\mat{X} \mat{X}} + \mathbf{N})^{-1})_{[i, j]})^2
\end{split}
\end{equation}
while we estimate the covariance by
\begin{equation}
\label{eq:covgpr}
\begin{split}
\cov{\GP[\Rb]}{\Rb} &= \var{\Rb} \sum_{j=1}^N \vect{k}_{{x}_k {x}_j} ((K_{\mat{X} \mat{X}} + \mathbf{N})^{-1})_{[k, j]} \\ 
&\approx (\frac{\bar{\varepsilon}_k^b}{t_{\nu}})^2 \sum_{j=1}^N \vect{k}_{{x}_k {x}_j} ((K_{\mat{X} \mat{X}} + \mathbf{N})^{-1})_{[k, j]}.
\end{split}
\end{equation}
Here, the notation $\mat{X}_{[i, j]}$ denotes the element of $\mat{X}$ at position $[i, j]$. Note that we compute the variance of the GP mean estimator $\var{\GP[\Rb]}$ in~\eqref{eq:vargp} which is not the same as the variance estimate $\GPvar$ of the GP surrogate. Similarly, we compute the covariance between two estimators $\cov{\GP[\Rb]}{\Rb}$ in~\eqref{eq:covgpr}.

Finally, we estimate the term $(\GP[\calRb] - \calRb)$ using a bootstrapping approach (cf. \cite{Efron2016, Efron1994}). Here we use that this term is the bias of $\GP$:
\begin{align}
\label{eq:GPbias}
\text{Bias}(\GP) = \expect{\mathcal{G}_k^b(x_k^{(i)}; \vectRb)} - \calRb = \mathcal{G}_k^b(x_k^{(i)}; \expect{\vectRb}) - \calRb = \GP[\vectcalRb] - \calRb.
\end{align}
Since $\mathcal{G}_k^b(x_k^{(i)}; \vectRb)$ estimates $\calRb$ we approximate \eqref{eq:GPbias} by
\begin{align}
\GP[\vectcalRb] -  \calRb \approx \expect{\GP[\vectbootRb]} -  \GP.
\end{align} 
Here, $\vectbootRb$ describes training sets $\hat{\theta}$ created from resampling with replacement from the original set $\theta$ and recomputing $\vectbootRb$ for all $X$. 

The exactness of the above quantities heavily depends on the approximation quality of the GP as well as the robustness measures $R^b$. Therefore, SNOWPAC provides the option to use a heuristic instead of computing $\gamma_k$ and the resulting noise. For this heuristic, the noise is reduced by using a similar linear combination as in \eqref{eq:adjusted_black_box_evaluations} leveraging the GP variance estimator:
\begin{align}
\label{eq:heuristicsmoothing}
\tilde{\varepsilon}_{k}^b\; &=\;\gamma_k t_{N-1,\nu} \GPvar^{\frac{1}{2}} + (1-\gamma_k)\bar{\varepsilon}_k^b.
\end{align}
The weight factor
$\gamma_k := e^{-\GPvar^{\frac{1}{2}}}$
is chosen to approach $1$ when the GP becomes more and more accurate as indicated by the vanishing variance of the GP approximation. 

\subsubsection{GP error balancing}
By combining the two surrogate models we balance two sources of approximation errors. On the one hand, there is the structural error in the approximation of the local surrogate models,~cf. (\ref{eq:error_bound_on_noisy_surrogate_models}), which is controlled by the size of the trust-region radius. On the other hand, we have the inaccuracy in the GP surrogate itself which is reflected by the variance $\GPvar$ of the GP.
Note that Algorithm~\ref{alg:noisy_trust_region_management} relates these two sources of errors by coupling the size of the trust-region radii to the size of the credible interval through (\ref{eq:adjusted_black_box_evaluations}), only allowing the trust-region radius to decrease if $\GPvar$ becomes small.

Finally, we ensure $\GPvar$ becomes smaller as $x_k$ approaches the optimal design, in three ways:
first, the increasing number of black box evaluations performed by the optimizer during the optimization process helps to increase the quality of the Gaussian process approximation \cite{Scheuerer2013, Wendland2004, Stuart2016}. However, these evaluations may be localized and geometrically not well distributed around the current iterate $x_k$. We therefore, second, draw additional points, 
$\hat{x} \sim \mathcal{N}\left(x_k, \zeta_2 \sqrt{\rho_k} I\right)$, with $\zeta_2 = \frac{3}{10}$ by default \todo{FM: constant},
whenever a trial point is rejected to improve the geometrical distribution of the regression points for the GP surrogates. The rejection of a trial point can happen because it may be infeasible under the current GP-corrected constraint approximation (\ref{eq:adjusted_black_box_evaluations}), or the step is rejected in {\tt STEP 3} in Algorithm~\ref{alg:simplified_NOWPAC}.
Third, in addition to enriching the set of regression points, SNOWPAC re-estimates the GP hyperparameters either after a user-prescribed number of black-box evaluations or after $\lambda_k\cdot n$ consecutive rejected or infeasible trial steps, where $\lambda_k$ is a user prescribed constant. This avoids problems with over-fitting \cite{Rasmussen2006, Cawley2007}. 

\subsection{Relaxed feasibility requirement}\label{sec:relaxed_feasibility_requirement}
An integral part of Algorithm \ref{alg:simplified_NOWPAC} is the feasibility requirement in {\tt STEP 2}. It guarantees feasibility of all intermediate design points $x_k$.  Checking feasibility in the presence of noise, however, is challenging. For example, it might happen that Algorithm~\ref{alg:SNOWPAC} accept an apparently feasible point given the current constraint approximations, which is in fact infeasible. We therefore have to generalize NOWPAC's capabilities to recover from infeasible points by introducing a feasibility restoration mode. The resulting algorithm has two operational modes,
\begin{align*}
\mbox{\bf (M1)} & \quad\mbox{objective minimization and}\\
\mbox{\bf (M2)} & \quad\mbox{feasibility restoration.}
\end{align*}
The algorithm operates in mode (M1) whenever the current point $x_k$ appears to be feasible under the current constraint approximations (\ref{eq:adjusted_black_box_evaluations}), and switches to mode (M2) if $x_k$ becomes infeasible. The switch between modes (M1) and (M2) is implemented by exchanging the underlying trust-region subproblem: 
in mode (M1) the standard subproblem
\begin{equation}\label{eq:M1_trial_point}
\begin{split}
&\;\;\;\;\; \; \; \min m_k^{{\tilde{R}}^f}(x_{k} + s_k),\\
&\mbox{s.t.} \quad \bar{m}_k^{\tilde{R}^{c_i}}(x_{k} + s_k) \leq 0,\ i=1 \ldots r\\
&\quad\quad\; \|s_k\| \leq \rho_k
\end{split}
\end{equation}
is solved to obtain a new trial point $x_{k} + s_k$. Here $\bar{m}_c^{\tilde{R}^{c_i}}$ denote the inner-boundary path augmented models of $\mathcal{R}^{c_i}$ as described in \eqref{eq:ibp_augmented_feasible_domain} using the updated evaluations $\tilde{R}_k^b$ from \eqref{eq:adjusted_black_box_evaluations}. The subproblem
\begin{equation}\label{eq:M1_criticality_measure}
\begin{split}
&\;\;\;\;\; \; \; \min \left\langle g_k^{{\tilde{R}}^f}, s_k\right\rangle,\\
&\mbox{s.t.} \quad \bar{m}_k^{{\tilde{R}}^{c_i}}(x_{k} + s_k) \leq 0,\ i=1 \ldots r\\
&\quad\quad\; \|s_k\| \leq \rho_k
\end{split}
\end{equation}
is used for computation of the criticality measure $\alpha_k$. 

In mode (M2) the subproblem
\begin{equation}\label{eq:M2_trial_point}
\begin{split}
&\;\;\;\;\; \; \; \min \sum\limits_{i\in\mathcal{I}_k}\left(m_k^{{\tilde{R}}^{c_i}}(x_{k} + s_k)^2 + \lambda_g m_k^{{\tilde{R}}^{c_i}}(x_{k} + s_k)\right),\\
&\mbox{s.t.} \quad \bar{m}_k^{{\tilde{R}}^{{c_i}}}(x_{k} + s_k) \leq \tau_i,\ i=1 \ldots r\\
&\quad\quad\; \|s_k\| \leq \rho_k
\end{split}
\end{equation}
is solved for the computation of a new trial point $x_{k} + s_k$, along with
\begin{equation}\label{eq:M2_criticality_measure}
\begin{split}
&\;\;\;\;\; \; \; \min \sum\limits_{i\in \mathcal{I}_k}\left( 2m_k^{{\tilde{R}}^{c_i}}(x_{k}) + \lambda_g \right)\left\langle g_k^{m_k^{{\tilde{R}}^{c_i}}}, s_k\right\rangle,\\
&\mbox{s.t.} \quad \bar{m}_k^{{\tilde{R}}^{c_i}}(x_{k} + s_k) \leq \tau_i,\ i=1 \ldots r\\
&\quad\quad\; \|s_k\| \leq \rho_k
\end{split}
\end{equation}
for computation of the corresponding criticality measure. Here, $\mathcal{I}_k = \{i \; : \; \tilde{R}_{k}^{c_i} > 0, \; i = 1, \ldots, r\}$, denotes the set of violated constraints. The slack variables $\tau := (\tau_1, \ldots, \tau_r)$ are set to $\tau_i = \max\{ \tilde{R}_{k}^{c_i}, \; 0\}$. We introduce the parameter $\lambda_g \geq 0$ in (\ref{eq:M2_trial_point}) and  (\ref{eq:M2_criticality_measure}) to guide the feasibility restoration towards the interior of the feasible domain. By default it is set to $\lambda_g = 10^{-4}$.

The respective mode also affects the acceptance ratio $r_k$. Therefore we adapt the algorithm as shown in Algorithm~\ref{alg:acceptance_ratio_mode}.
\begin{algorithm}
\begin{algorithmic}
\STATE{Input:  Current design $x_k$ and trial point $x_{\text{trial}} = x_k + s_k$.\;}
\IF{Mode (M1)}\STATE{$r_k = \frac{\tilde{R}^f(x_k) - \tilde{R}^f(x_{\text{trial}})}{m_k^{\tilde{R}^f}(x_k) - m_k^{\tilde{R}^f}(x_{\text{trial}})}$.\;}
	\ELSE
	\STATE{$r_k = \frac{\sum\limits_{i\in\mathcal{I}_k}\left(\tilde{R}^{c_i}(x_{k})^2 + \lambda_g \tilde{R}^{c_i}(x_{k})\right) - \sum\limits_{i\in\mathcal{I}_k}\left(\tilde{R}^{c_i}(x_{\text{trial}})^2 + \lambda_g \tilde{R}^{c_i}(x_{\text{trial}})\right)}
	{\sum\limits_{i\in\mathcal{I}_k}\left(m_k^{{\tilde{R}}^{c_i}}(x_{k})^2 + \lambda_g m_k^{{\tilde{R}}^{c_i}}(x_{k})\right) - \sum\limits_{i\in\mathcal{I}_k}\left(m_k^{{\tilde{R}}^{c_i}}(x_{\text{trial}})^2 + \lambda_g m_k^{{\tilde{R}}^{c_i}}(x_{\text{trial}})\right)}$.\;}
	\ENDIF
\end{algorithmic}
\caption{Compute acceptance ratio $r_k$ in normal and feasibility restoration mode.}
\label{alg:acceptance_ratio_mode}
\end{algorithm}

\subsection{The stochastic trust-region algorithm (S)NOWPAC}\label{sec:algorithm_SNOWPAC}
In this section we state the final algorithm of (S)NOWPAC which is summarized in Algorithm~\hyperref[alg:SNOWPAC]{\SNOWPACnumber}. The general procedure follows closely the steps in Algorithm~\ref{alg:simplified_NOWPAC} and includes the generalizations we introduced in Sections~\ref{sec:noise_adapted_trust_region_management}, \ref{sec:gaussian_process_supported_noise_correction} and \ref{sec:relaxed_feasibility_requirement} to handle noisy black box evaluations. A summary of all default values for internal parameters we use in our implementation of (S)NOWPAC is given in Table~\ref{tab:tuning_parameters_for_snowpac}. 
(S)NOWPAC is available under the BSD 2-Clause license on Github\footnote{\url{https://github.com/snowpac/snowpac}}. 

\begin{table}[!htb]
\begin{center}
\caption{Internal parameters of (S)NOWPAC and their default values}\label{tab:tuning_parameters_for_snowpac}
\begin{tabular}{|c|c|c|}\hline
description & parameter & default value\\\hline\hline
factor for lower bound on trust-region radii & $\lambda_{t}$ & $\sqrt{2}$\\\hline
poisedness threshold &$\Lambda$& $100$\\\hline
gradient contribution to feasibility restoration &$\lambda_{g}$& $10^{-4}$\\\hline
factor for GP region & $\zeta_1$ & $3$ \\\hline
constant for normal distribution to enrich GP & $\zeta_2$ & $\frac{3}{10}$ \\\hline
\end{tabular}
\end{center}
\end{table}


\begin{algorithm}[!h]
\begin{algorithmic}
\STATE{Construct the initial fully linear models $m_0^{R^f}(x_0 + s)$, $m_0^{R^c}(x_0 + s)$.\;}
\STATE{Set $x^* = x_0$ and $R^{f*} = R_N^f(x_0)$.\;}
\STATE{Switch Mode if feasibility restoration is necessary.\;} 
\FOR{$k = 0, 1, \ldots, n_{max}$}
\STATE{\hspace{-\algorithmicindent}// \textit{STEP 0:} Criticality step}
\WHILE{$\alpha_k(\rho_k)$ is too small}
    \IF{  $\rho_k < \rho_{min}$ } \STATE{STOP\;}\ENDIF
	\STATE{Call Algorithm \ref{alg:noisy_trust_region_management} with $a = \omega$.\;}
	\STATE{Evaluate $R_N^{f}(x)$ and $R_N^{c}(x)$ for a randomly sampled $x$ in a neighbourhood of $B(x_k, \rho_k)$.\;}
	\STATE{Update Gaussian processes and black box evaluations.\;}
	\STATE{Construct surrogate models $m_k^{\tilde{R}^f}(x_k + s)$, $m_k^{\tilde{R}^c}(x_k + s)$.\;}
	\STATE{Switch Mode according to observed feasibility of trial point.\;}
\ENDWHILE

\STATE{\hspace{-\algorithmicindent}// \textit{STEP 1:} Step calculation}
\STATE{Compute a trial step $x_{\text{trial}} = x_k + s_k$ via~\eqref{eq:M1_trial_point} or~\eqref{eq:M2_trial_point}.\;}
\STATE{Evaluate $R_N^{f}(x_{\text{trial}})$ and $R_N^{c}(x_{\text{trial}})$.\;}
\STATE{Update Gaussian processes and black box evaluations.\;}
\STATE{\hspace{-\algorithmicindent}// \textit{STEP 2:} Check feasibility of trial point}
\IF{$R_N^{c_i}(x_{\text{trial}}) > \tau_i$ for an $i = 1, \ldots, r$}
\STATE{Call Algorithm \ref{alg:noisy_trust_region_management} with $a = \omega$.\;}
\STATE{Evaluate $R_N^{f}(x)$ and $R_N^{c}(x)$ for a randomly sampled $x \in B(x_k, \rho_k)$.\;}
\STATE{Update Gaussian processes, black box evaluations and surrogate models.\;}
\STATE{Switch Mode to observed feasibility of trial point.\;}
\ENDIF
\STATE{\hspace{-\algorithmicindent}// \textit{STEP 3:} Acceptance of trial point and update trust-region}
\STATE{Compute $r_k$ according to Algorithm \ref{alg:acceptance_ratio_mode}.\;}
\IF {$r_k > \eta_0$}
\STATE{Set $x_{k+1} = x_{\text{trial}}$\;}
\STATE{Call Algorithm \ref{alg:noisy_trust_region_management} with $a = \min\{1, \gamma_{inc}\}$\;}
\ELSE
\STATE{Set $x_{k+1} = x_k$, $m_{k+1}^{\tilde{R}^f} = m_k^{\tilde{R}^f}$ and $m_{k+1}^{\tilde{R}^c} = m_k^{\tilde{R}^c}$\;}
\STATE{Call Algorithm \ref{alg:noisy_trust_region_management} with $a = \gamma_{dec}$\;}
\STATE{Evaluate $R_N^{f}(x)$ and $R_N^{c}(x)$ for a randomly sampled $x \in B(x_k, \rho_k)$.\;}
\ENDIF
\STATE{Update Gaussian processes and black box evaluations.\;}
\STATE{Update surrogate models $m_k^{\tilde{R}^f}(x_{k+1} + s)$, $m_k^{\tilde{R}^c}(x_{k+1} + s)$.\;}
\ENDFOR
\end{algorithmic}
\caption{(S)NOWPAC\label{alg:SNOWPAC}}
\end{algorithm}

\section{Numerical examples} \label{sec:numerical_results}
We first discuss a two-dimensional test problem in Section~\ref{sec:low_dim_numerical_example} to build intuition about the optimization process and the effect of the Gaussian process to reduce the noise. 
Thereafter, in Section~\ref{sec:numerical_results_exp} we discuss numerical results for (S)NOWPAC on nonlinear optimization problems from the CUTEst benchmark suite, in particular, benchmark examples from \cite{Hock1981, Schittkowski1987, Schittkowski2008}. We use three different formulations with various combinations of robustness measures from Section~\ref{sec:robustness_formulation} and the data profiles proposed in \cite{More2009} to compare (S)NOWPAC with cBO, COBYLA, NOMAD as well as the stochastic approximation methods SPSA and KWSA. Since COYBLA and NOMAD are not designed for stochastic optimization they will perform better for smaller noise levels. We therefore vary the sample sizes to discuss their performance based on different magnitudes of the noise in the sample approximations of the robust objective function and constraints. 
For the results in this work we employ a stationary square-exponential kernels 
\begin{equation}
\label{eq:stat_gp_kernel}
K^b(x_p, x_q) = {\sigma^b}^2 \prod\limits_{i=1}^n \text{exp}({-\frac12 \left(\frac{\|x_p-x_q\|}{l^b_i}\right)^2})
\end{equation} 
for the construction of the GP surrogates
with standard deviations $\sigma^b$ and length scales $l^b_1, \ldots, l_n$---the hyperparameters of the GP. The hyperparameters $\sigma^b$ and $l^b_1, \ldots, l^b_n$ are found by maximizing the marginal likelihood of the estimator following the automatic relevance determination (ARD) \cite{Rasmussen2006}. The hyperparameters are estimated after a predetermined number of optimization steps. Other kernels may be employed to account for potentially available additional information. 

\subsection{A two dimensional test example}\label{sec:low_dim_numerical_example}
\begin{figure}[!htb]\centering
\begin{center}
\includegraphics[bb=0 0 560 420, width=0.45\textwidth]{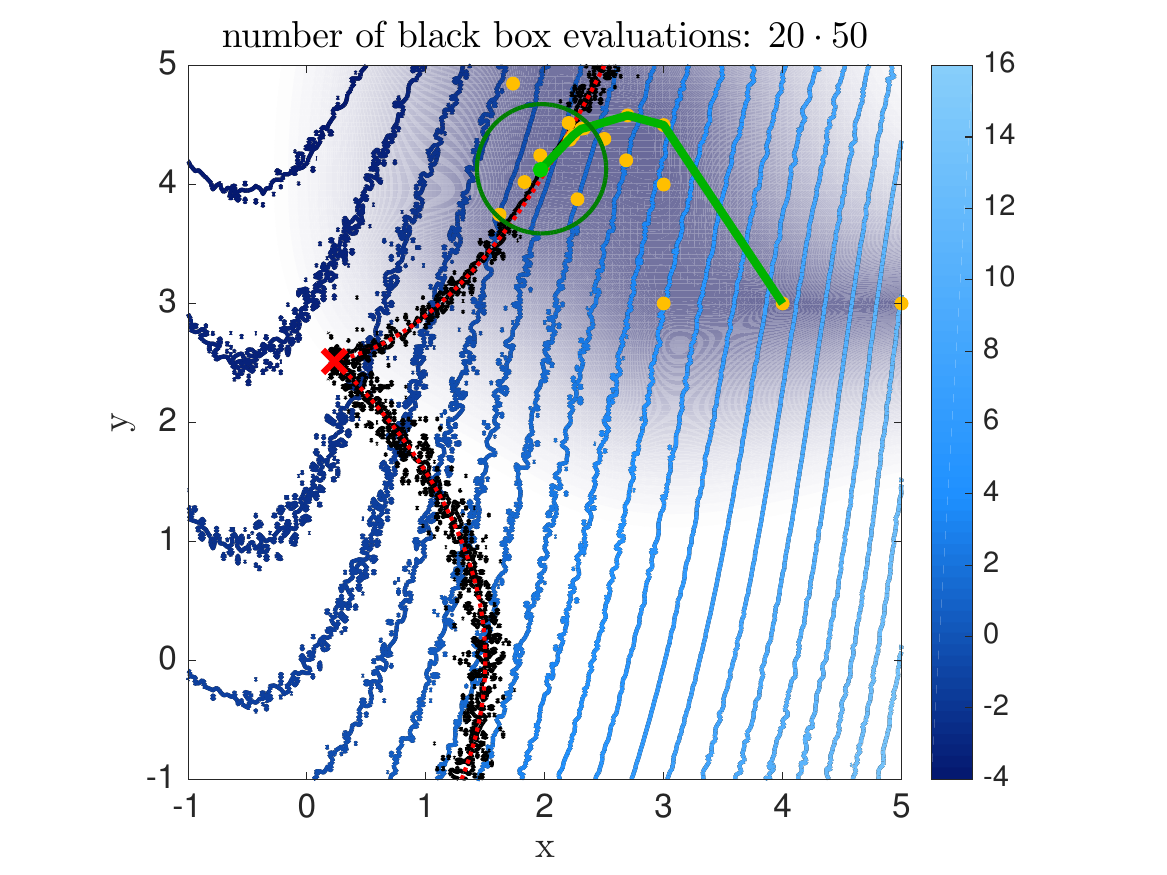}
\includegraphics[bb=0 0 560 420, width=0.45\textwidth]{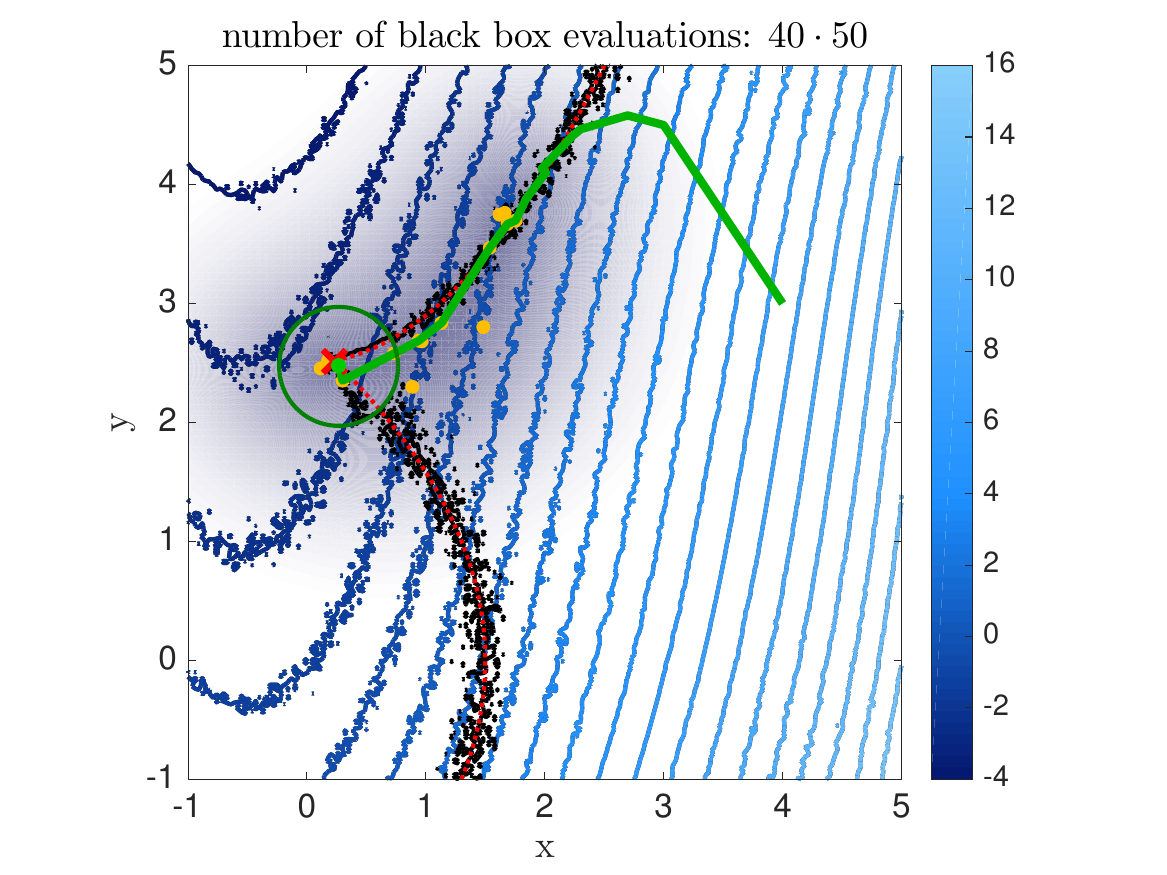}
\includegraphics[bb=0 0 560 420, width=0.45\textwidth]{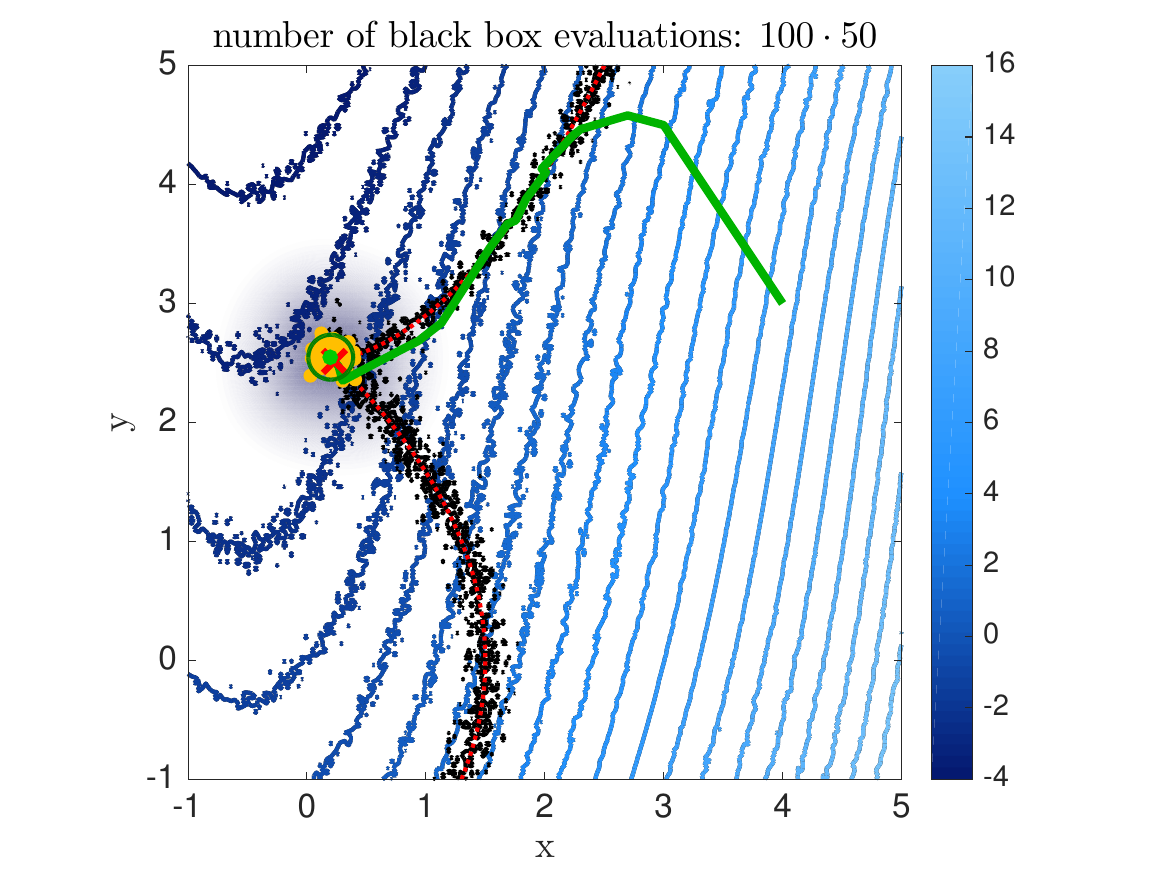}
\includegraphics[bb=0 0 560 420, width=0.45\textwidth]{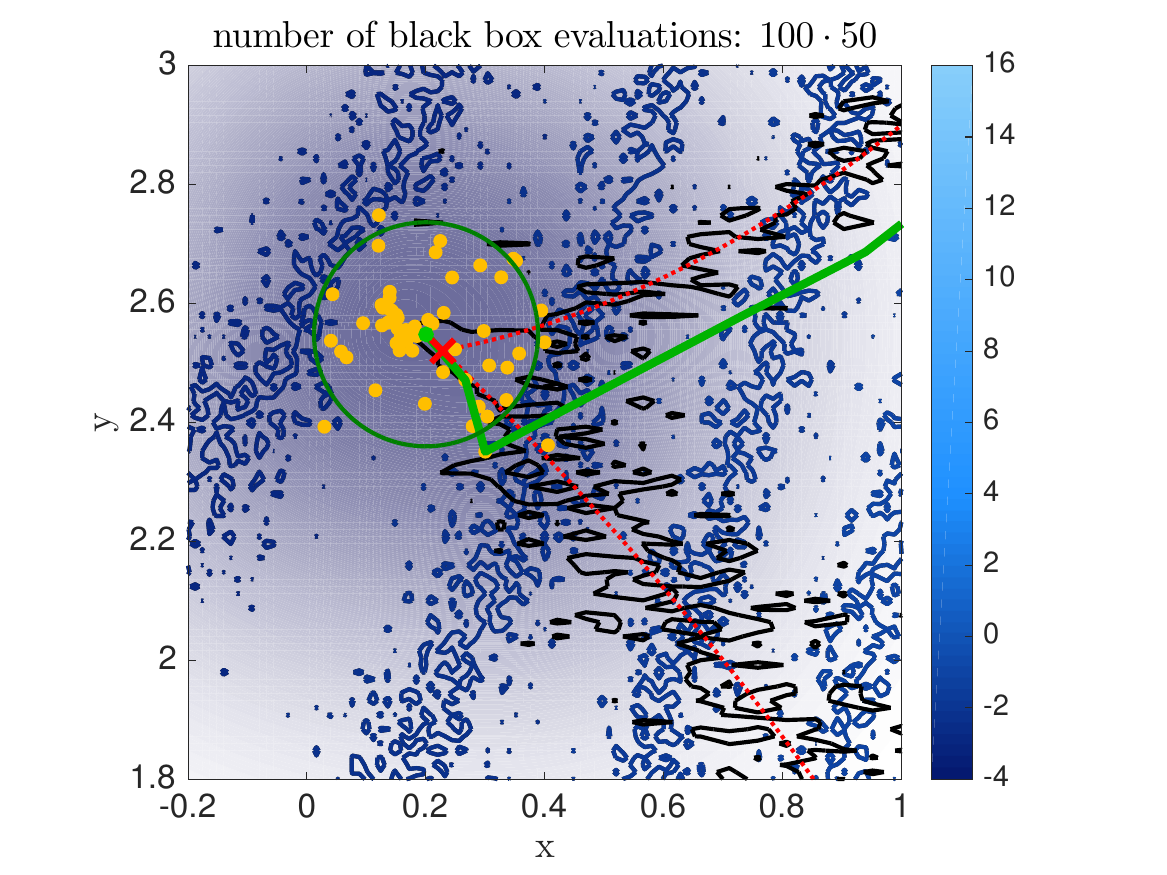}
\end{center}
\caption{Realizations of the contour plots of the noisy objective function and constraints for optimization problem (\ref{eq:first_numerical_example}). The exact constraints are indicated by a dotted red line and the exact optimal point is marked with a red cross. The plots show the best point (green dot) and the optimization path (green line) after $20$, $40$ and $100$ evaluations of the robustness measures; the lower right plot is zoomed in to the neighborhood of the optimal point. The corresponding trust-regions are indicated by green circlec. Within the trust-regions the local smoothing effect of the Gaussian process corrected objective function and constraints can be seen. The gray cloud indicates the size weighing factor $\gamma_s^i$ from (\ref{eq:adjusted_black_box_evaluations}); the darker the area the more weight is given to the Gaussian process mean. The Gaussian regression points are indicated by yellow dots.}\label{fig:resp_surf_numerical_example}
\end{figure}

We consider the optimization problem
\begin{equation}\label{eq:first_numerical_example}
\begin{split}
&\min \mathbb{E}\left[\sin(x-1+\theta_1) + \sin\left(\frac12y-1+\theta_1\right)^2\right] + \frac12\left(x+\frac12\right)^2 - y\\
&\quad\mbox{s.t.} \quad \mathbb{E}\left[-4x^2(1+\theta_2) -10\theta_3\right]  \qquad\;\;\;\leq 25 - 10y
\\
&\quad\;\qquad\mathbb{E}\left[  - 2y^2(1+\theta_4)- 10(\theta_4 + \theta_2) \right] \leq 20x - 15
\end{split}
\end{equation}
with $\theta = (\theta_1, \ldots, \theta_4) \sim \mathcal{U}[-1,1]^4$ and the starting point $x_0 = (4, 3)$. For the approximation of the expected values we use $N = 50$ samples of $\theta$ and we estimate the magnitudes of the noise terms as described in Section~\ref{sec:robustness_sample_approximation}. The noise in the objective function and constraints can be seen in Figure~\ref{fig:resp_surf_numerical_example}. The feasible domain is to the right of the exact constraints which are indicated by dotted red lines. We see that the noise is the largest in the region around the optimal solution (red cross). 

To show the effect of the noise reduction we introduced in Section~\ref{sec:gaussian_process_supported_noise_correction}, we plot the objective function and the constraints corrected by the respective Gaussian process surrogates around the current design point at $20$ (upper left), $40$ (upper right) and $100$ (lower plots) evaluations of the robustness measures. We see that the noise is reduced which enables (S)NOWPAC to efficiently approximate the optimal solution. 

Note that the GP-corrected approximate feasible domains within the trust-region show significantly less noise than outside of the trust-region. Moreover, we see that the optimizer eventually gathers more and more black box evaluations, yielding an increasingly better noise reduction. 
Looking at the noisy constraint contours at $40$ evaluations, we see that the quantification of feasibility based on the Gaussian process supported black box evaluations is not always reliable. This underlines the necessity of the feasibility restoration mode we introduced in Section~\ref{sec:relaxed_feasibility_requirement}, which allows the optimizer to recover feasibility from points that appear infeasible. 

\subsection{Optimization performance on benchmark test set}\label{sec:numerical_results_exp}

Its utilization of Gaussian process surrogate models relates (S)NOWPAC to the successful class of Bayesian optimization techniques \cite{Mockus1974, Mockus1989}, and its extensions for nonlinear optimization using either an augmented Lagrangian approach \cite{Gramacy2015} or expected constrained improvement in the constrained Bayesian optimization (cBO) \cite{Gardner2014}. As opposed to Bayesian optimization, (S)NOWPAC introduces Gaussian process surrogates to smooth local trust-region steps instead of aiming at global optimization. We will demonstrate that the combination of fast local optimization with a second layer of {\it smoothing} Gaussian process models makes (S)NOWPAC an efficient and accurate optimization technique. Additionally, we compare the performance of (S)NOWPAC to the optimization codes COBYLA and NOMAD as well as to the stochastic approximation methods SPSA and KWSA.

We test the performances of all optimizers on the Schittkowski optimization benchmark set \cite{Hock1981, Schittkowski2008}, which is part of the CUTEst benchmark suit for nonlinear constraint optimization.
The dimensions of the feasible domains within our test set range from $2$ to $16$ with a number of constraints ranging from $1$ to $10$. 
Since the problems are deterministic, we add noise to the objective functions, $f(x) + \theta_1$ and constraints, $c(x) + \theta_2$ with $(\theta_1, \theta_2) \sim \mathcal{U}[-1, 1]^{1+r}$ and solve the following three classes of robust optimization problems: \todo{FM: Only linear noise, could we run into issues here from review?}
\begin{enumerate}
\item Minimization of the average objective function subject to the constraints being satisfied in expectation:
\begin{equation}\label{eq:robust_expected_value_Schittkowski_problems}
\begin{split}
&\;\;\;\min \mathcal{R}_0^f (x)\\
&\mbox{s.t.} \quad \mathcal{R}_0^c (x) \leq 0.
\end{split}
\end{equation}
\item Minimization of the average objective function subject to the constraints being satisfied in $95\%$ of all cases:
\begin{equation}\label{eq:robust_quantile_value_Schittkowski_problems}
\begin{split}
&\;\;\;\min \mathcal{R}_0^f (x) \\
&\mbox{s.t.} \quad \mathcal{R}_4^{c,0.95}(x) \leq 0.
\end{split}
\end{equation}
\item Minimization of the $95\%$-CVaR of the objective function subject to the constraints being satisfied on average:
\begin{equation}\label{eq:robust_cvar_value_Schittkowski_problems}
\begin{split}
&\;\;\;\min \mathcal{R}_5^{f,0.95}(x) \\
&\mbox{s.t.} \quad \mathcal{R}_0^c(x) \leq 0,
\end{split}
\end{equation}
\end{enumerate}

\begin{figure}[!htb]
\begin{center}
\includegraphics[width=0.45\textwidth]{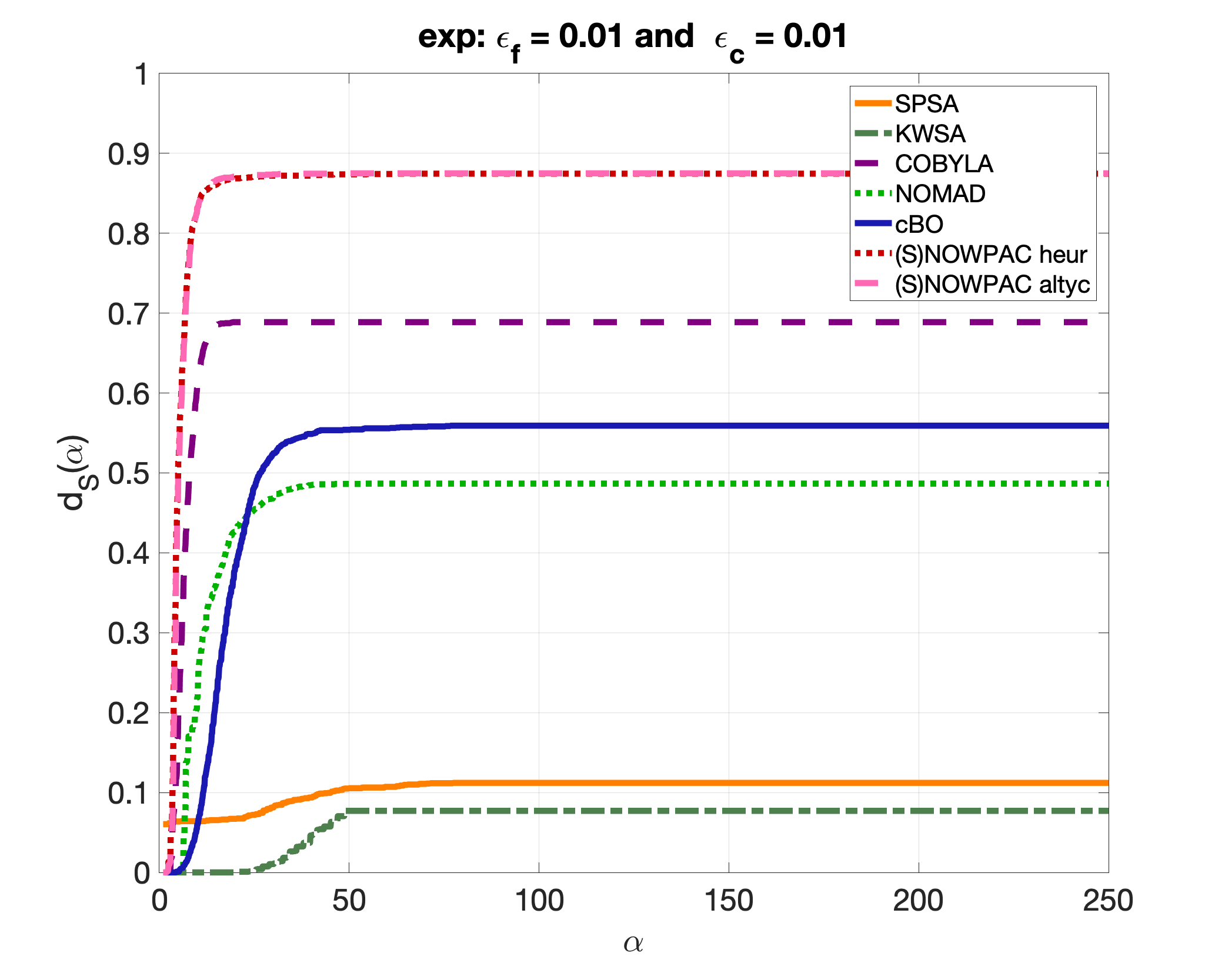}
\includegraphics[width=0.45\textwidth]{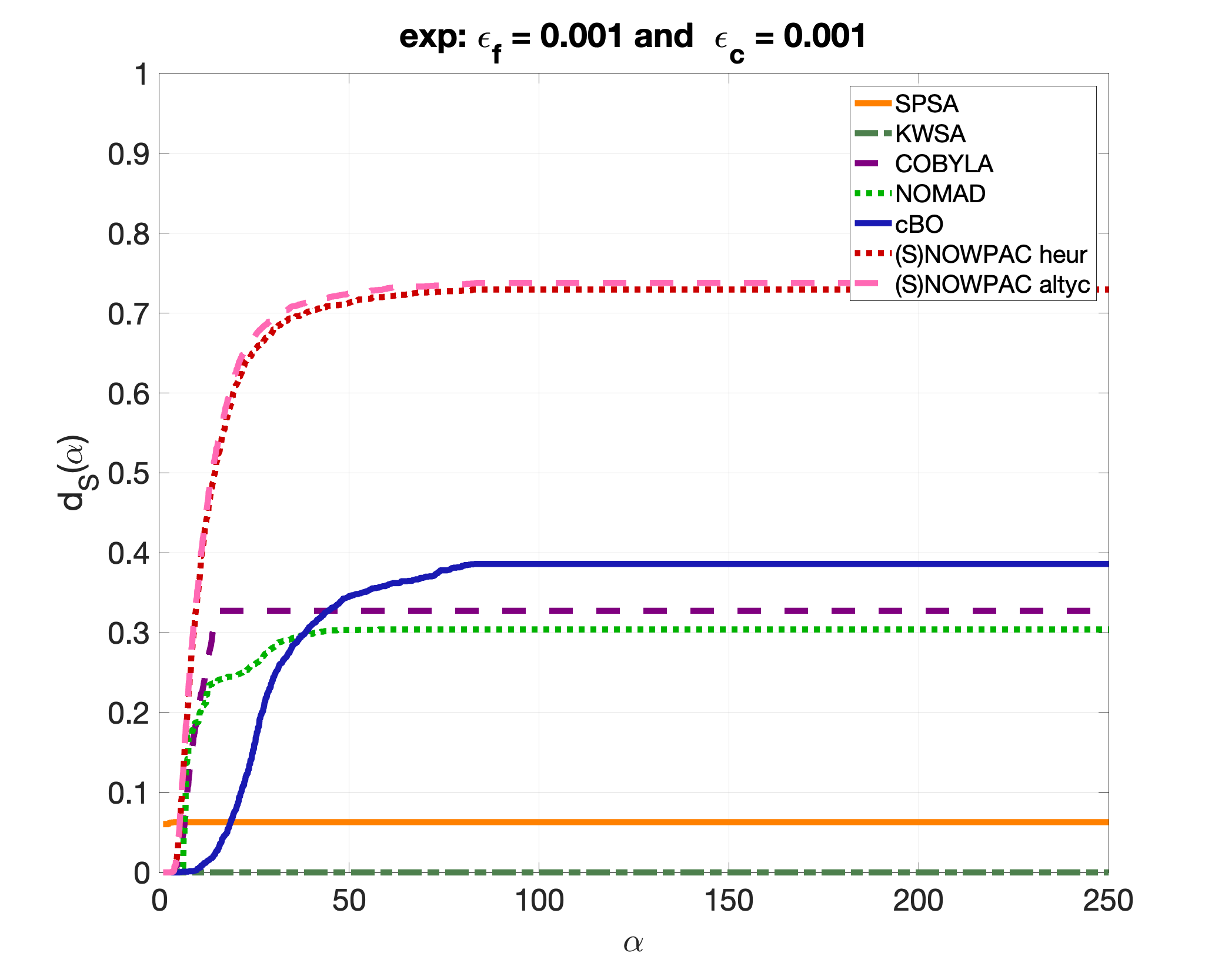}\\[2mm]
\includegraphics[width=0.45\textwidth]{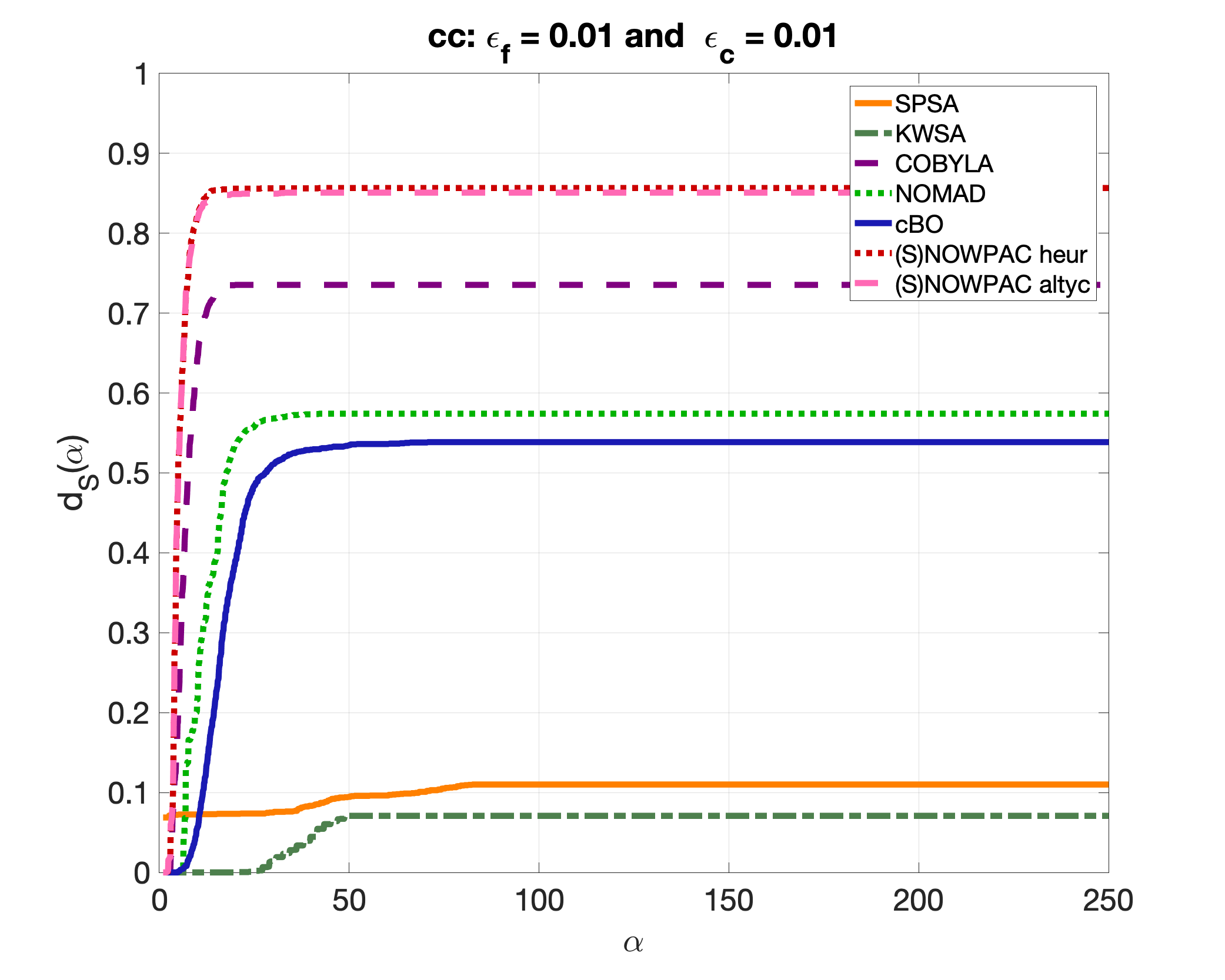}
\includegraphics[width=0.45\textwidth]{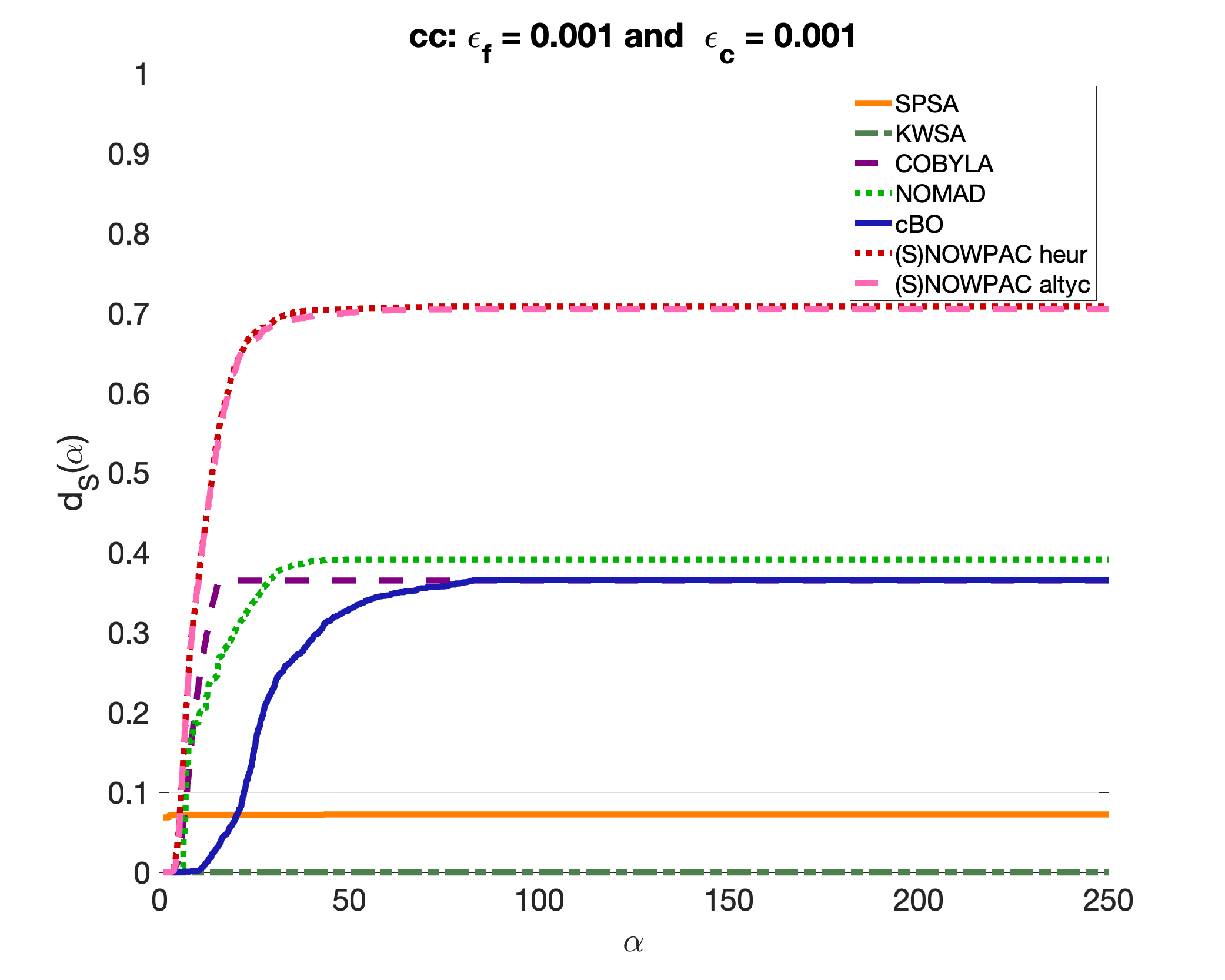}
\includegraphics[width=0.45\textwidth]{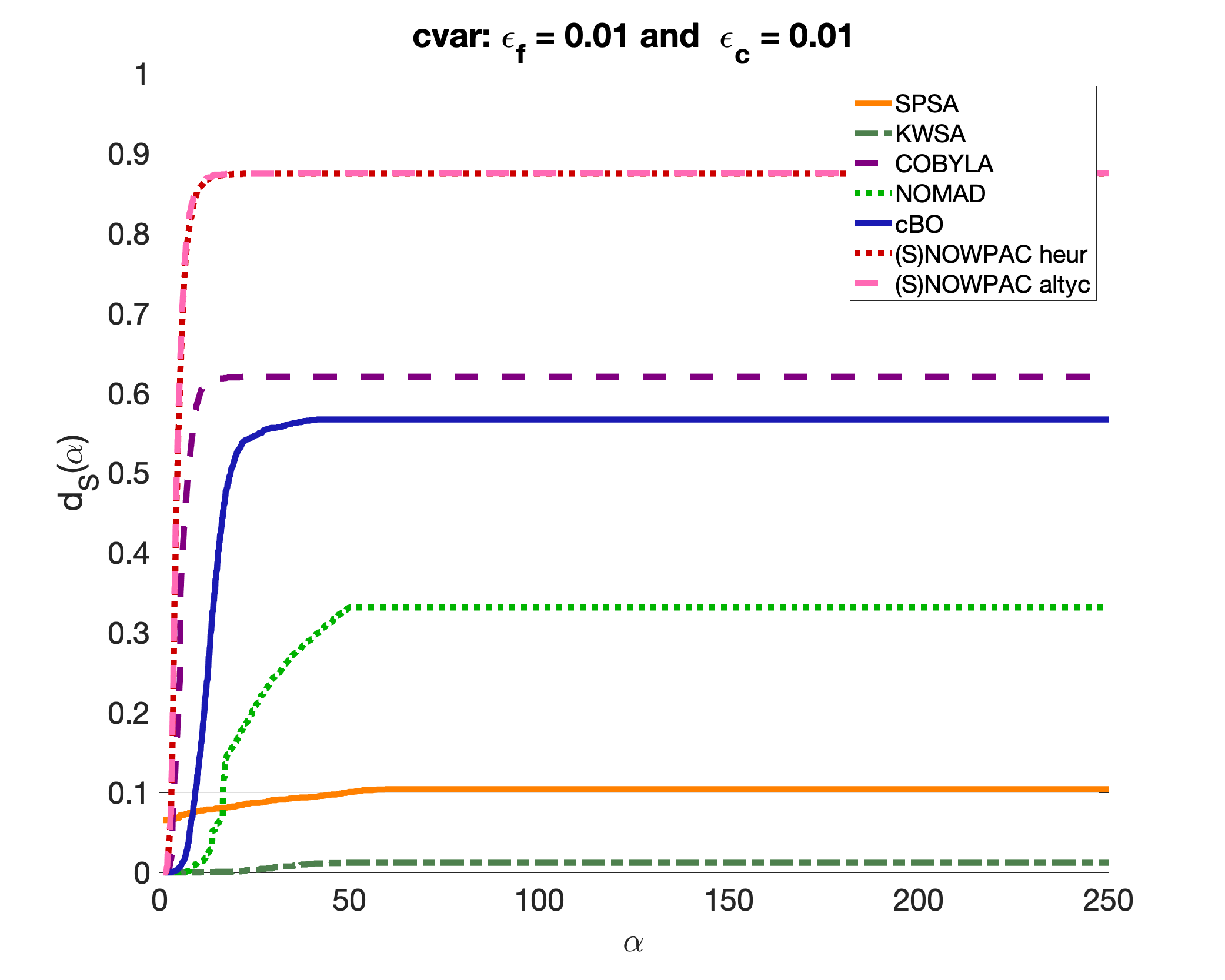}
\includegraphics[width=0.45\textwidth]{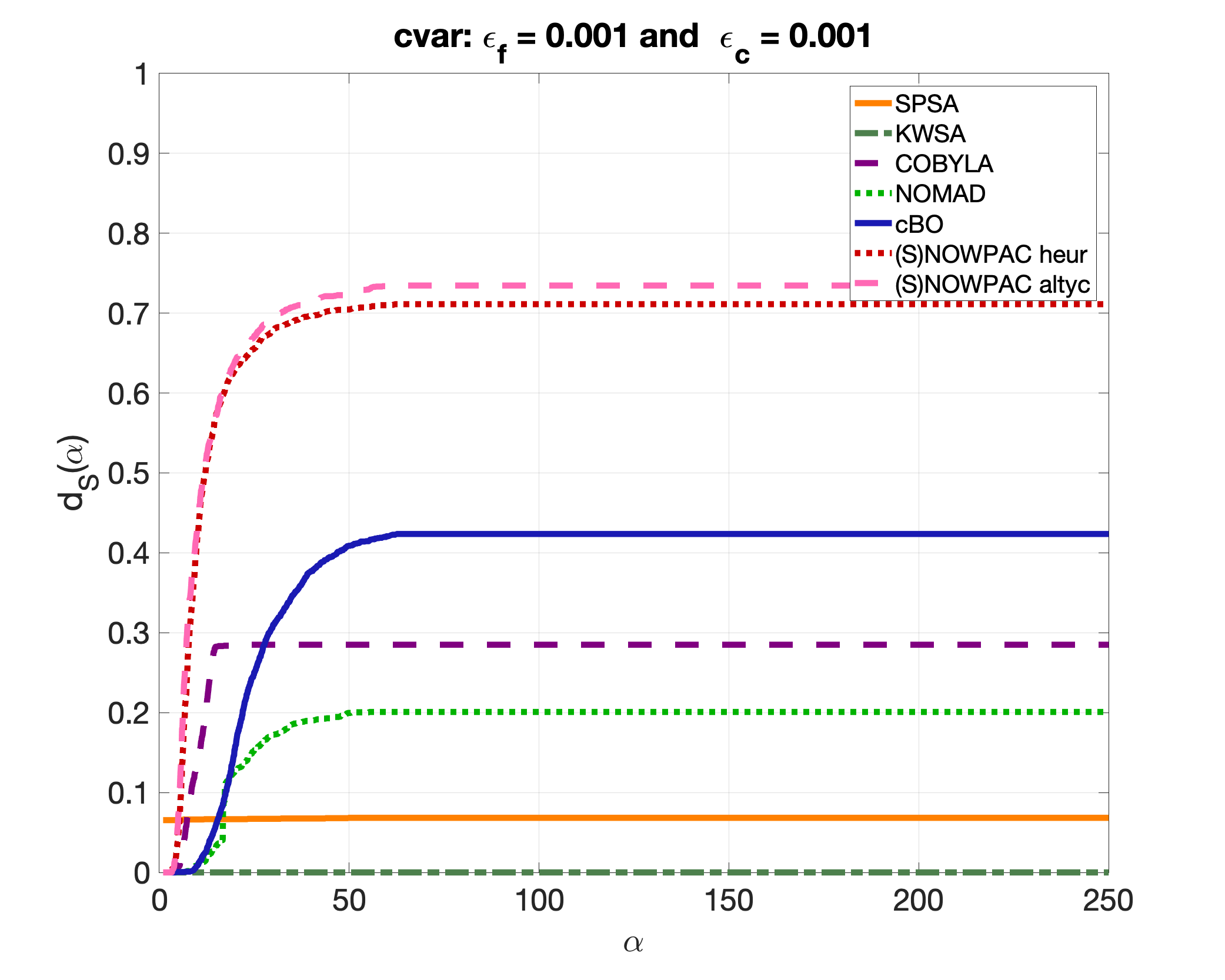}
\end{center}
\caption{Data profiles for (S)NOWPAC altyc (pink), (S)NOWPAC heur (red), cBO (blue), COBYLA (purple), NOMAD (green), SPSA (orange) and KWSA (dark green) of $2400$ runs of the benchmark problems. The results for (\ref{eq:robust_expected_value_Schittkowski_problems}), (\ref{eq:robust_quantile_value_Schittkowski_problems}) and (\ref{eq:robust_cvar_value_Schittkowski_problems}) are plotted in the first, second and third row respectively. The profiles shown are based on the exact values for the objective function and constraints evaluated at the intermediate points computed by the respective optimizers. The data profiles are shown for varying thresholds $\epsilon_f \in \{10^{-2}, 10^{-3}\}$ and $\epsilon_c \in \{10^{-2}, 10^{-3}\}$ on the objective values and the constraint violation respectively.}\label{fig:data_profiles}
\end{figure}

For the performance comparison we use a total number of $3 \cdot 8 \cdot 100 = 2400$ optimization runs ($3$ different number of Monte Carlo sampling sizes $N \in \{200, 1000, 2000\}$, $8$ benchmark problems with $100$ repeated optimization runs) and denote the benchmark set by $\mathcal{P}$. To obtain the data profiles we determine the minimal number $t_{p,S}$ of optimization steps a solver $S$ requires to solve problem $p \in \mathcal{P}$ under the accuracy requirement
\[
\frac{\left| \mathcal{R}^f(x_k) -  \mathcal{R}^f(x^\ast) \right|}{
\max\{1, \left|\mathcal{R}^f(x^\ast)\right| \}} \leq \varepsilon_f \quad \mbox{and}\quad
\max\limits_{i=1}^r\left\{ \left[\mathcal{R}^{c_i}(x_k)\right]^+  \right\} \leq \varepsilon_c.
\]
Hereby we limit the maximal number of optimization steps to $250$ and set $t_{p,S} = \infty$ if the accuracy requirement is not met after $250\cdot N$ black box evaluations. To decide whether the accuracy requirement is met, we use the {\it exact} objective and constraint values of the robustness measures which we obtained in a post-processing step. Specifically, we use the data profile
\[
d_{S}(\alpha) = \frac{1}{2400}\left|\left\{
p \in \mathcal{P} \; : \; \frac{t_{p,S}}{n_p + 1}\leq \alpha
\right\}\right|,
\]
where $n_p$ denotes the number of design parameters in problem $p$. We remark that, although this allows us to eliminate the influence of the noise on the performance evaluation, it is information that is not available in general. For this reason, we also include a more detailed analysis of individual optimization results below. Figure~\ref{fig:data_profiles} shows the data profiles for different error thresholds $\epsilon_f \in \{10^{-2}, 10^{-3}\}$ and $\epsilon_c \in \{10^{-2}, 10^{-3}\}$ and for (S)NOWPAC altyc (pink), (S)NOWPAC heur (red), cBO (blue), COBYLA (purple), NOMAD (green), SPSA (orange) and KWSA (dark green) respectively. Here, "altyc" stands for the analytic approach described in Section \ref{sec:GPsmoothing} while "heur" uses the heuristic approach introduced in \eqref{eq:heuristicsmoothing}. 

We see that both (S)NOWPAC approaches 
solve the most test problems within the given budget of black box evaluations. Looking at the performance for small values of $\alpha$ we also see that (S)NOWPAC exhibits a comparable or superior performance, indicating fast initial descent which is highly desirable in particular if the evaluations of the robustness measures is computationally expensive.

The performance of cBO suffers in the higher dimensional benchmark problems. Here in particular the global optimization strategy of cBO naturally requires more function evaluations. Furthermore, we used the stationary kernel \eqref{eq:stat_gp_kernel}, which may not properly reflect the properties of the objective functions and constraints. A problem dependent choice of kernel function might help to reduce this problem, however, this information is often hard to obtain in black box optimization. With the localized usage of Gaussian process approximations, as we introduced in Section~\ref{sec:gaussian_process_supported_noise_correction}, (S)NOWPAC reduces the problem of violated stationarity assumptions on the objective function and constraints.

As expected, COBYLA and NOMAD perform well for larger thresholds that are of the same magnitudes as the noise term in some test problems. The noise reduction in (S)NOWPAC using the Gaussian process support helps to approximate the optimal solution more accurately, resulting in better performance results. The Stochastic Approximation approaches SPSA and KWSA, despite a careful choice of hyper-parameters, do not perform well on the benchmark problems. This can be explained by the limited number of overall optimization iterations not being sufficient to achieve a good approximation of the optimal solution using inaccurate gradients.

\begin{figure}[!htb]
\begin{center}
\includegraphics[width=0.3\textwidth]{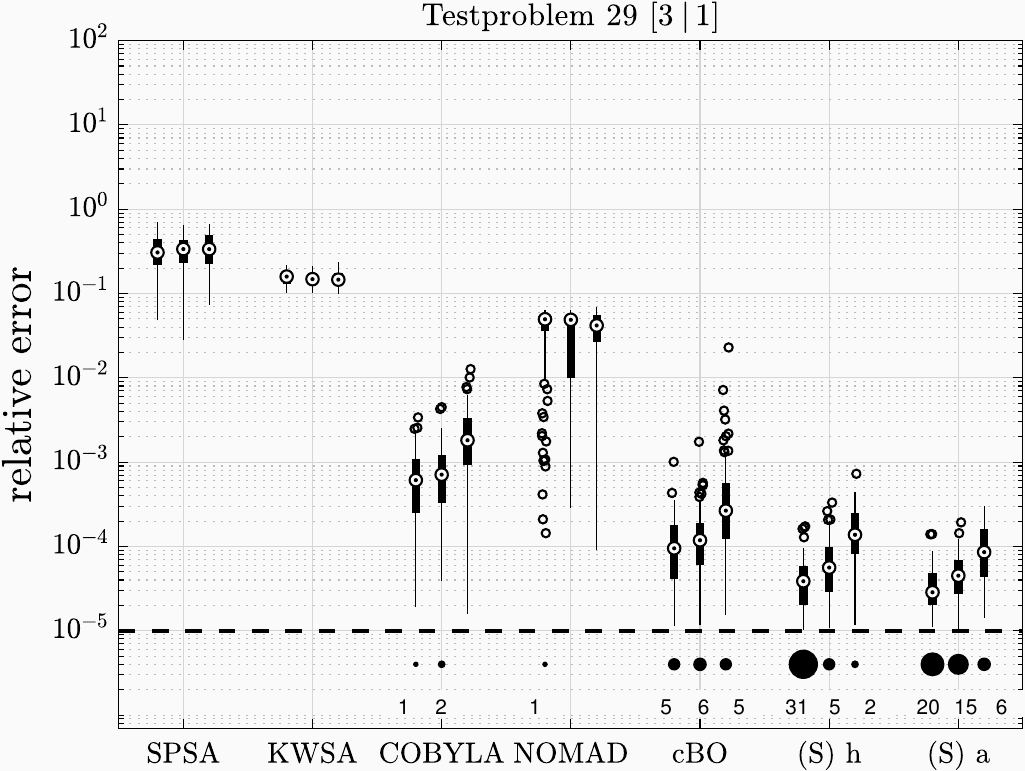}\hspace{2mm}
\includegraphics[width=0.3\textwidth]{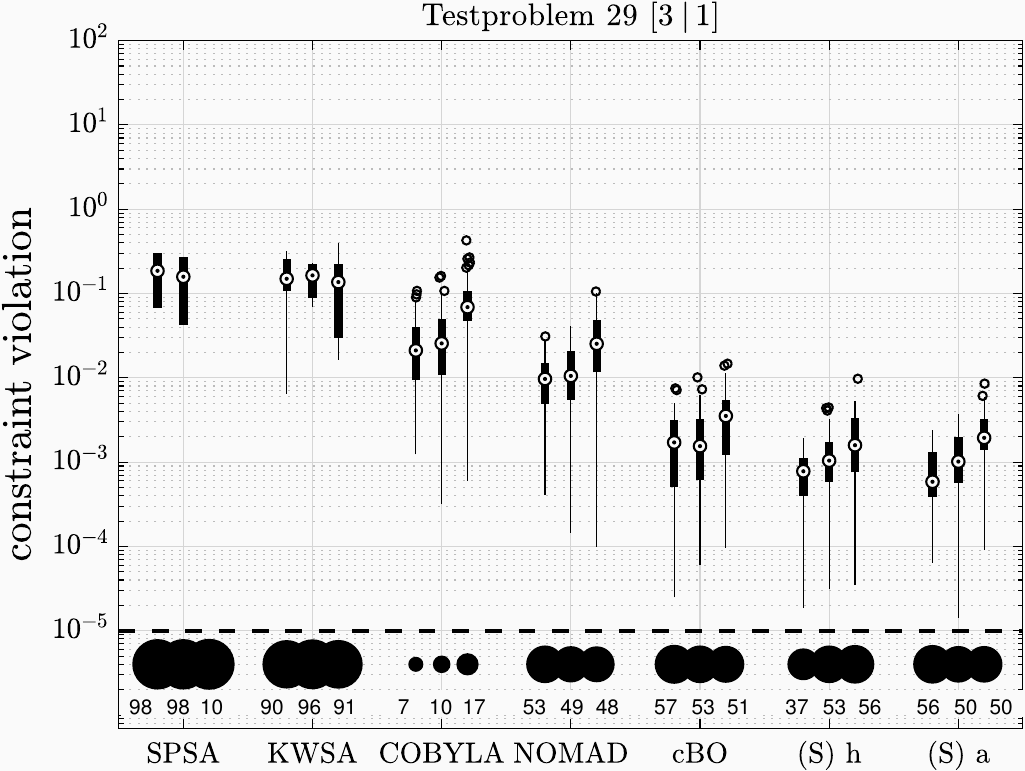}\hspace{2mm}
\includegraphics[width=0.3\textwidth]{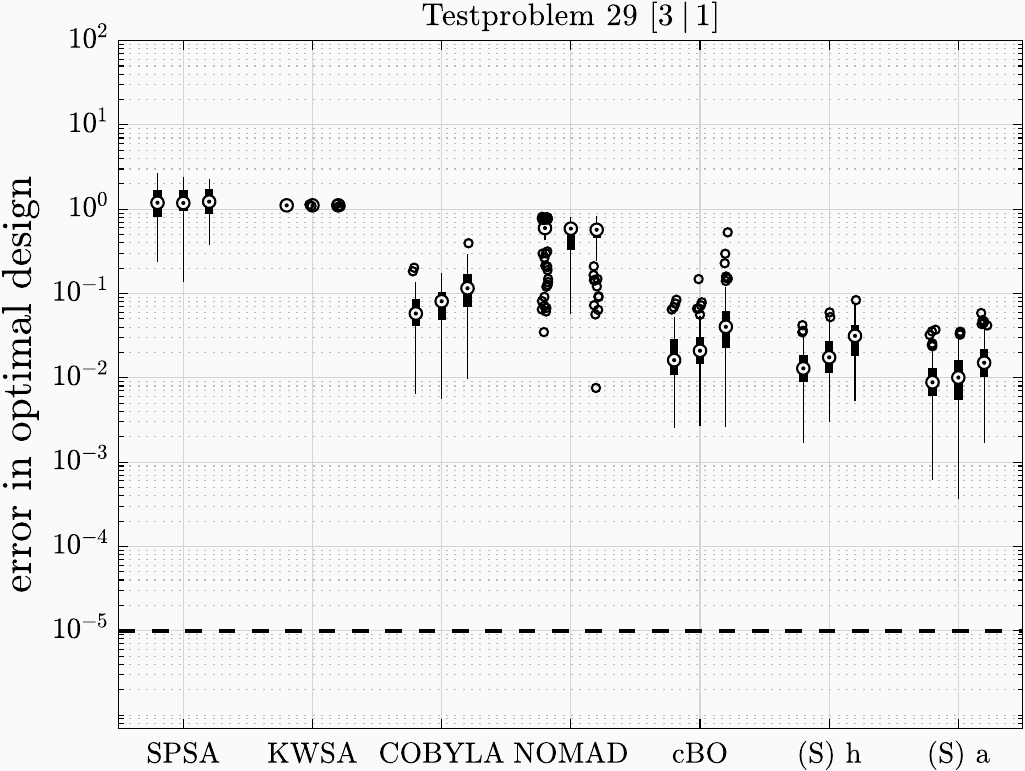}\\
\includegraphics[width=0.3\textwidth]{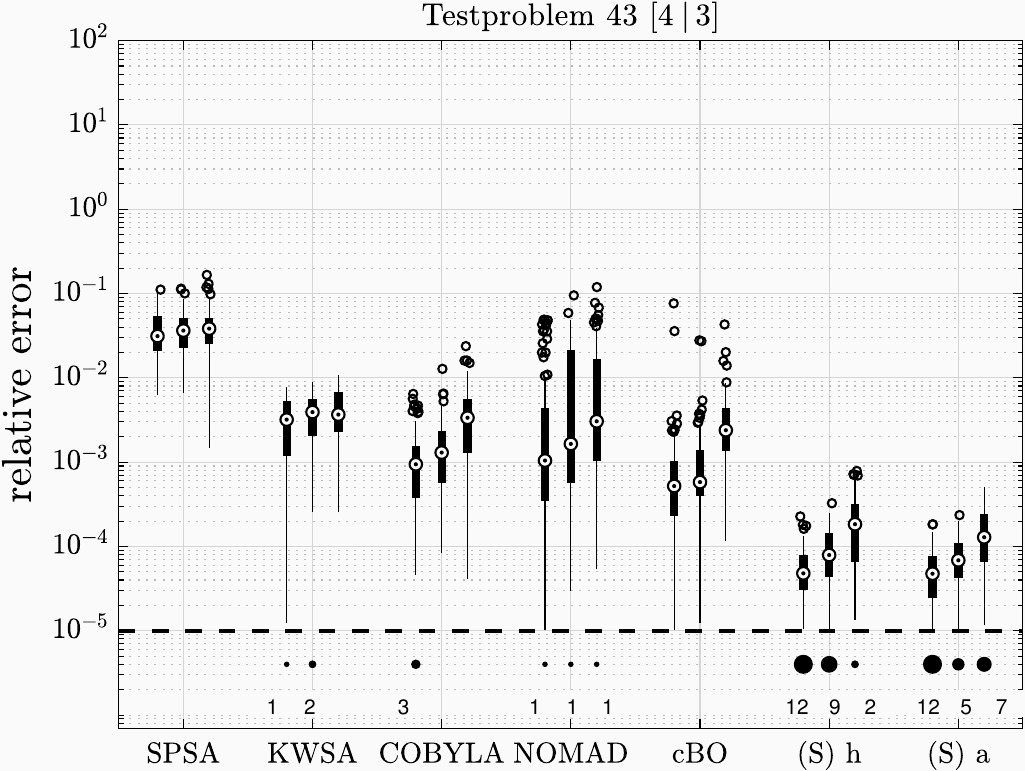}\hspace{2mm}
\includegraphics[width=0.3\textwidth]{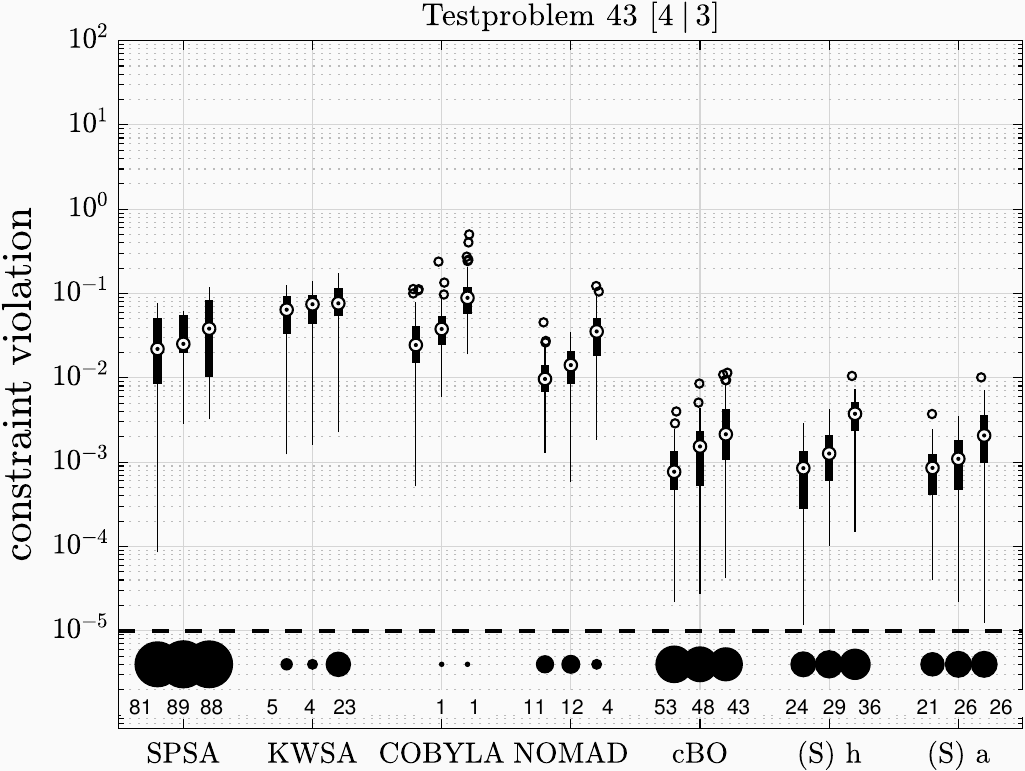}\hspace{2mm}
\includegraphics[width=0.3\textwidth]{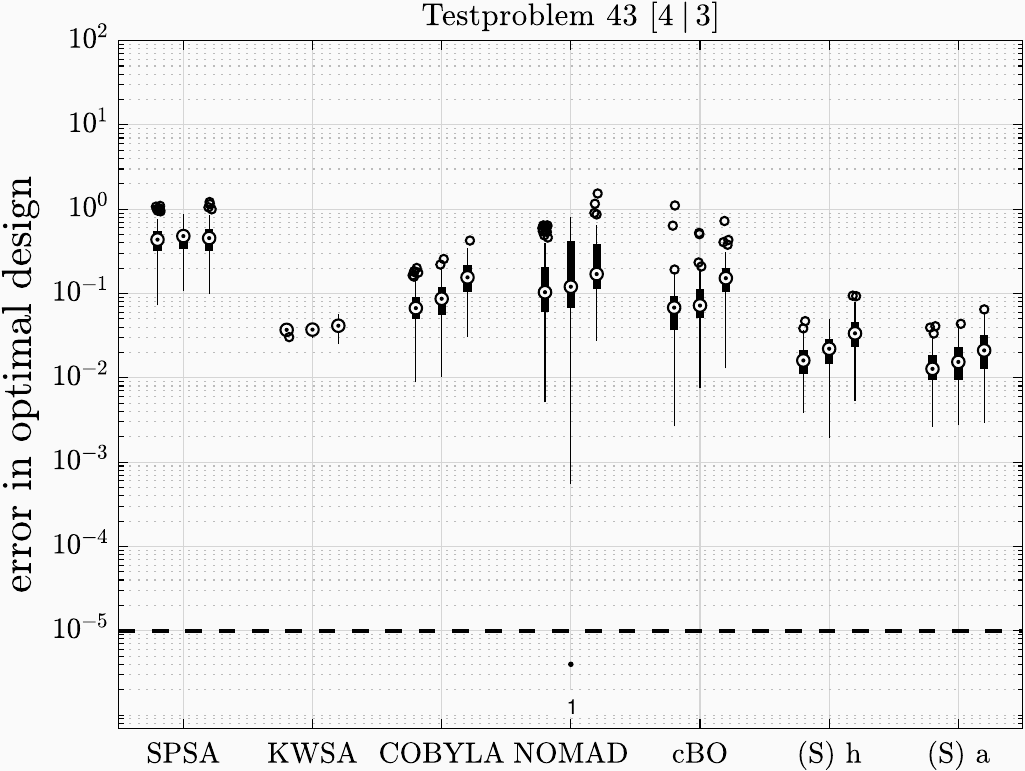}\\
\includegraphics[width=0.3\textwidth]{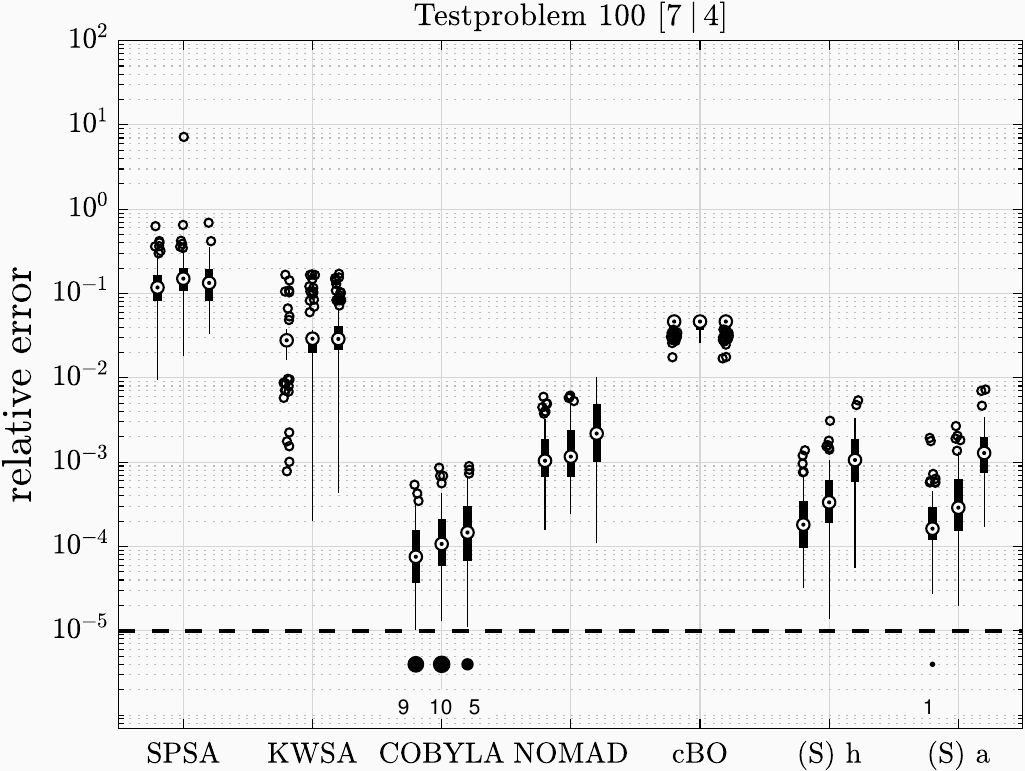}\hspace{2mm}
\includegraphics[width=0.3\textwidth]{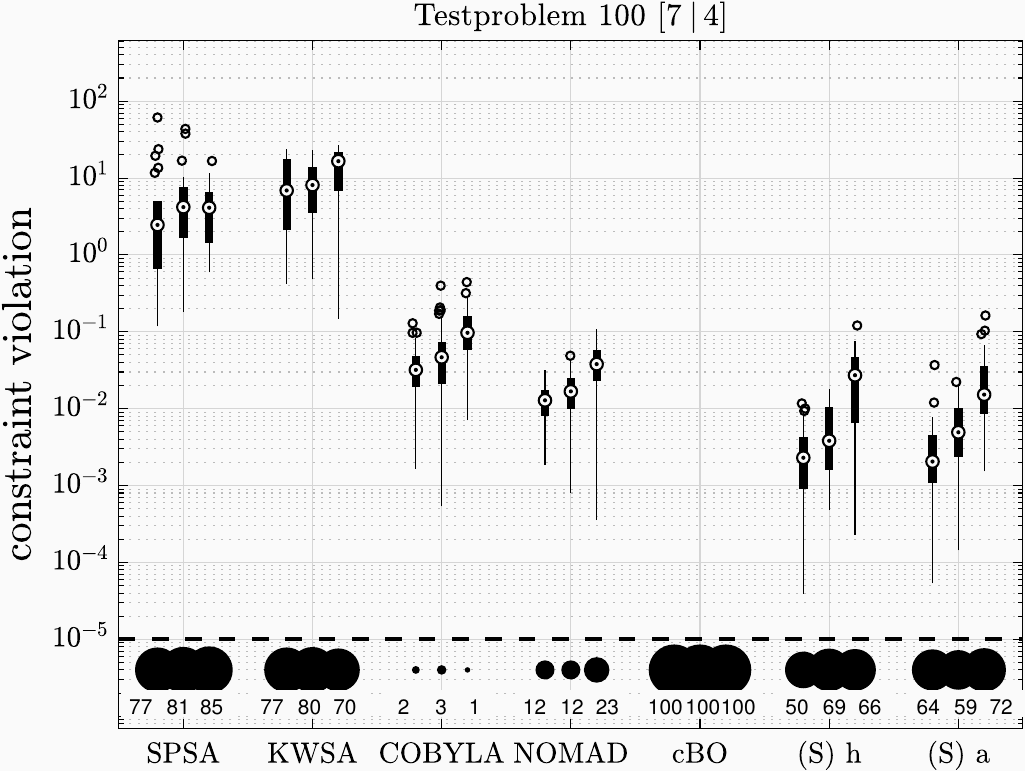}\hspace{2mm}
\includegraphics[width=0.3\textwidth]{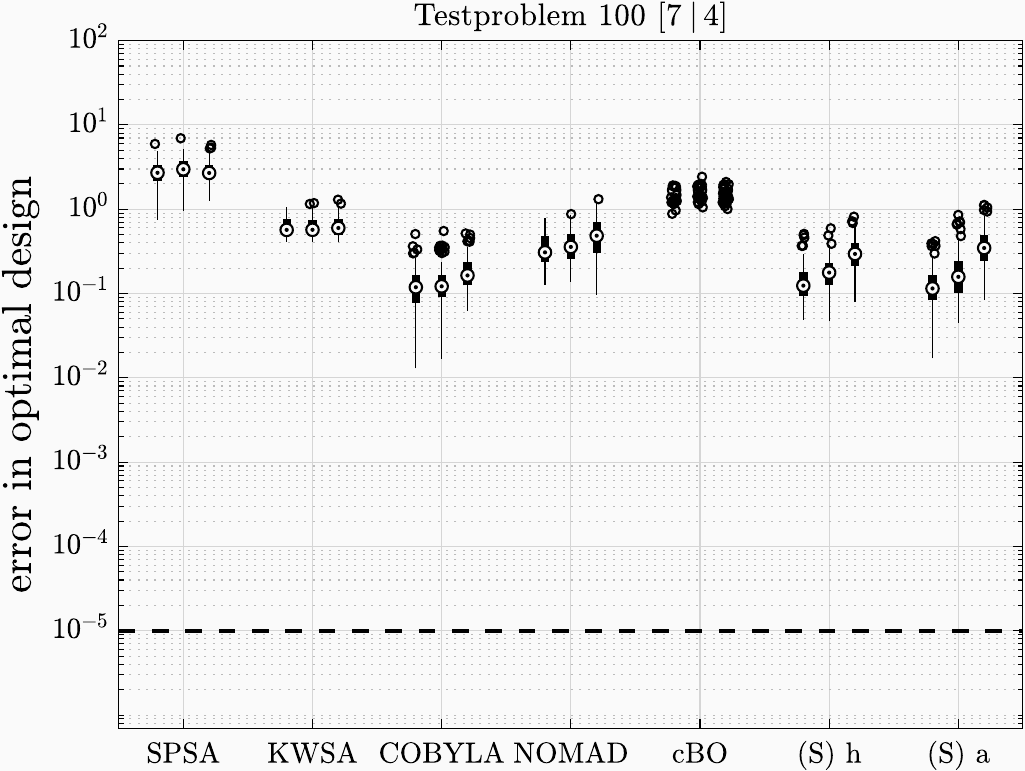}\\
\includegraphics[width=0.3\textwidth]{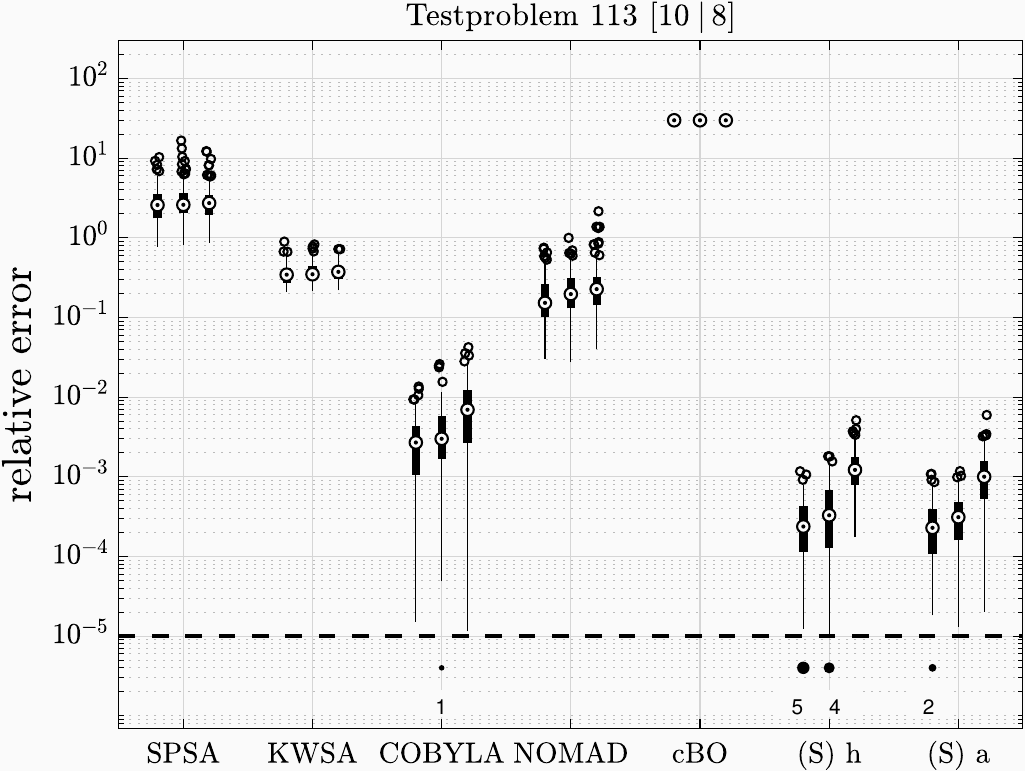}\hspace{2mm}
\includegraphics[width=0.3\textwidth]{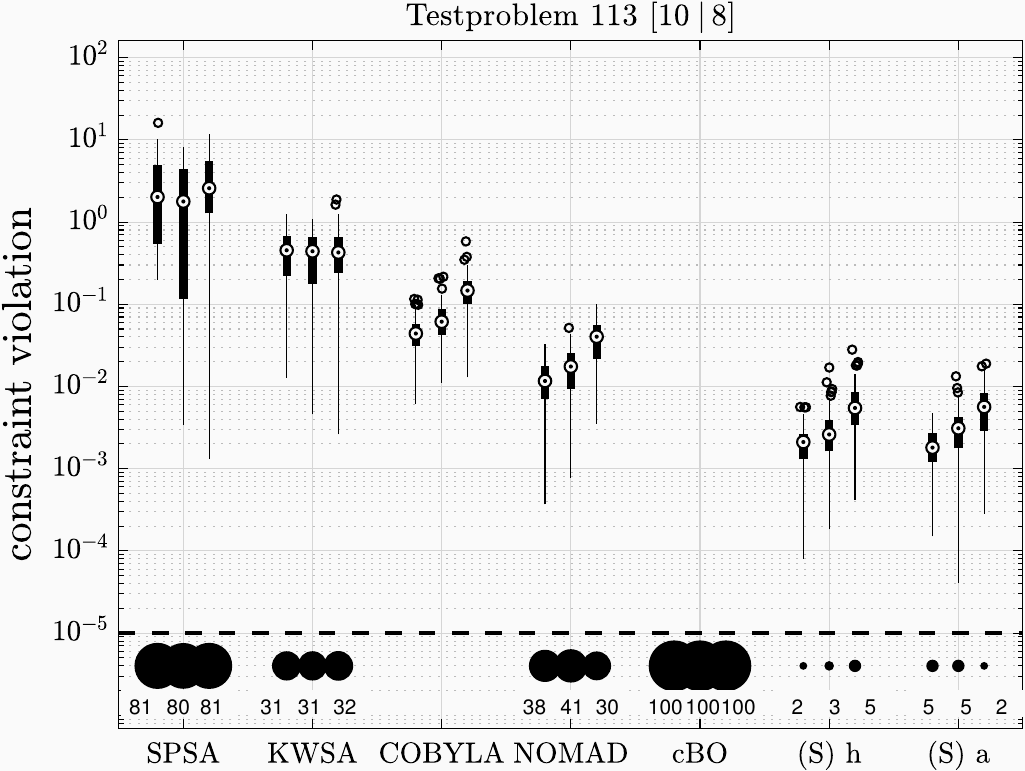}\hspace{2mm}
\includegraphics[width=0.3\textwidth]{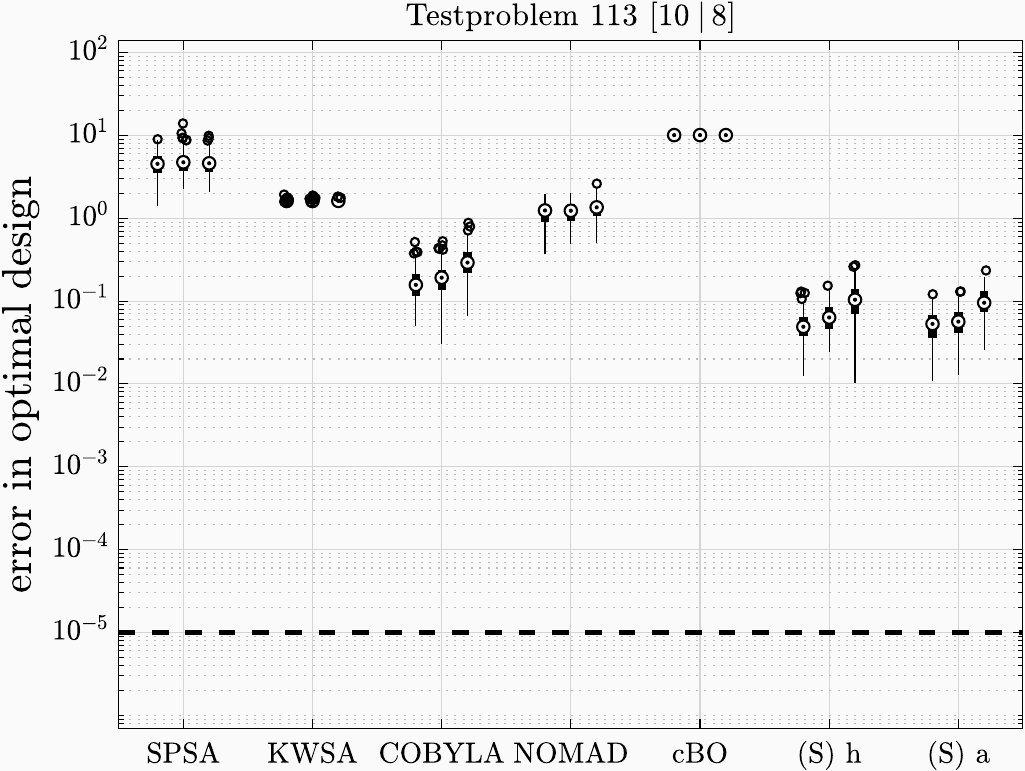}
\caption{Box plots of the errors in the approximated optimal objective values (left plots), the constraint violations (middle plots) and the $l_2$ distance to the exact optimal solution (right plots) of $100$ repeated optimization runs for the Schittkowski test problems number $29$, $43$, $100$, and $113$ for (\ref{eq:robust_expected_value_Schittkowski_problems}). The plots show results of the exact objective function and constraints evaluated at the approximated optimal design computed by (S)NOWPAC analytic ((S) a), (S)NOWPAC heuristic ((S) h), cBO, COBYLA, NOMAD, SPSA and KWSA. All errors or constraint violations below $10^{-5}$ are stated separately below the $10^{-5}$ threshold and the box plots only contain data above this threshold.}\label{fig:schittkowski_testset1a}
\end{center}
\end{figure}

\begin{figure}[!htb]
\begin{center}
\includegraphics[width=0.3\textwidth]{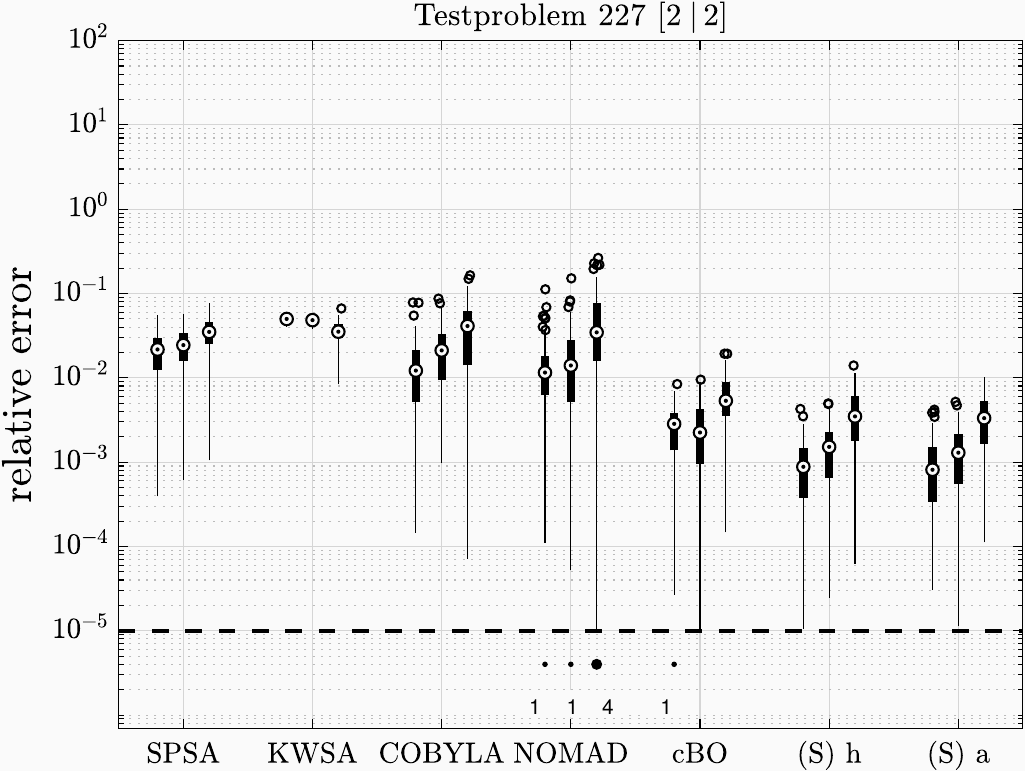}\hspace{2mm}
\includegraphics[width=0.3\textwidth]{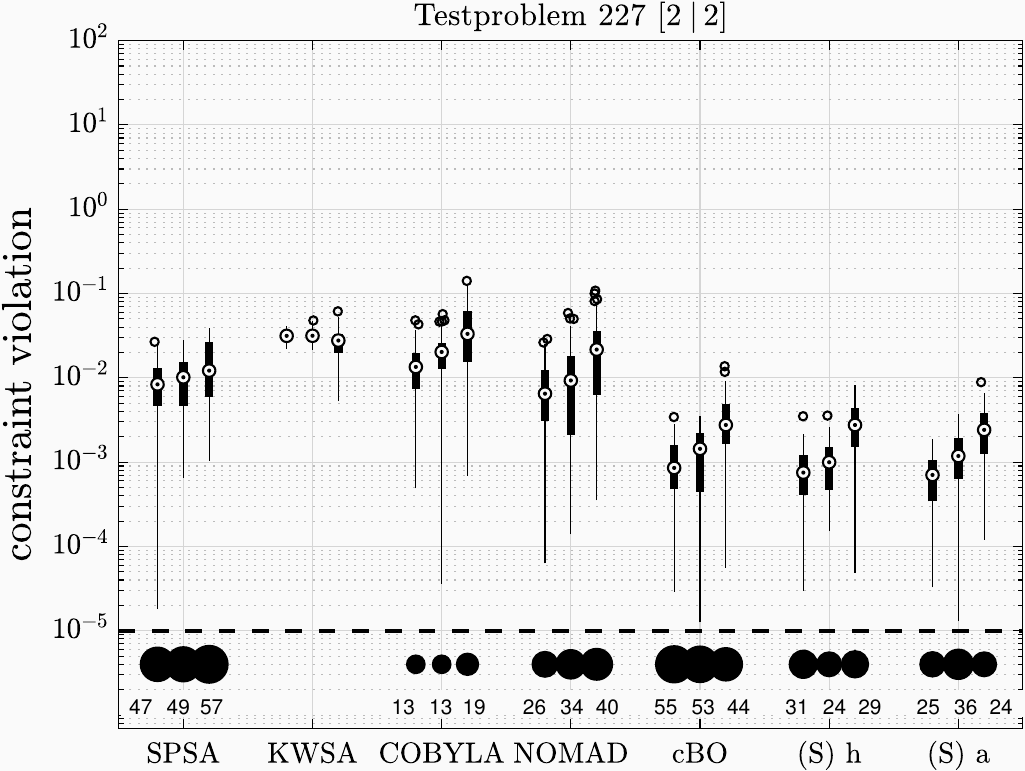}\hspace{2mm}
\includegraphics[width=0.3\textwidth]{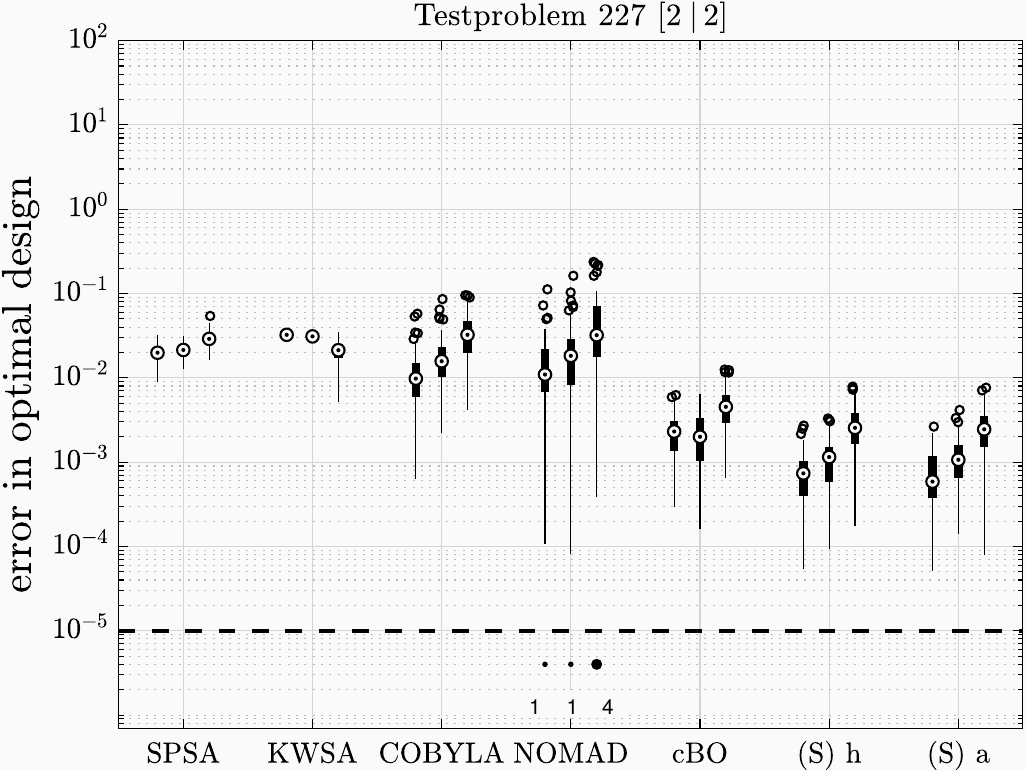}\\[3mm]
\includegraphics[width=0.3\textwidth]{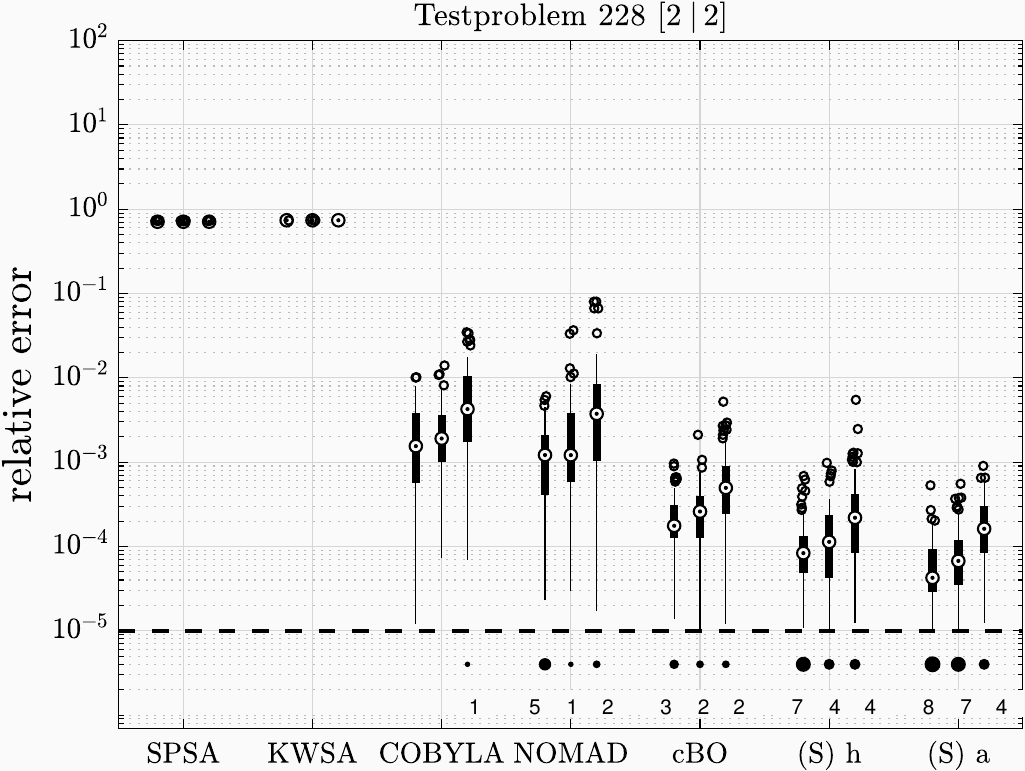}\hspace{2mm}
\includegraphics[width=0.3\textwidth]{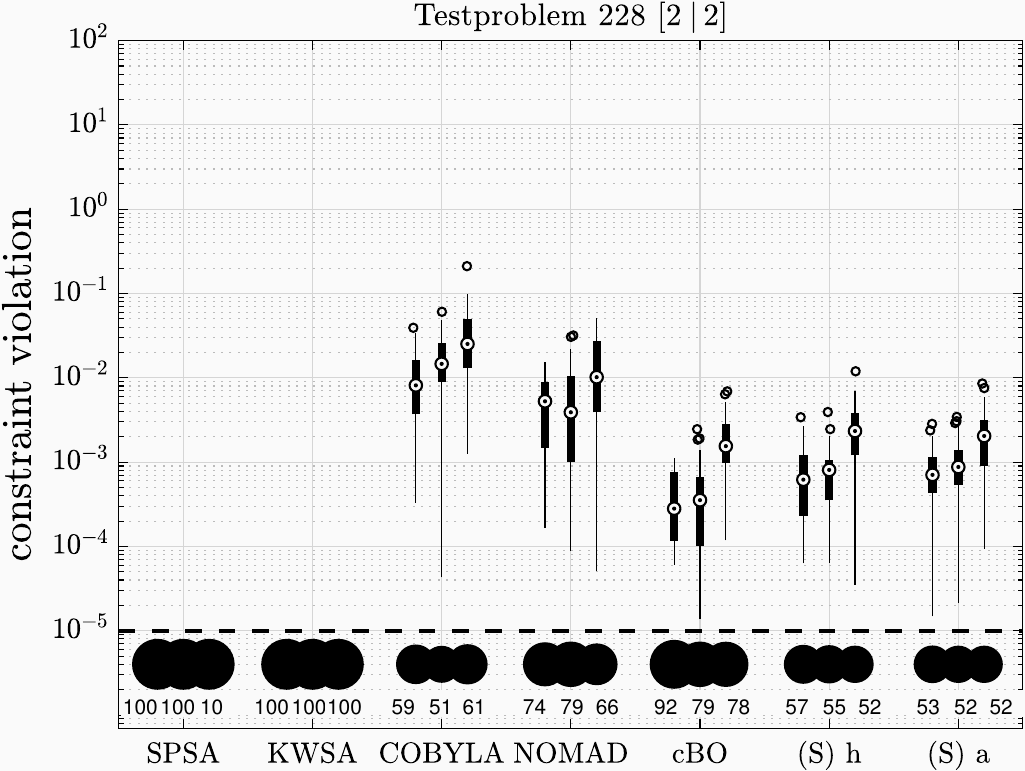}\hspace{2mm}
\includegraphics[width=0.3\textwidth]{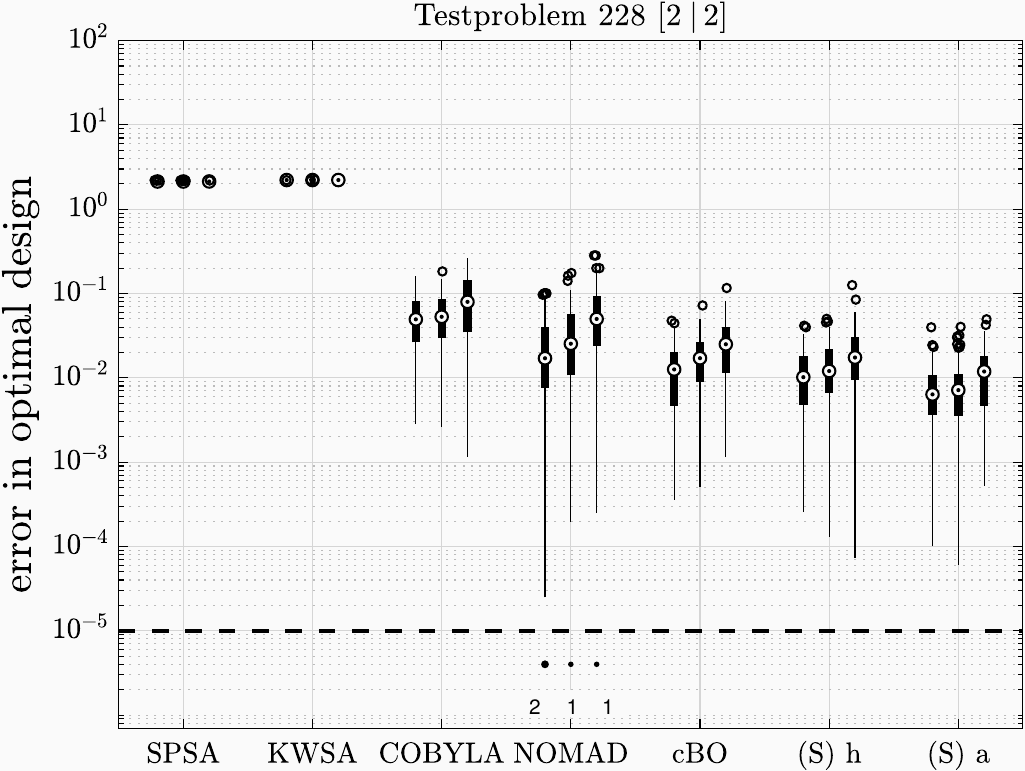}\\[3mm]
\includegraphics[width=0.3\textwidth]{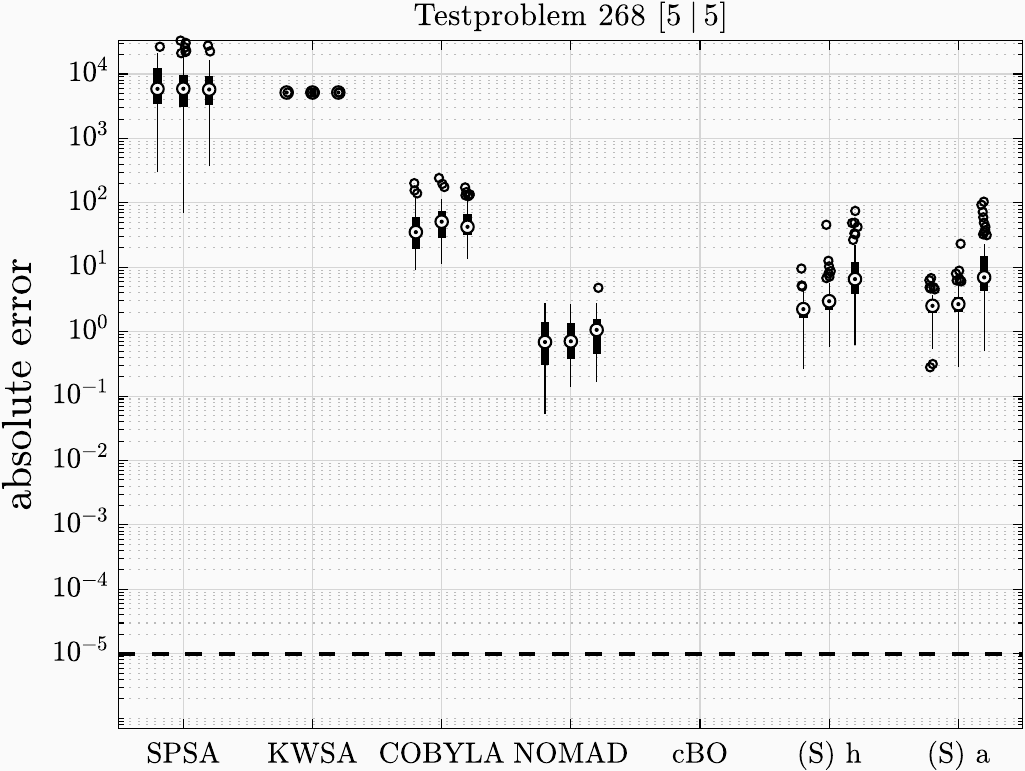}\hspace{2mm}
\includegraphics[width=0.3\textwidth]{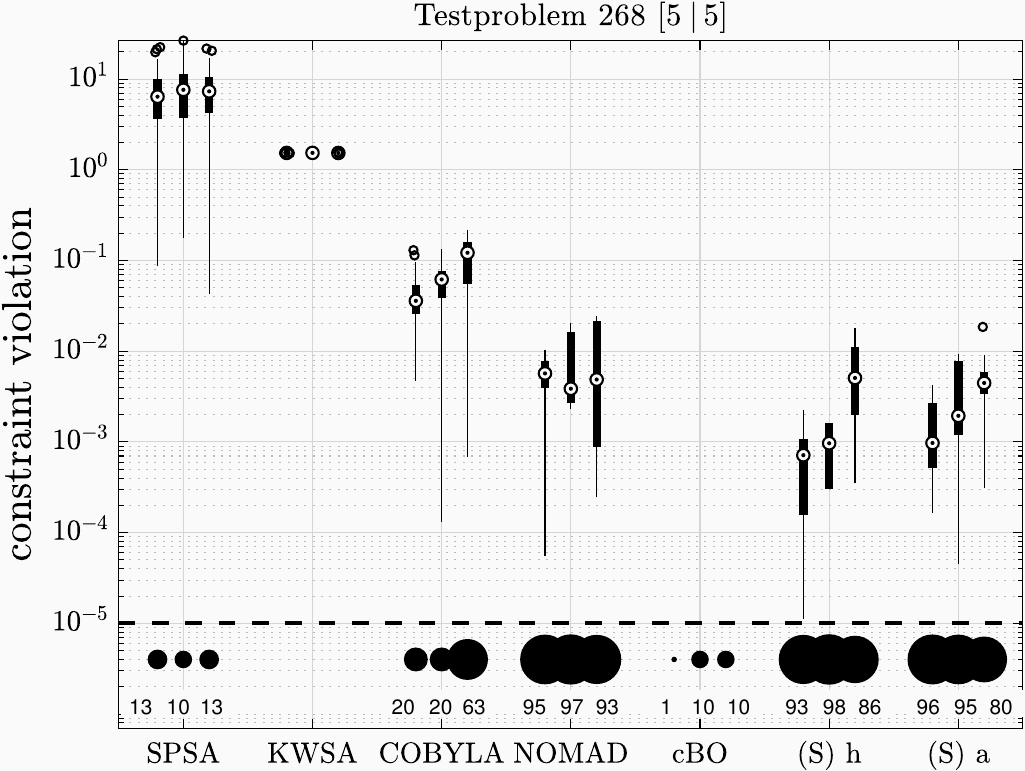}\hspace{2mm}
\includegraphics[width=0.3\textwidth]{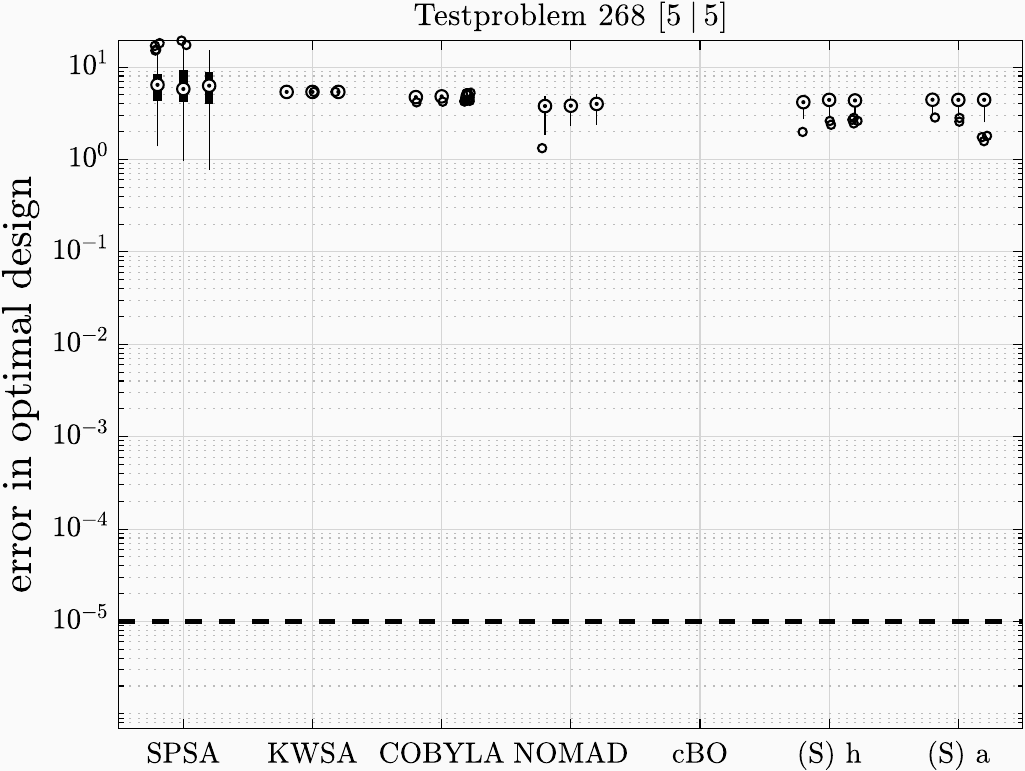}\\[3mm]
\includegraphics[width=0.3\textwidth]{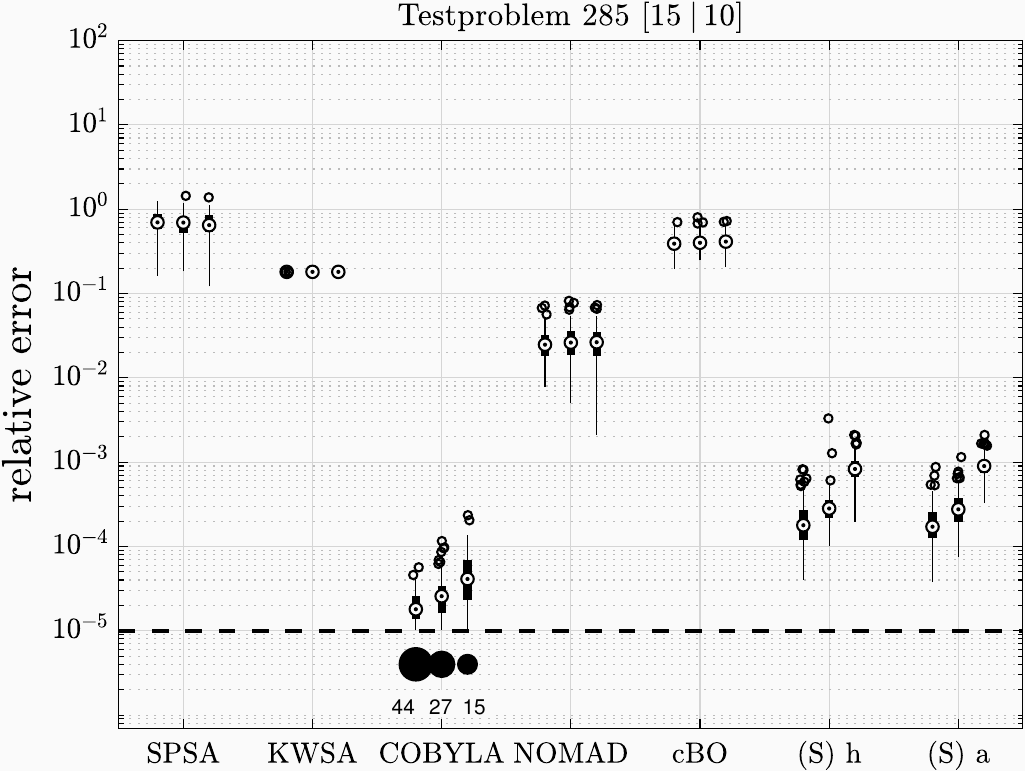}\hspace{2mm}
\includegraphics[width=0.3\textwidth]{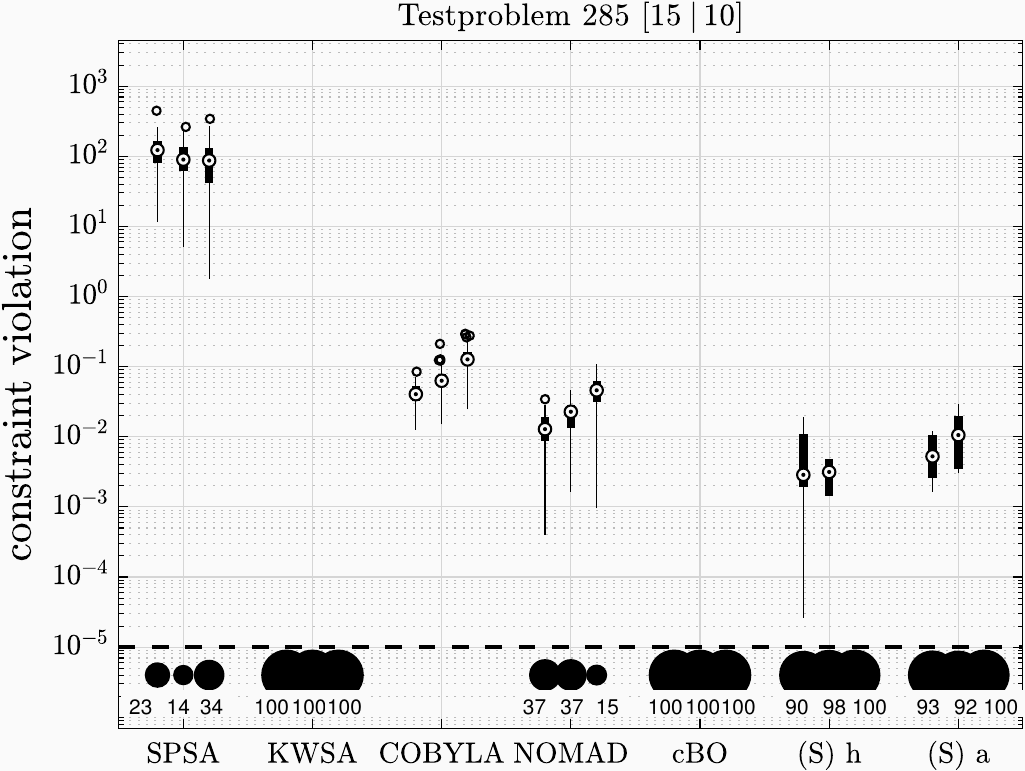}\hspace{2mm}
\includegraphics[width=0.3\textwidth]{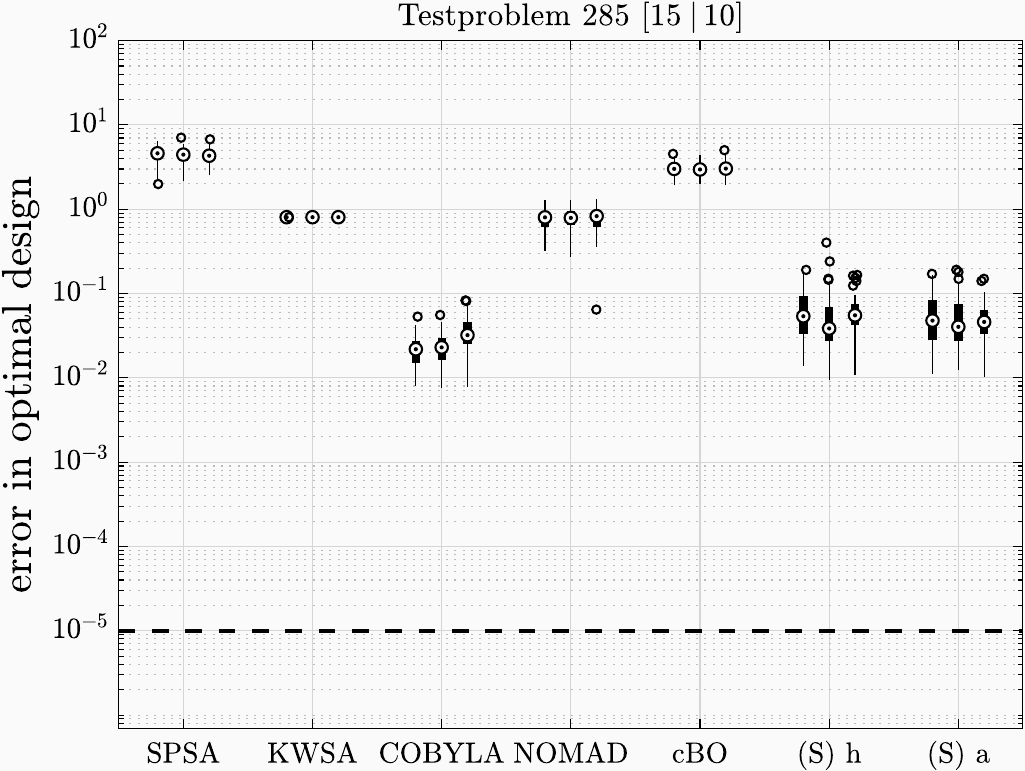}
\caption{Box plots of the errors in the approximated optimal objective values (left plots), the constraint violations (middle plots) and the $l_2$ distance to the exact optimal solution (right plots) of $100$ repeated optimization runs for the Schittkowski test problems number $227$, $228$, $268$, and $285$ for (\ref{eq:robust_expected_value_Schittkowski_problems}). The plots show results of the exact objective function and constraints evaluated at the approximated optimal design computed by (S)NOWPAC, cBO, COBYLA, NOMAD, SPSA and KWSA. All errors or constraint violations below $10^{-5}$ are stated separately below the $10^{-5}$ threshold and the box plots only contain data above this threshold.}\label{fig:schittkowski_testset1b}
\end{center}
\end{figure}

We cannot see, however, a big difference between the two (S)NOWPAC approaches. For further insight we show a detailed accuracy comparison of the individual optimization results at termination, i.e., $250 \cdot N$ black box evaluations in Figs.~\ref{fig:schittkowski_testset1a} - \ref{fig:schittkowski_testset1b}. Here,
we show the accuracy of the optimization results at the approximated optimal points at termination of the optimizers. The plots show the errors in the objective values, the constraint violations and the errors in the approximated optimal designs found by the optimizers at termination respectively. Since the optimal solution for test problem $268$ is zero, we show the absolute error for this test problem. We use MATLAB's box plots to summarize the results for $100$ optimization runs for each benchmark problem for different sample sizes $N \in \{2000, 1000, 200\}$ from left to right separately for each individual robust formulation (\ref{eq:robust_expected_value_Schittkowski_problems})-(\ref{eq:robust_cvar_value_Schittkowski_problems}).
The exact evaluation of the robust objective function and constraints at the approximated optimal designs are shown to eliminate the randomness in the qualitative accuracy of the optimization results. 

We see that (S)NOWPAC most reliably finds accurate approximations to the exact optimal solutions. Note that all optimizers benefit from increasing the number of samples for the approximation of the robustness measures. In (S)NOWPAC, however, the Gaussian process surrogates additionally exploit information from neighbouring points to further reduce the noise, allowing for a better accuracy in the optimization results. 
Additionally, the designs computed by (S)NOWPAC and cBO match well for low-dimensional problems $29$, $227$, $228$, but the accuracy of the results computed by cBO begins to deteriorate in dimensions larger than $4$. This has two reasons: firstly, the global search strategy aims at variance reduction within the whole search domain. This requires more function evaluations than local search. Secondly, the global nature of the Gaussian processes requires a suitable choice of kernels that fits to the properties of the optimization problems, i.e. non-stationarity of the optimization problem, which is not the case in all benchmark problems. Additionally, global maximization of the expected constrained improvement function in every step of the optimization procedure becomes very costly and becomes significant for more than $250$ design points where the Gaussian process evaluation becomes a dominant source of computational effort. To reduce computational costs, approximate Gaussian processes can be employed, an improvement that both, cBO and (S)NOWPAC, benefit from. In (S)NOWPAC first implementations of approximate Gaussian process methods are available to the user—namely SoR, DTC and FITC\cite{Quinonero-Candela2005}. Finally, despite tuning the hyper-parameters for the SPSA and KWSA approaches, the results of these optimizers are not satisfactory in most test examples.

The middle plots in Fig.~\ref{fig:schittkowski_testset1a} - \ref{fig:schittkowski_testset1b} show the maximal constraint violations at the approximated optimal designs. Here, (S)NOWPAC's constraint handling, see \cite{Augustin2014}, in combination with the feasibility restoration mode from Section~\ref{sec:relaxed_feasibility_requirement} allows the computation of approximate optimal designs that exhibit only small constraint violations well below the noise level. 
Additionally, the right plots in Figures~\ref{fig:schittkowski_testset1a} - \ref{fig:schittkowski_testset1b} show the error in the approximated optimal designs. We see that (S)NOWPAC yields either comparable or significantly better results than all the other optimization procedures. Additional benchmarks results for a tighter error tolerance of $\epsilon_f = 10^{-4}$ and $\epsilon_c = 10^{-4}$ can be found in ~\autoref{appssec:benchmarkresultsfortolerance} and additional result for the problem formulations \eqref{eq:robust_quantile_value_Schittkowski_problems} and \eqref{eq:robust_cvar_value_Schittkowski_problems} can be found in the~\autoref{appssec:furtherbenchmarkresults}. \todo{FM: Do not forget appendix}

Finally, comparing the analytic and heuristic approach we see that the analytic approach shows a lower relative error in the objective and the optimal design, especially for lower dimensional problems. This is due to the fact that the Gaussian Process surrogates works especially well in lower dimensions and therefore the optimal smoothing parameter is well approximated. The improvement is, e.g., visible in test problem $29$ and $228$ of Figure \ref{fig:schittkowski_testset1a} and Figure \ref{fig:schittkowski_testset1b}. The analytic smoothing, nevertheless, also shows similar or even better results for high dimensional problems $100$, $113$ and $285$. Combined with the results mentioned above this shows the validity of both approaches and the user can decide which one to use based on the problem and the available computational resources.

\section{Conclusions}\label{sec:conclusions}
We proposed a new stochastic optimization framework SNOWPAC based on the derivative-free trust-region method NOWPAC. The resulting optimization procedure is capable of handling noisy black box evaluations of the objective function and the constraints, which is of particular interest for, but not limited to, robust stochastic optimization problems as discussed in Section~\ref{sec:robustness_formulation}.

Existing approaches for handling noisy constraints either rely on increasing accuracy of the black box evaluations or on Stochastic Approximation \cite{Wang2003}. Increasing the accuracy of individual evaluations of the robustness measures may not be an efficient usage of computational effort as in local approaches individual black box evaluations are often discarded. We therefore introduced Gaussian process surrogates to reduce the noise in the black box evaluations by re-using all available information. This is in contrast to Stochastic Approximation techniques \cite{Kushner1997} which only work with local gradient approximations, disregarding available information. Despite the rich convergence theory for Stochastic Approximation approaches, their practical application often strongly depends of the choice of technical parameters for step and stencil sizes as well as a penalty scheme for handling constraints. Bayesian optimization techniques, in contrast make full use of all available data, resulting in computationally expensive optimization methods, in particular in higher dimensions. (S)NOWPAC combines the advantages of both worlds by utilizing fast local optimization with Gaussian process corrected black box evaluations. We showed in Section~\ref{sec:numerical_results} that the overall performance of (S)NOWPAC is superior to existing optimization approaches by showing improved results for the same computational budget.

In our future work we will investigate convergence properties of our proposed stochastic derivative-free trust-region framework towards a first order critical points.
\section*{Acknowledgements}
This work was partially supported by BP under the BP-MIT Conversion Research Program.

\bibliographystyle{unsrtnat}
\bibliography{references}

\newpage
\section{Appendix}
\subsection{Quantile sampling estimator}
\label{appssec: quantilesamplingestimator}
In this excursion we discuss smoothness properties of the robustness measures $\mathcal{R}_3^{b,\beta}$ and $\mathcal{R}_4^{b,\beta}$. In Example~\ref{thm:example_chance_constraints} we show that $\mathcal{R}_3^{b,\beta}$ often exhibits large curvatures or even non-smoothness in $x$, creating a challenge for approximating this robustness measure using surrogate models. We therefore use the quantile reformulation $\mathcal{R}_3^{b,\beta}$ over the probabilistic constraints $\mathcal{R}_3^{b,\beta}$.

\begin{example}[Non-smoothness of $\mathcal{R}_3^{b,\beta}$]\label{thm:example_chance_constraints} Let us consider the two robust constraints
\begin{align*}
\mathcal{R}_3^{c_1, \beta}(x) &=  \mathbb{E}_{\theta}\left[ \mathds{1}\left\{ x \, : \, \exp\left(\frac{x}2\right)- 16(x-2)^2 \theta^2 + x-1\right\}\right]  -0.1\\
\mathcal{R}_3^{c_{2},\beta}(x) &= \mathbb{E}_{\theta}\left[\mathds{1}\{ x \, : \, 30 x + \theta) \geq 0\}\right] - 0.1
\end{align*}
with $\theta \sim \mathcal{N}(0,1)$ and $\beta = 0.9$. We compute the sample average estimator using $1000$ samples and plot the robustness measures $\mathcal{R}_3^{c_1,\beta}$ (top left) and $\mathcal{R}_3^{c_{2},\beta}(x)$ (bottom left) in Figure~\ref{fig:smoothness_of_robustness_measures}.
\begin{figure}[!htb]
\includegraphics[width=0.45\textwidth]{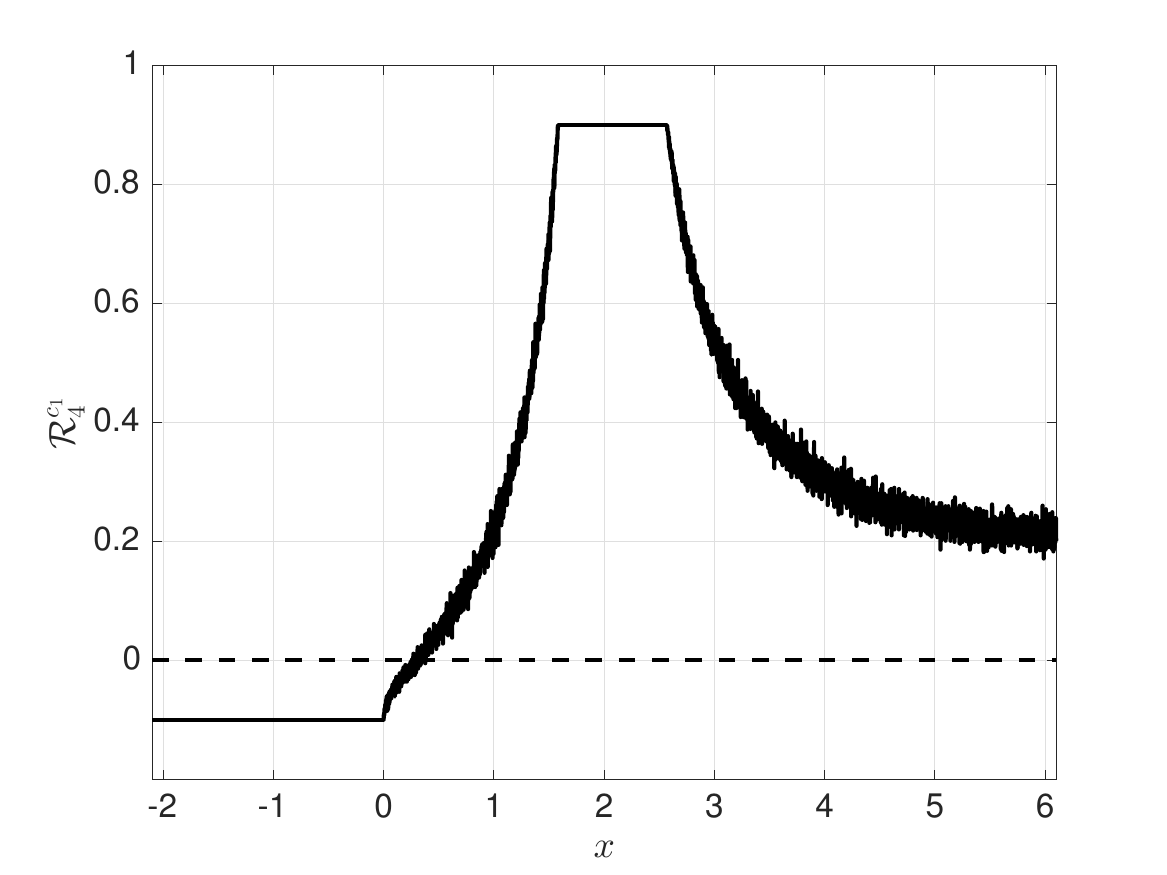}\hfill
\includegraphics[width=0.45\textwidth]{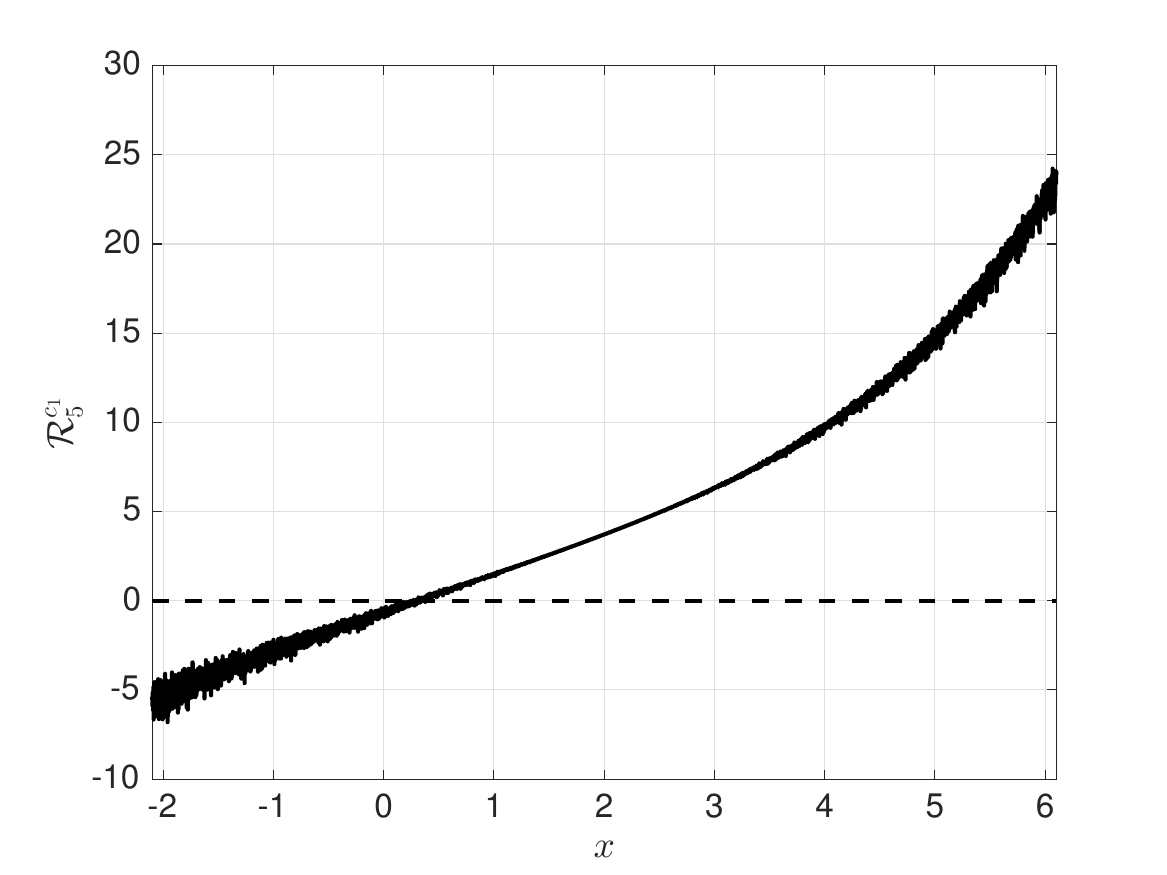}\\
\includegraphics[width=0.45\textwidth]{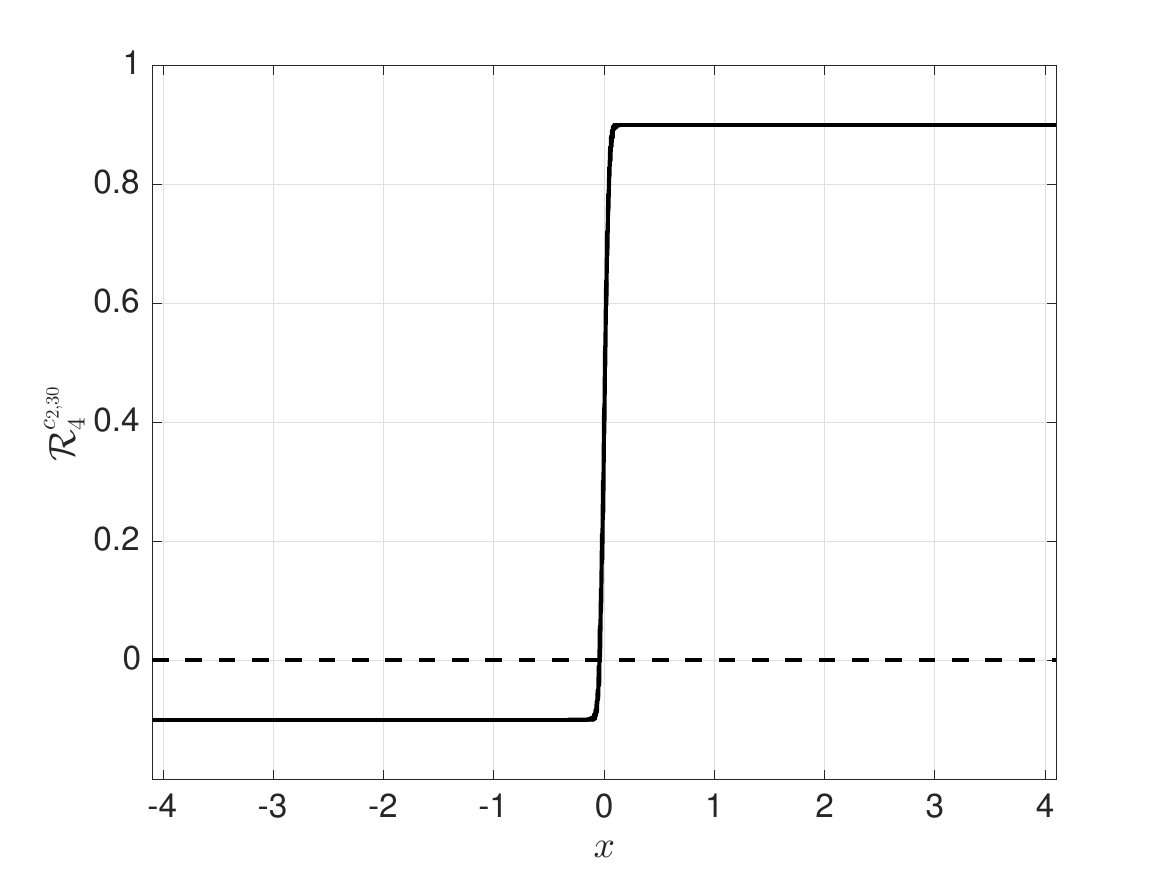}\hspace{12mm}
\includegraphics[width=0.45\textwidth]{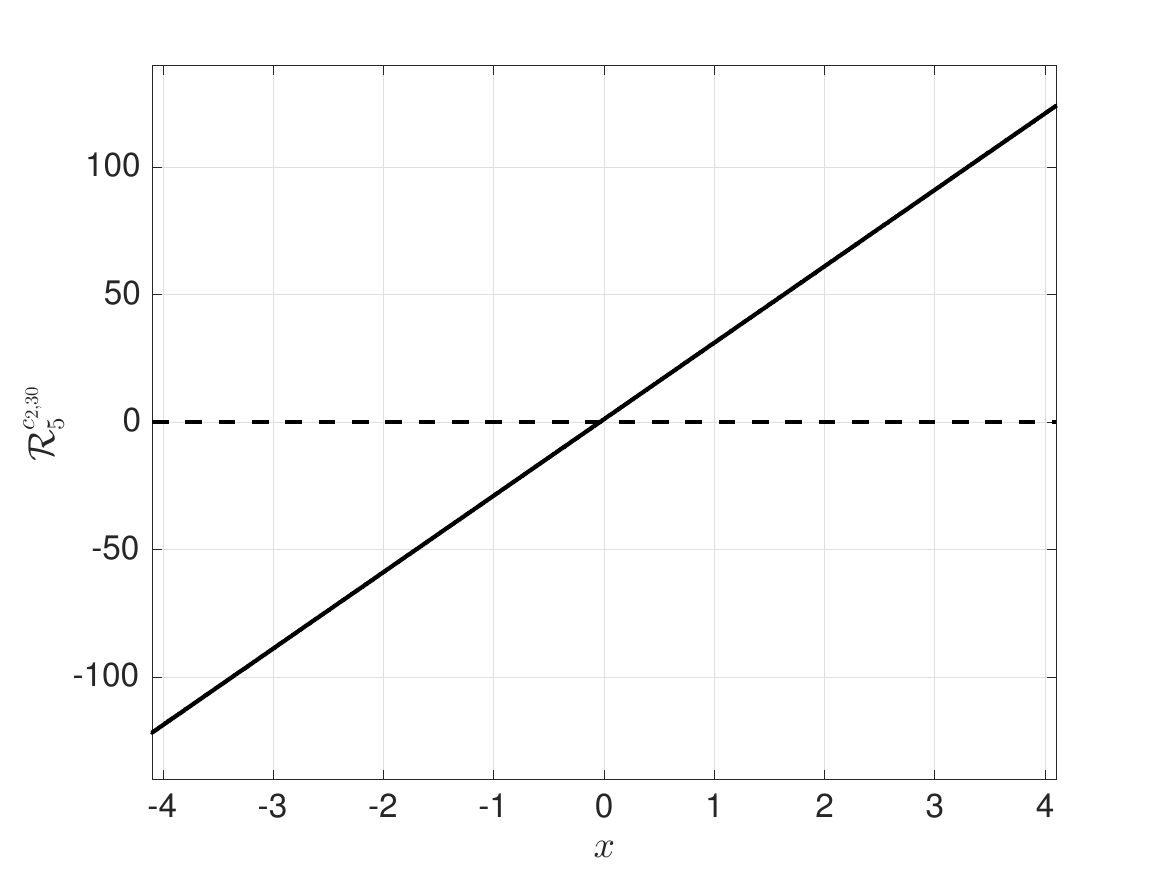}
\caption{Sample approximation of $\mathcal{R}_3^{c_1,0.9}$ (upper left) and $\mathcal{R}_4^{c_1,0.9}$ (upper right) based on resampling $1000$ samples at each $x$. The thresholds $0$ are plotted as dashed lines. The lower plots show $\mathcal{R}_3^{c_2,0.9}$ (left) and $\mathcal{R}_4^{c_2,0.9}$ (right).}\label{fig:smoothness_of_robustness_measures}
\end{figure}
Besides the sample noise we observe that the response surface of $\mathcal{R}_3^{c_1,\beta}$ has kinks at $x \approx 0$, $x \approx 1.5$ and $x \approx 2.5$ which violates the smoothness assumptions on the constraints; for an in depths discussion about smoothness properties of probability distributions we refer to \cite{Kibzun1998, Uryasev1995, Uryasev2000}. Apart from the kinks, even in cases where $\mathcal{R}_3^{c,\beta}$ is arbitrarily smooth, cf. $\mathcal{R}_3^{c_2,\beta}$, it may be a close approximation to a discontinuous step function. The quantile formulations of the probabilistic constraints, $\mathcal{R}_4^{c_1,\beta}$ (top right) and $\mathcal{R}_4^{c_2,\beta}$ (bottom right) in Figure~\ref{fig:smoothness_of_robustness_measures}, on the other hand exhibit smooth behavior.
$\hfill\Diamond$
\end{example}
To approximate the quantile function $\mathcal{R}_4^{b,\beta}$ we can not rely on the standard Monte Carlo estimator for approximating $\mathcal{R}_4^{b,\beta}$ anymore. Instead we follow \cite{Zielinski2004} and use the order statistic $b_{1:N}^x \leq \cdots \leq b_{N:N}^x$, $b_{i:N}^x \in \{b(x, \theta_i)\}_{i=1}^N$ to compute an approximation $b_{\bar{\beta}:N}^x$ of the quantile $b_{\beta}(x)$. More specifically we choose the standard estimator $b_{\bar{\beta}:N}^x \approx b_{\beta}(x)$ with
\[
\bar{\beta} = \left\{\begin{array}{ll}
N\beta, &\mbox{if $N\beta$ is an integer and $\beta < 0.5$}\\ 
N\beta+1, &\mbox{if $N\beta$ is an integer and $\beta > 0.5$}\\ 
\frac N2 + \mathds{1}(U \leq 0), &\mbox{if $N\beta$ is an integer and $\beta = 0.5$}\\
\lfloor N\beta\rfloor + 1, &\mbox{if $N\beta$ is not an integer} 
\end{array}\right.
\]
and $U \sim \mathcal{U}[0, 1]$, yielding 
\begin{equation}\label{eq:simple_quantile_estimator}
b_{\beta}(x) = R_3^{b, \beta}(x) + \varepsilon_x = b_{\bar{\beta}:N}^x + \varepsilon_x.
\end{equation}
Since the order statistic satisfies 
\[
\mu_b[b_{l:N}^x \leq b_{\beta}(x) \leq b_{u:N}^x] \geq \sum\limits_{i=l}^{u-1}\binom{N}{i}\beta^i(1-\beta)^{N-i} =: \pi(l, u, N, \beta)
\]
we use it to define a highly probable confidence interval
$
\left[ b_{l:N}^k , b_{u:N}^k\right];
$ 
see \cite{David2003}.
In the same way as for the sample averages we obtain a highly probable upper bound $\bar{\varepsilon}_x$ on $\varepsilon_x$ by choosing 
\[
\bar{\varepsilon}_x := \max\left\{b_{\bar{\beta}:N}^x - b_{(\bar{\beta}-i):N}^x, 
b_{(\bar{\beta}+i):N}^x - b_{\bar{\beta}:N}^x\right\}
\]
for an $i \in \{1, \ldots, N\}$ such that $\pi\left(\bar{\beta}-i, \bar{\beta}+i, N, \beta\right) \geq \nu$ for the confidence level $\nu \in \,]0, 1[$. 
We refer to \cite{Zielinski2009} for a detailed discussion about optimal quantile estimators.\\

\subsection{Benchmark results for tolerance $1e-4$}
\label{appssec:benchmarkresultsfortolerance}
\begin{figure}[!htb]
\begin{center}
\includegraphics[width=0.45\textwidth]{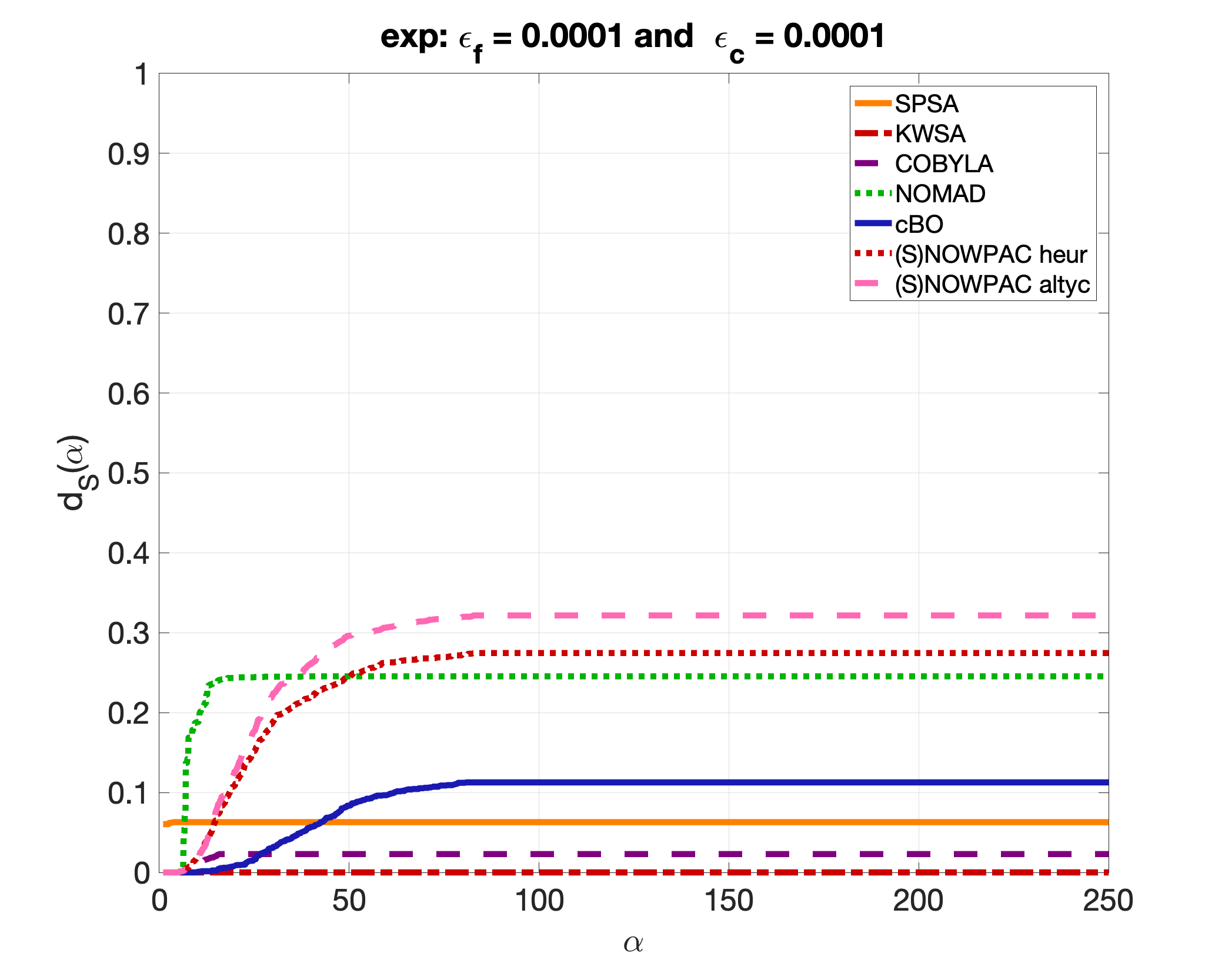}
\includegraphics[width=0.45\textwidth]{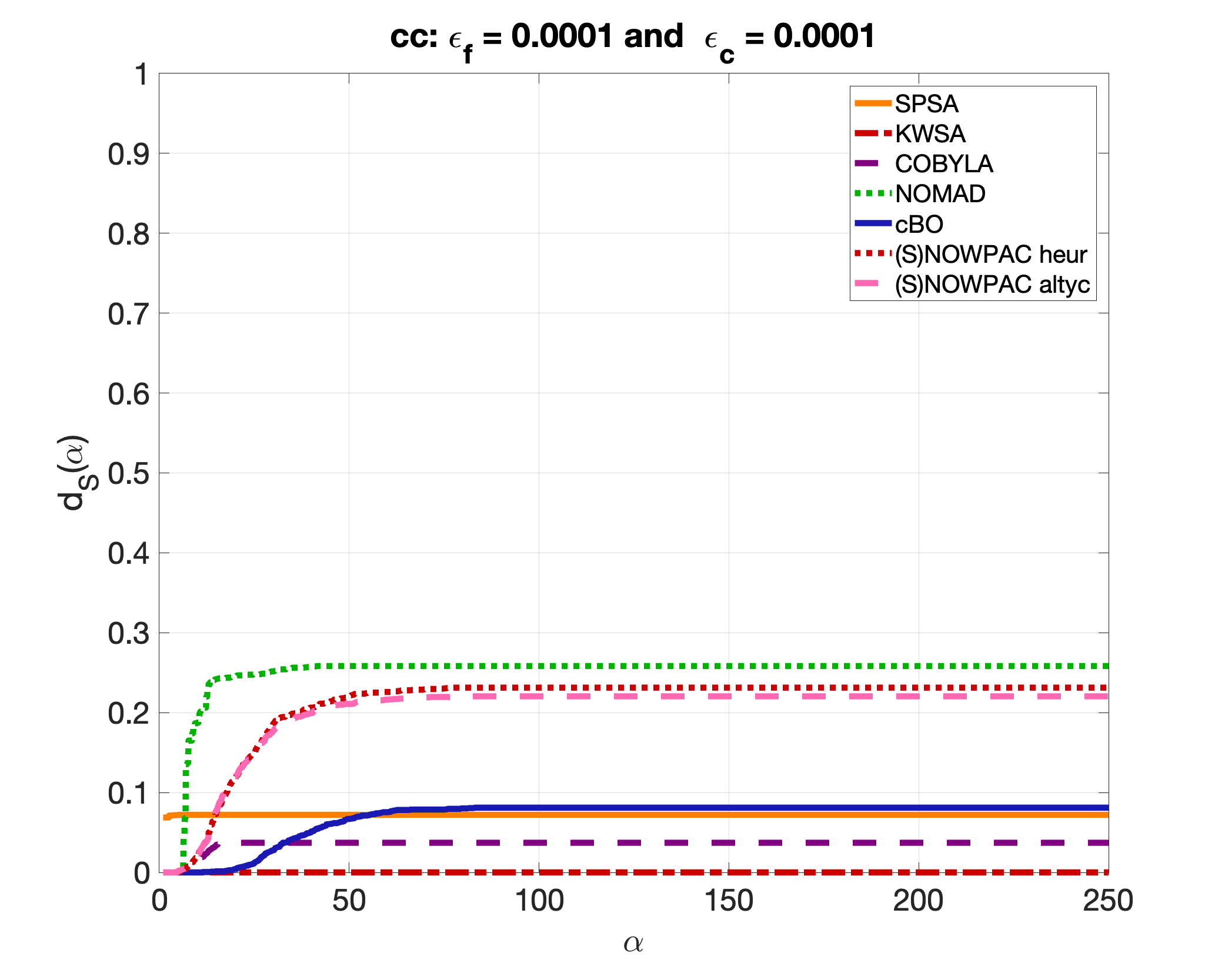}
\includegraphics[width=0.45\textwidth]{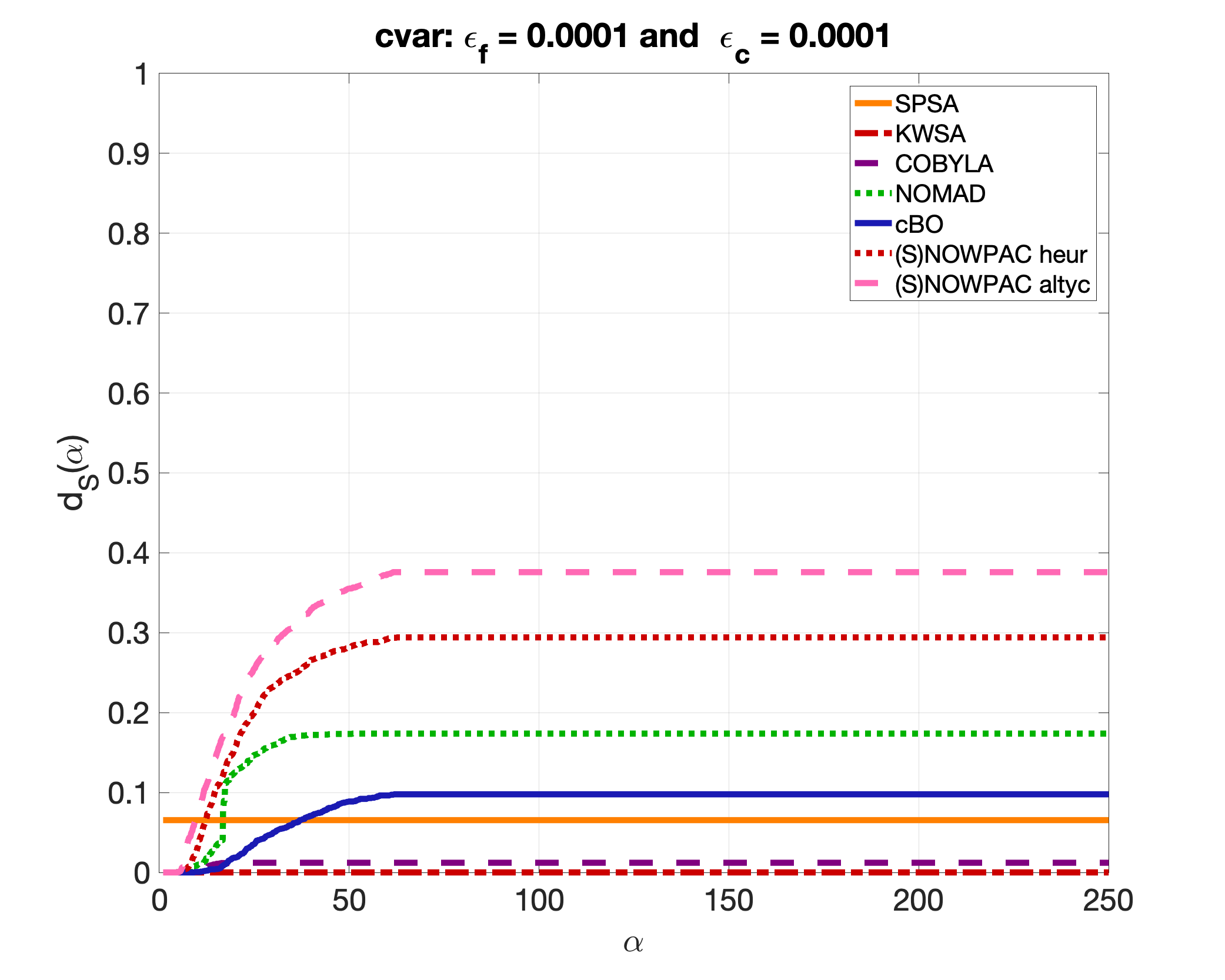}
\end{center}
\caption{Data profiles for (S)NOWPAC altyc (pink), (S)NOWPAC heur (red), cBO (blue), COBYLA (purple), NOMAD (green), SPSA (orange) and KWSA (dark green) of $2400$ runs of the benchmark problems. The results for (\ref{eq:robust_expected_value_Schittkowski_problems}), (\ref{eq:robust_quantile_value_Schittkowski_problems}) and (\ref{eq:robust_cvar_value_Schittkowski_problems}) are plotted in the first, second and third row respectively. The profiles shown are based on the exact values for the objective function and constraints evaluated at the intermediate points computed by the respective optimizers. The data profiles are shown for varying thresholds $\epsilon_f = 10^{-4}$ and $\epsilon_c = 10^{-4}$ on the objective values and the constraint violation respectively.}\label{fig:data_profiles2}
\end{figure}

The results are visualized in~\autoref{fig:data_profiles2}. Due to the small tolerance we now see that only a minority of runs reach the required threshold. Again we see a good performance of (S)NOWPAC where we now see the analytic approach perform the best for (\ref{eq:robust_expected_value_Schittkowski_problems}) and (\ref{eq:robust_cvar_value_Schittkowski_problems}). Since we are in a region of small tolerance and therefore small noise, also the bootstrapping approximation improves and hence improves the final result of the optimization.

\subsection{Further benchmark results for \eqref{eq:robust_quantile_value_Schittkowski_problems} and \eqref{eq:robust_cvar_value_Schittkowski_problems}}
\label{appssec:furtherbenchmarkresults}
In this section we show the full set of benchmark results for the following two robust optimization formulations:
\begin{enumerate}
\item Minimization of the average objective function subject to the constraints being satisfied in $95\%$ of all cases:
\begin{equation}\label{eq:robust_quantile_value_Schittkowski_problems}
\begin{split}
&\;\;\;\min \mathcal{R}_0^f (x) \\
&\mbox{s.t.} \quad \mathcal{R}_4^{c,0.95}(x) \leq 0.
\end{split}
\end{equation}
\item Minimization of the $95\%$-CVaR of the objective function subject to the constraints being satisfied on average:
\begin{equation}\label{eq:robust_cvar_value_Schittkowski_problems}
\begin{split}
&\;\;\;\min \mathcal{R}_5^{f,0.95}(x) \\
&\mbox{s.t.} \quad \mathcal{R}_0^c(x) \leq 0,
\end{split}
\end{equation}
\end{enumerate}

\begin{figure}[!htb]
\begin{center}
\includegraphics[width=0.3\textwidth]{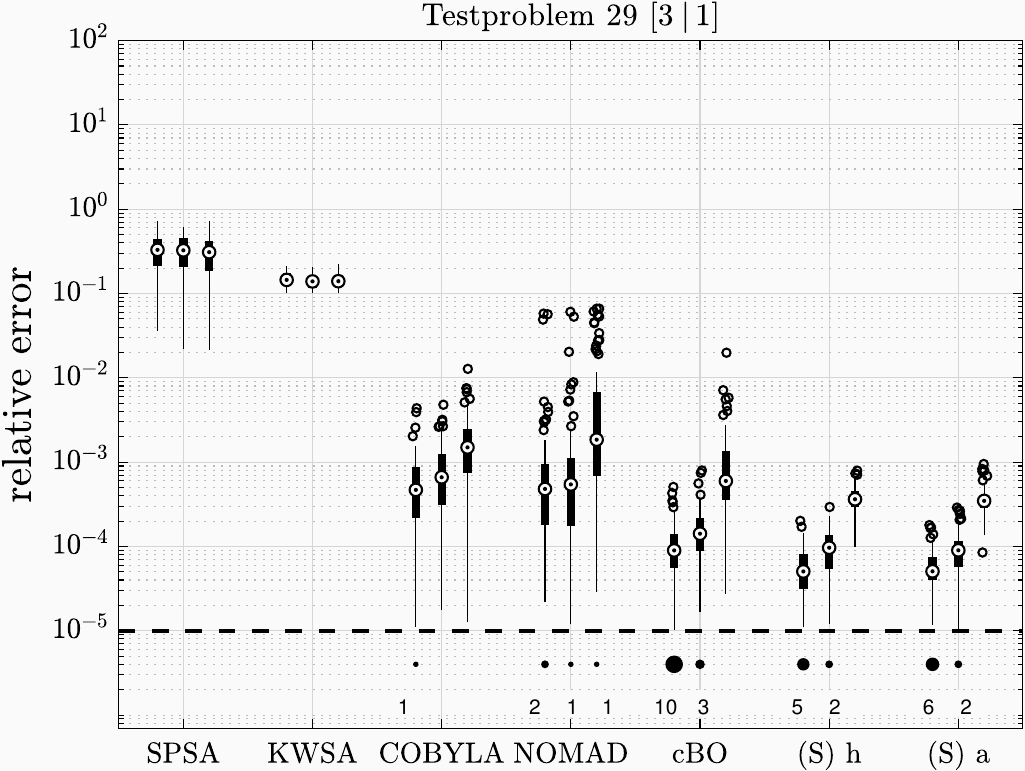}\hspace{2mm}
\includegraphics[width=0.3\textwidth]{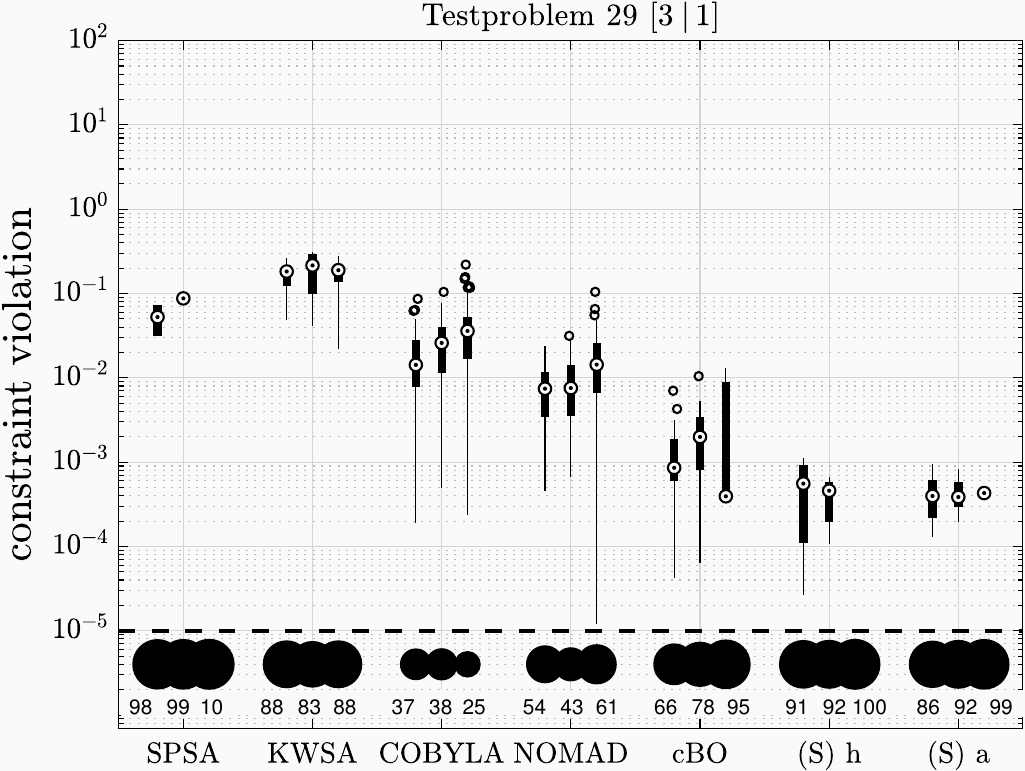}\hspace{2mm}
\includegraphics[width=0.3\textwidth]{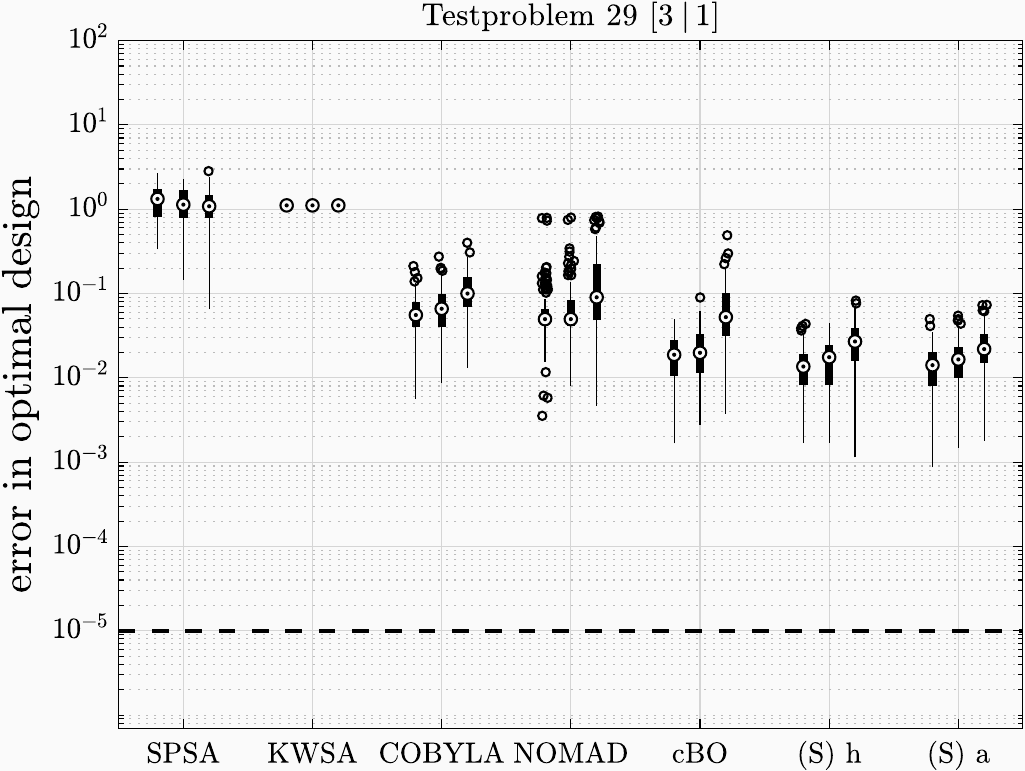}\\[3mm]
\includegraphics[width=0.3\textwidth]{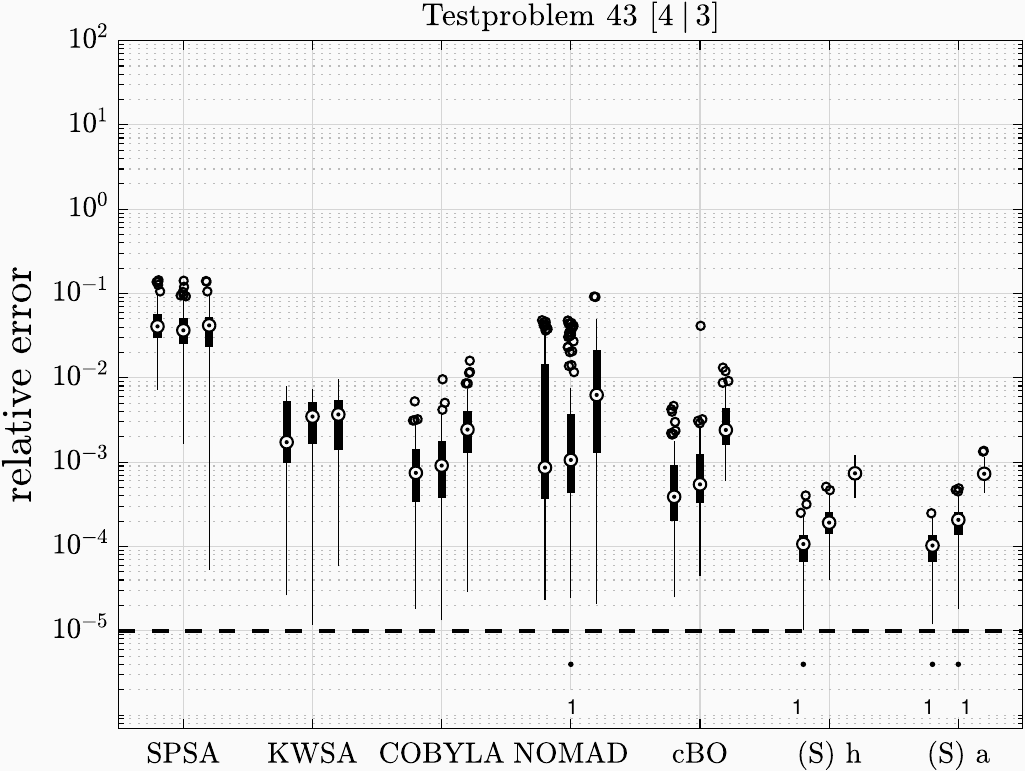}\hspace{2mm}
\includegraphics[width=0.3\textwidth]{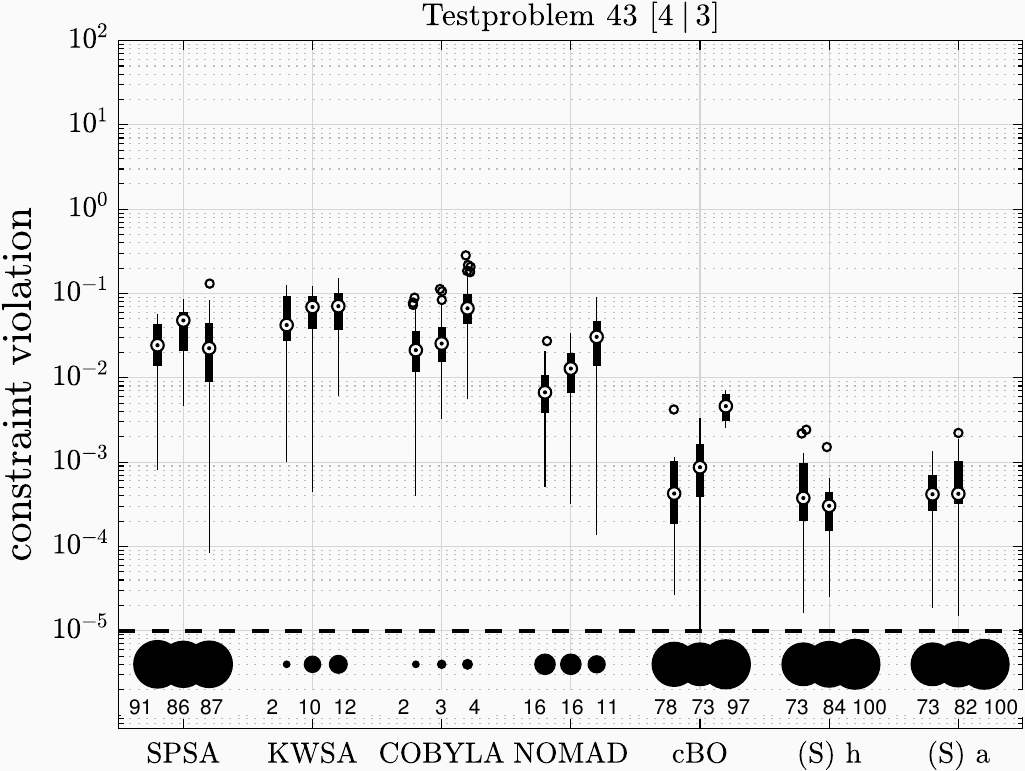}\hspace{2mm}
\includegraphics[width=0.3\textwidth]{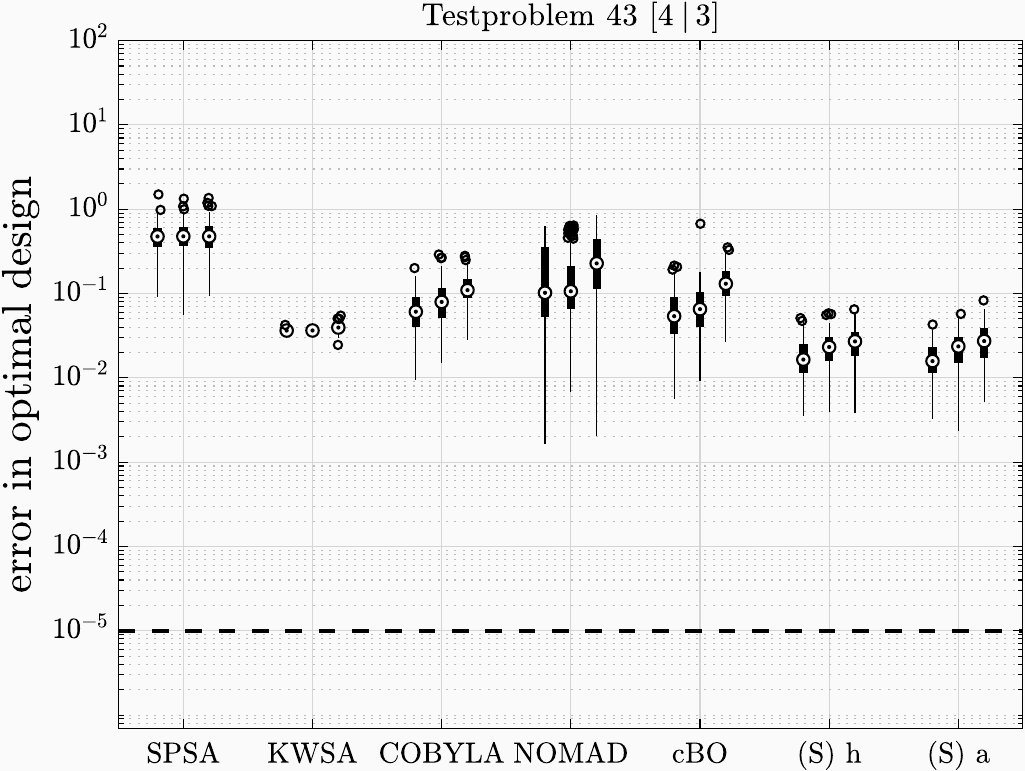}\\[3mm]
\includegraphics[width=0.3\textwidth]{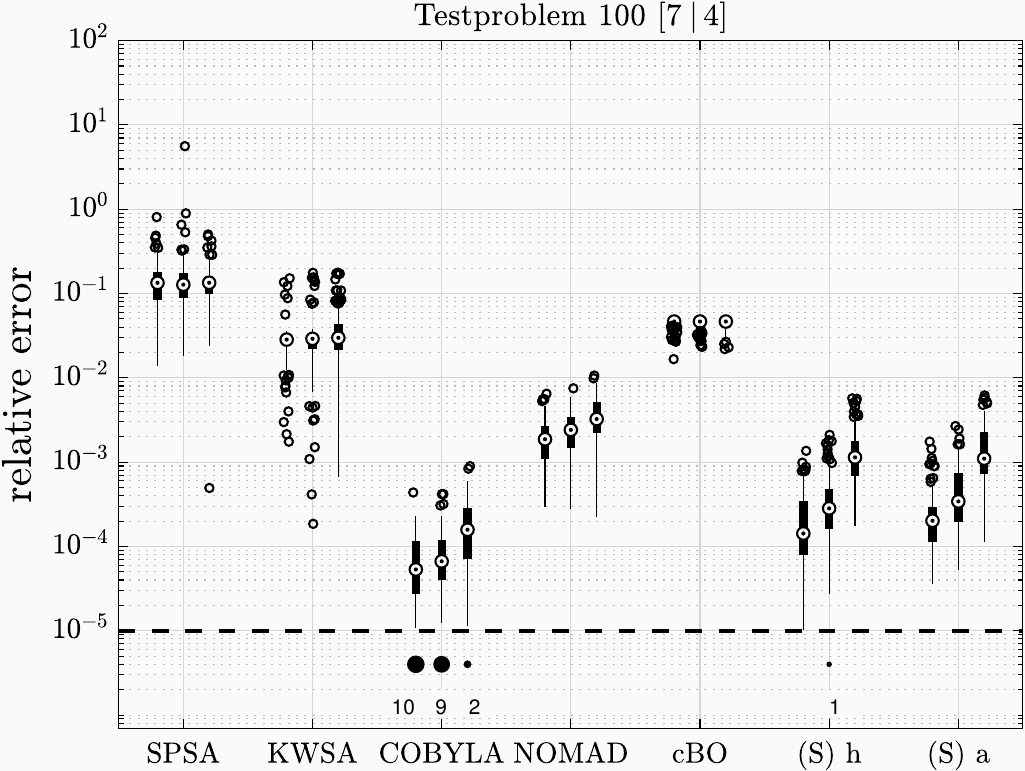}\hspace{2mm}
\includegraphics[width=0.3\textwidth]{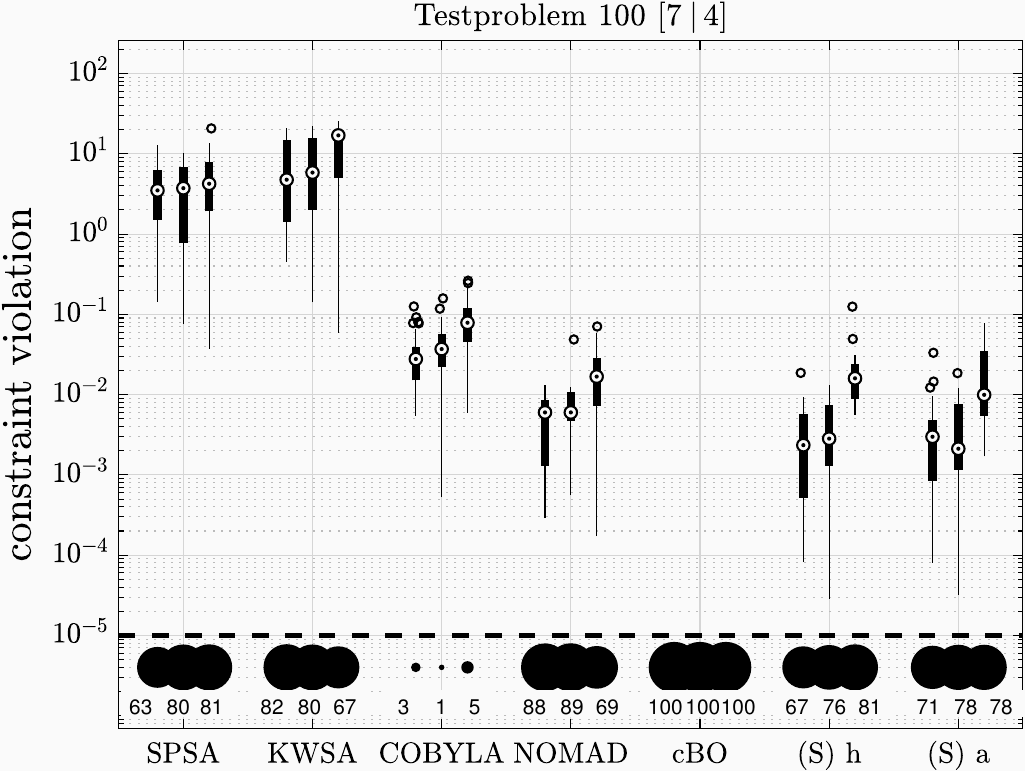}\hspace{2mm}
\includegraphics[width=0.3\textwidth]{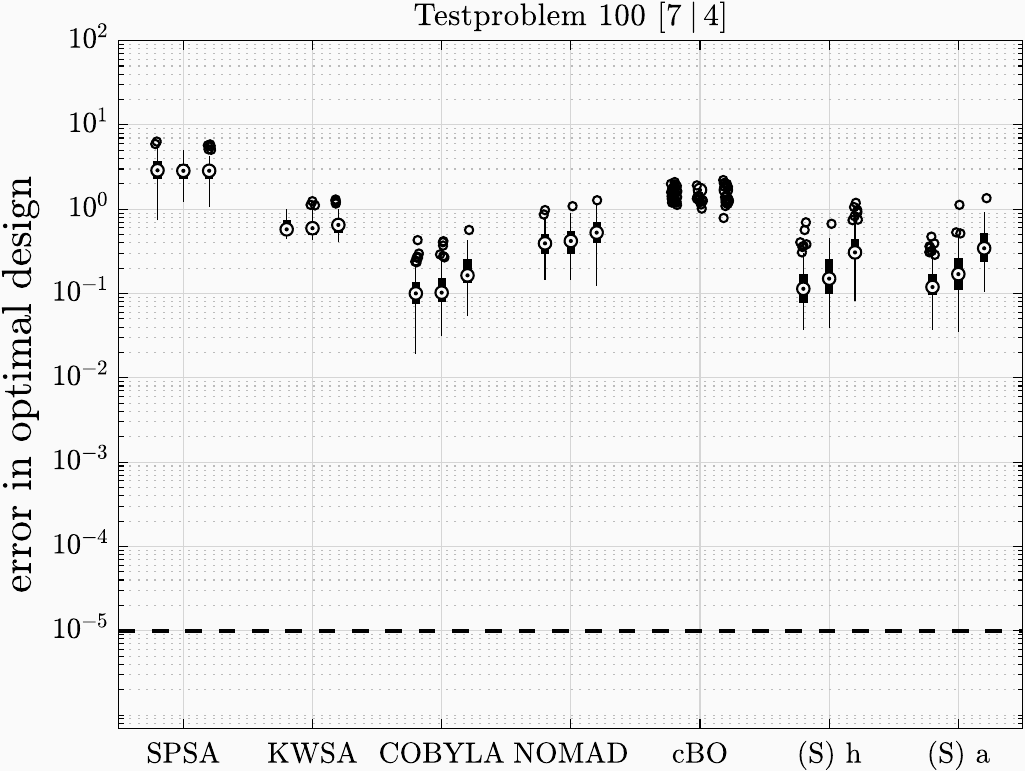}\\[3mm]
\includegraphics[width=0.3\textwidth]{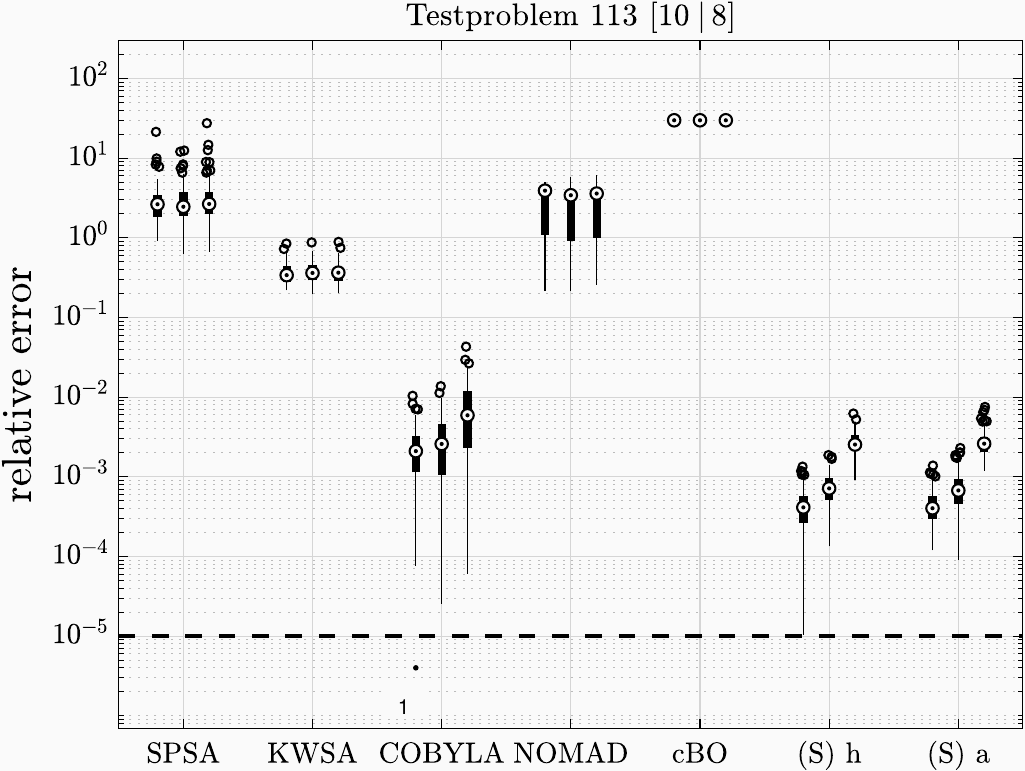}\hspace{2mm}
\includegraphics[width=0.3\textwidth]{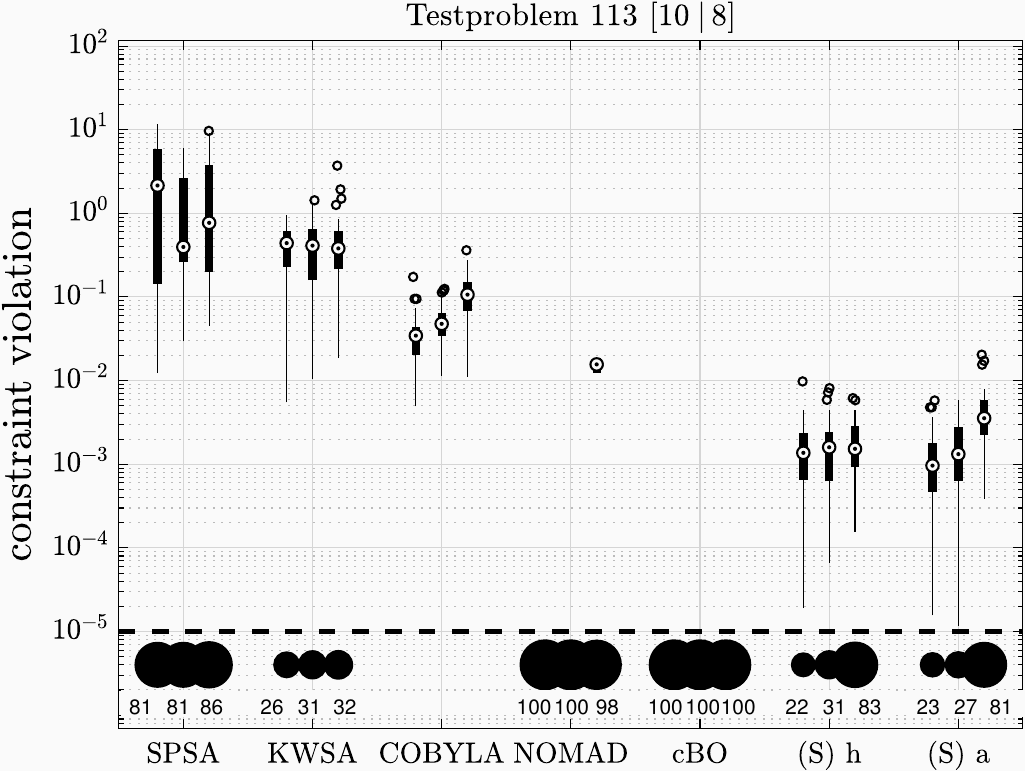}\hspace{2mm}
\includegraphics[width=0.3\textwidth]{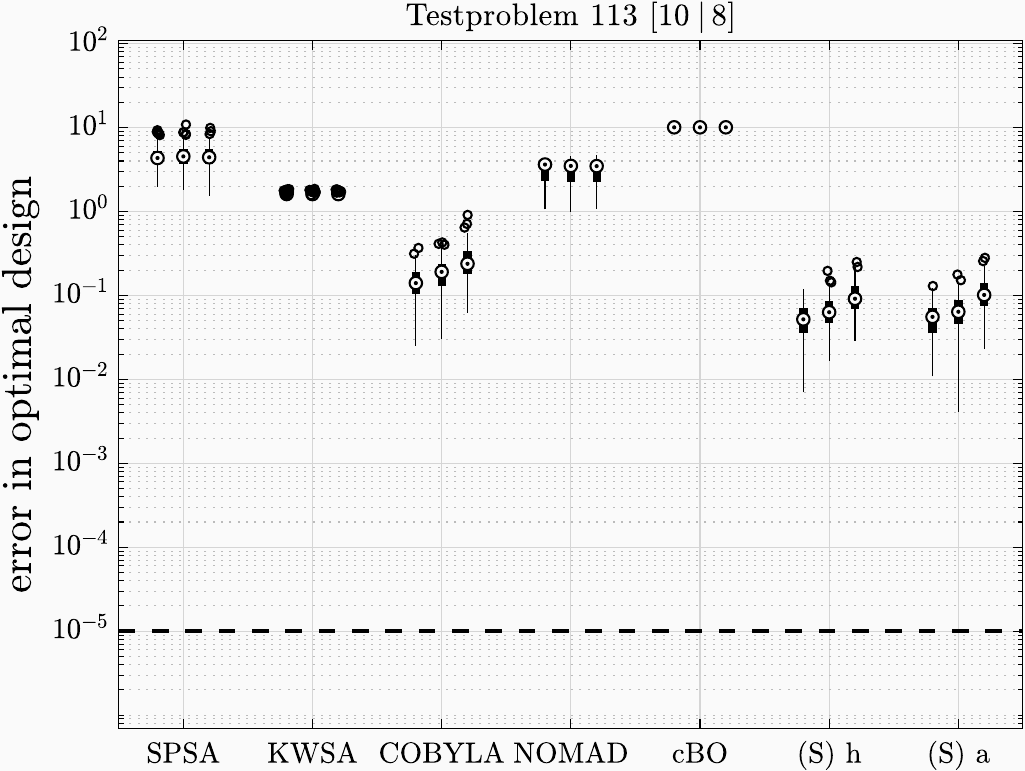}
\caption{Box plots of the errors in the approximated optimal objective values (left plots), the constraint violations (middle plots) and the $l_2$ distance to the exact optimal solution (right plots) of $100$ repeated optimization runs for the Schittkowski test problems number $29$, $43$, $100$, and $113$ for (\ref{eq:robust_quantile_value_Schittkowski_problems}). The plots show results of the exact objective function and constraints evaluated at the approximated optimal design computed by (S)NOWPAC, cBO, COBYLA, NOMAD, SPSA and KWSA. Thereby all errors or constraint violations below $10^{-5}$ are stated separately below the $10^{-5}$ threshold and the box plots only contain data above this threshold.}\label{fig:schittkowski_testset2a}
\end{center}
\end{figure}
\begin{figure}[!htb]
\begin{center}
\includegraphics[width=0.3\textwidth]{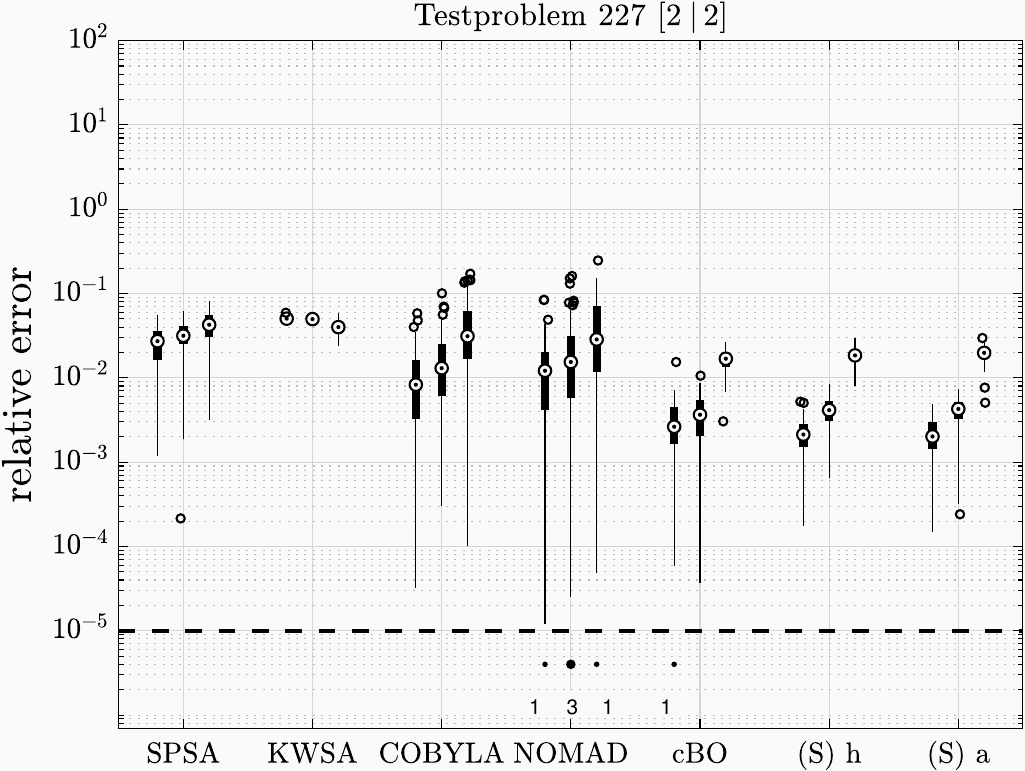}\hspace{2mm}
\includegraphics[width=0.3\textwidth]{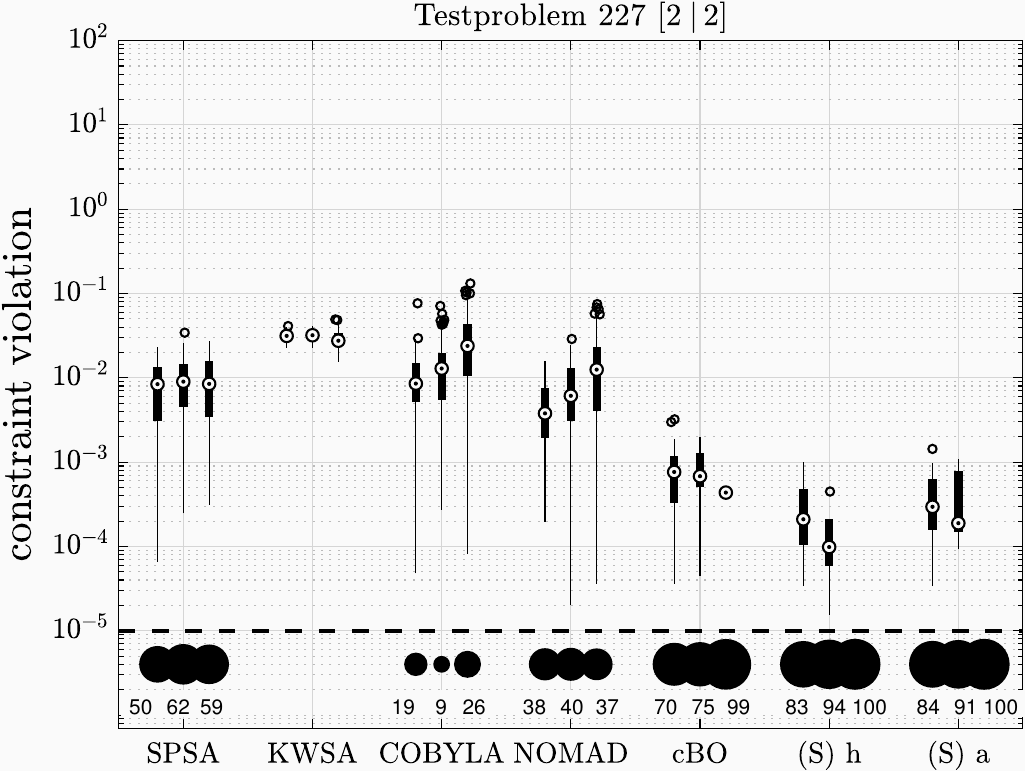}\hspace{2mm}
\includegraphics[width=0.3\textwidth]{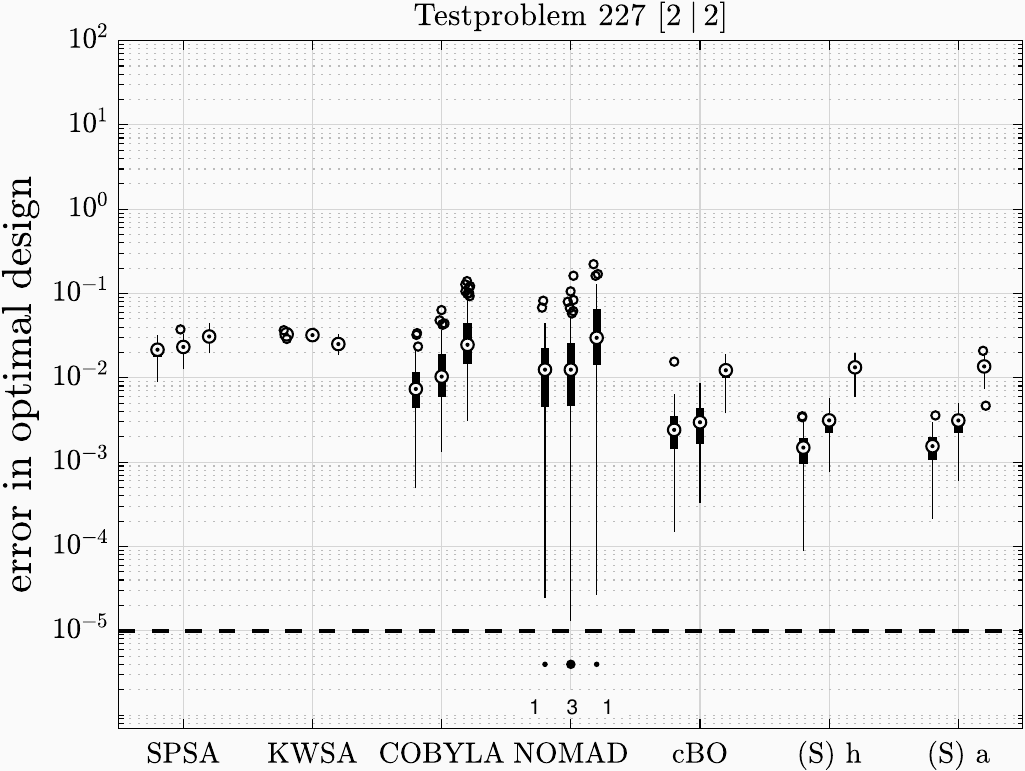}\\[3mm]
\includegraphics[width=0.3\textwidth]{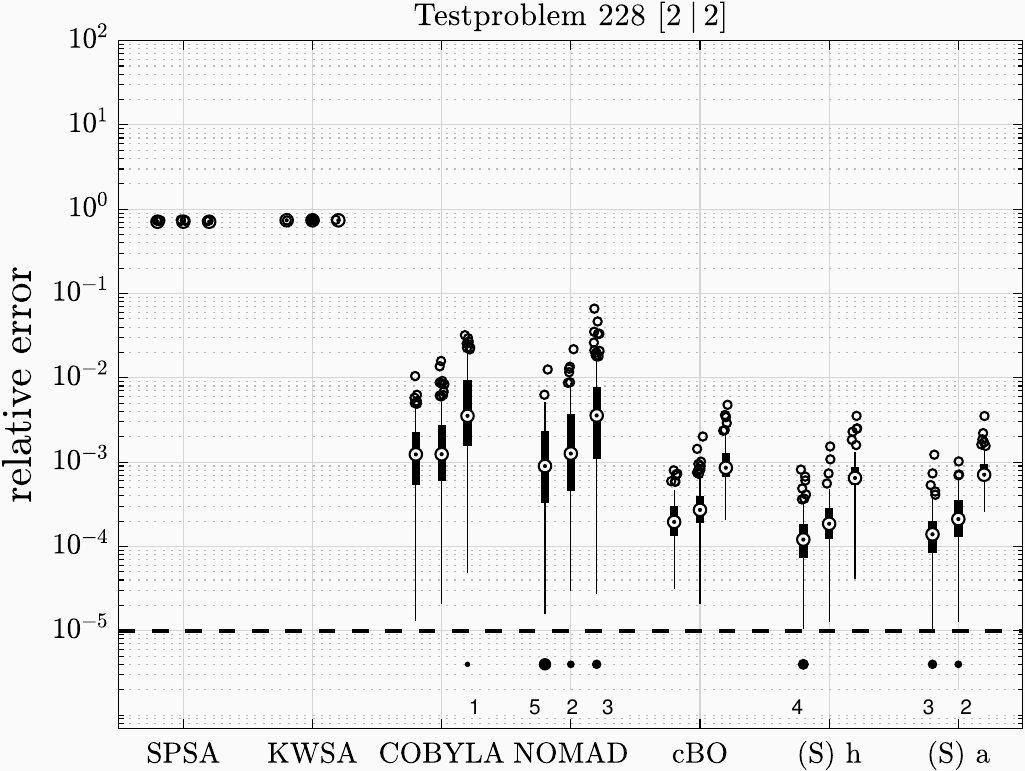}\hspace{2mm}
\includegraphics[width=0.3\textwidth]{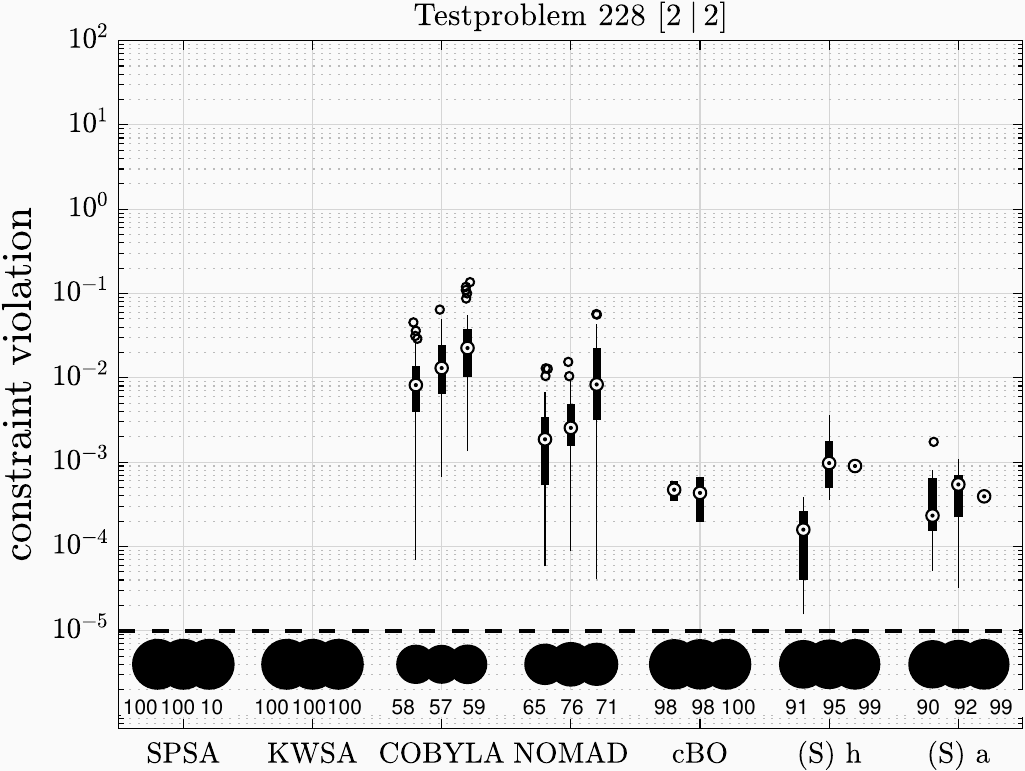}\hspace{2mm}
\includegraphics[width=0.3\textwidth]{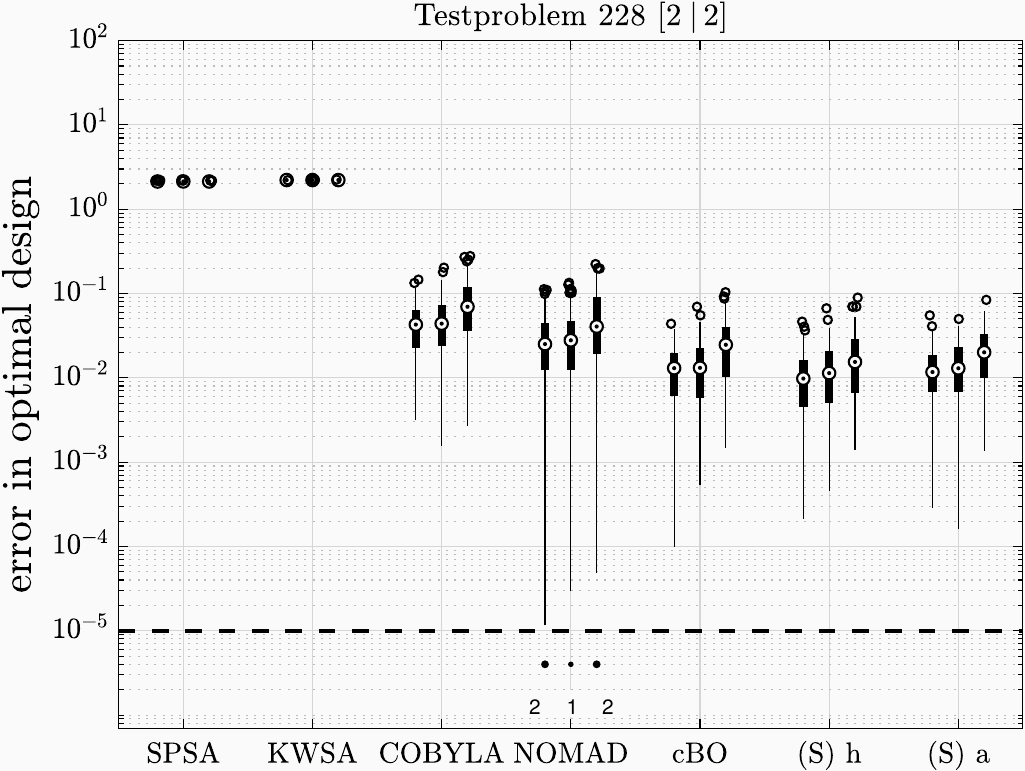}\\[3mm]
\includegraphics[width=0.3\textwidth]{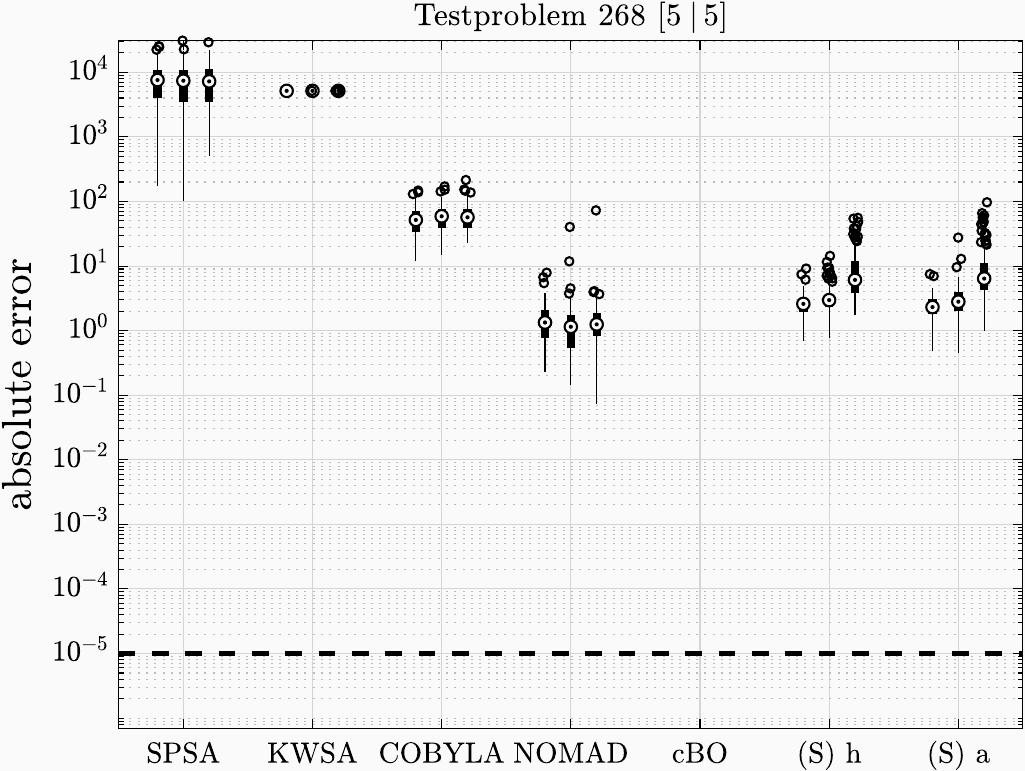}\hspace{2mm}
\includegraphics[width=0.3\textwidth]{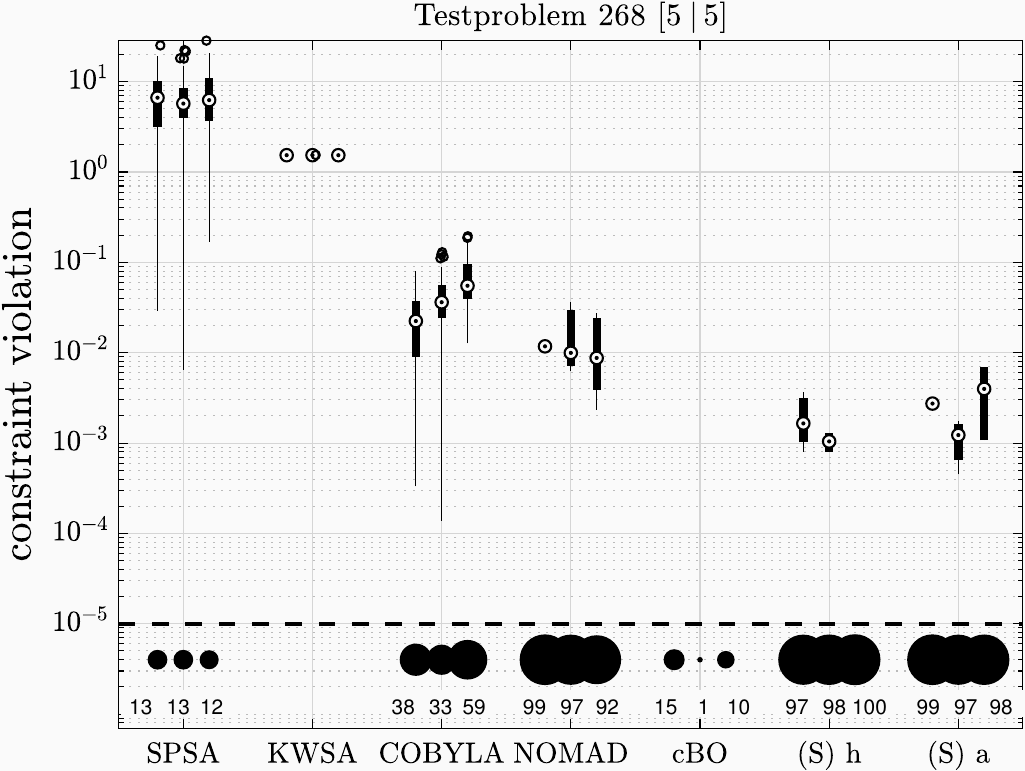}\hspace{2mm}
\includegraphics[width=0.3\textwidth]{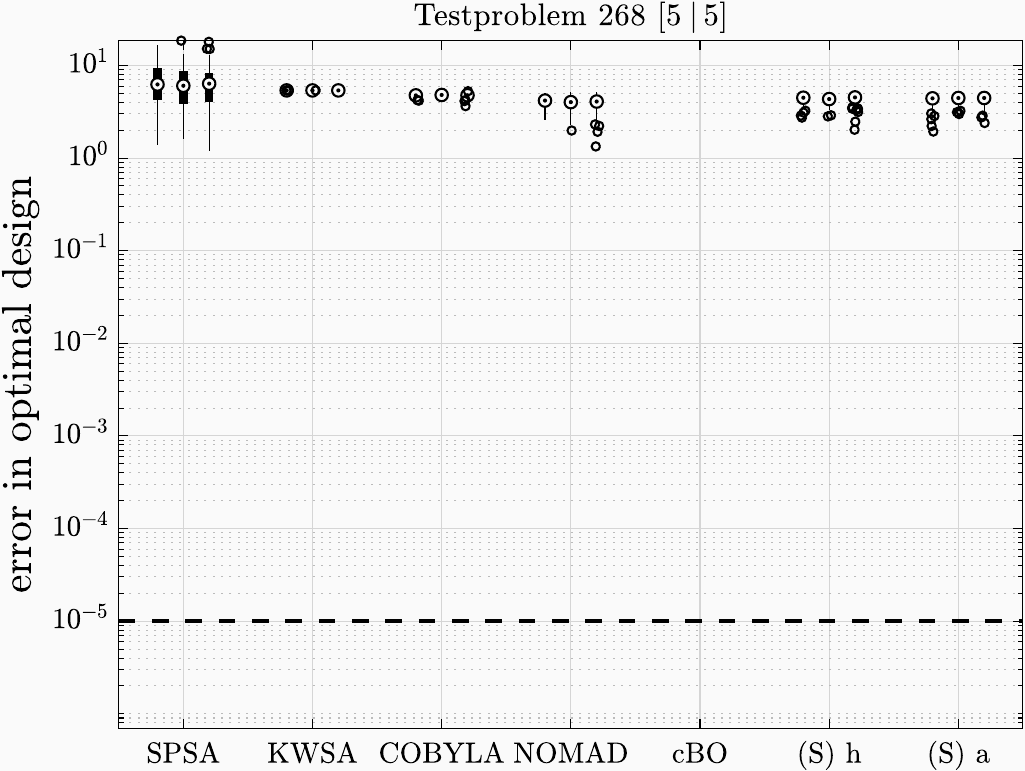}\\[3mm]
\includegraphics[width=0.3\textwidth]{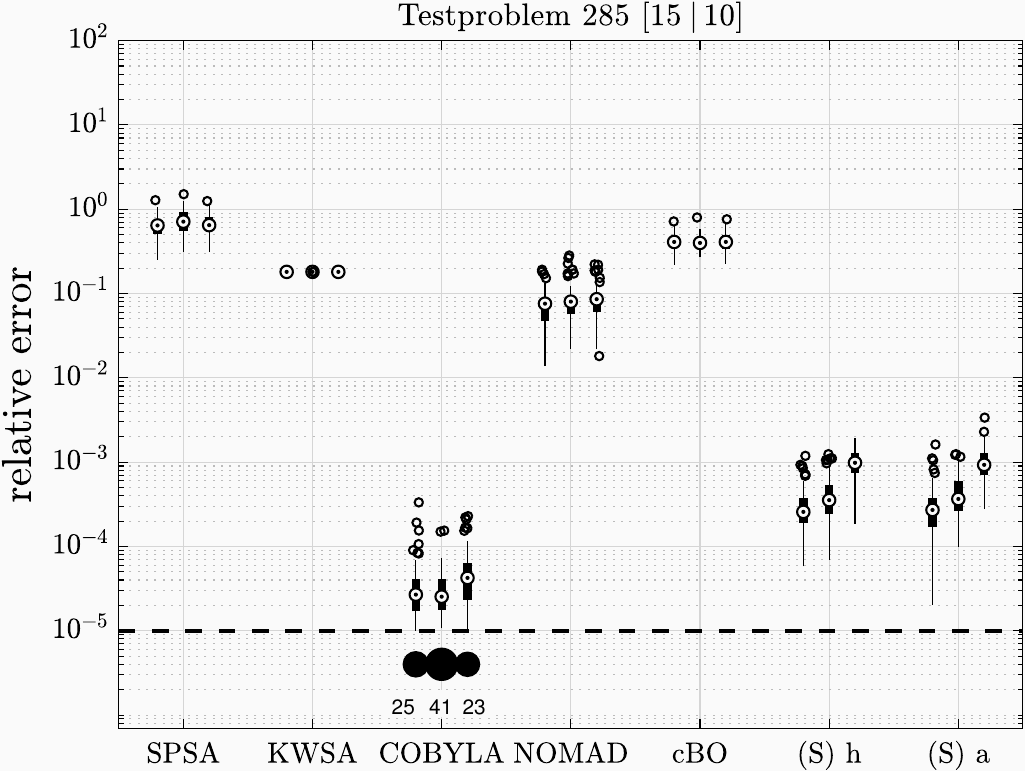}\hspace{2mm}
\includegraphics[width=0.3\textwidth]{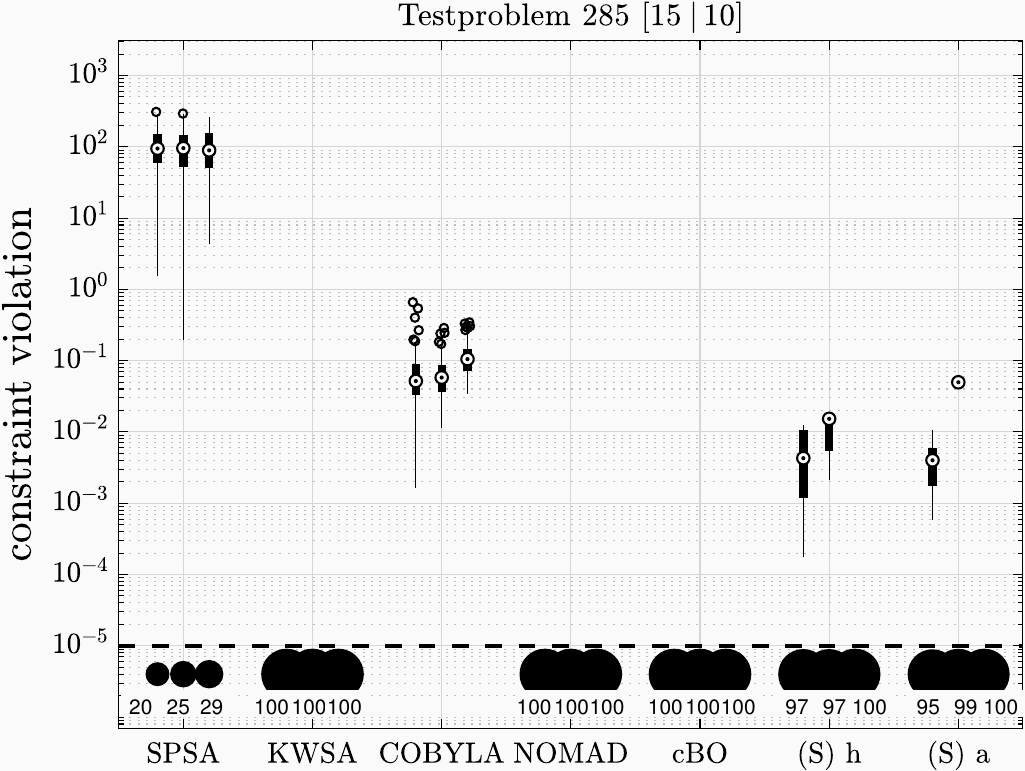}\hspace{2mm}
\includegraphics[width=0.3\textwidth]{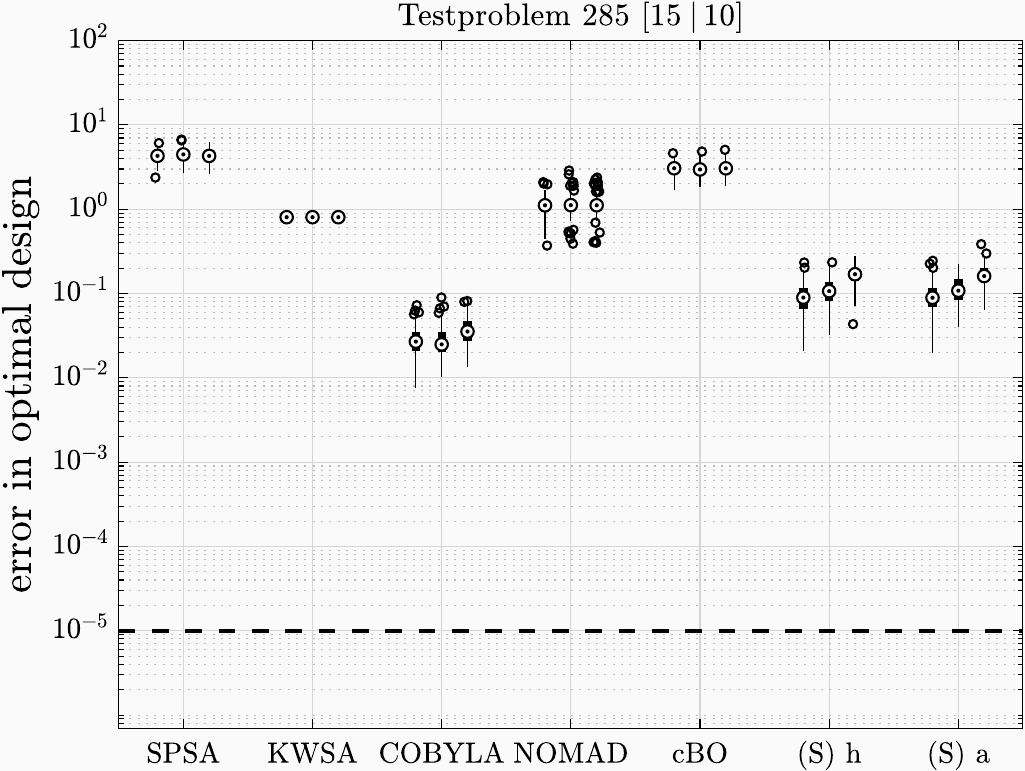}
\caption{Box plots of the errors in the approximated optimal objective values (left plots), the constraint violations (middle plots) and the $l_2$ distance to the exact optimal solution (right plots) of $100$ repeated optimization runs for the Schittkowski test problems number $227$, $228$, $268$, and $285$ for (\ref{eq:robust_quantile_value_Schittkowski_problems}). The plots show results of the exact objective function and constraints evaluated at the approximated optimal design computed by (S)NOWPAC, cBO, COBYLA, NOMAD, SPSA and KWSA. Thereby all errors or constraint violations below $10^{-5}$ are stated separately below the $10^{-5}$ threshold and the box plots only contain data above this threshold.}\label{fig:schittkowski_testset2b}
\end{center}
\end{figure}

\begin{figure}[!htb]
\begin{center}
\includegraphics[width=0.3\textwidth]{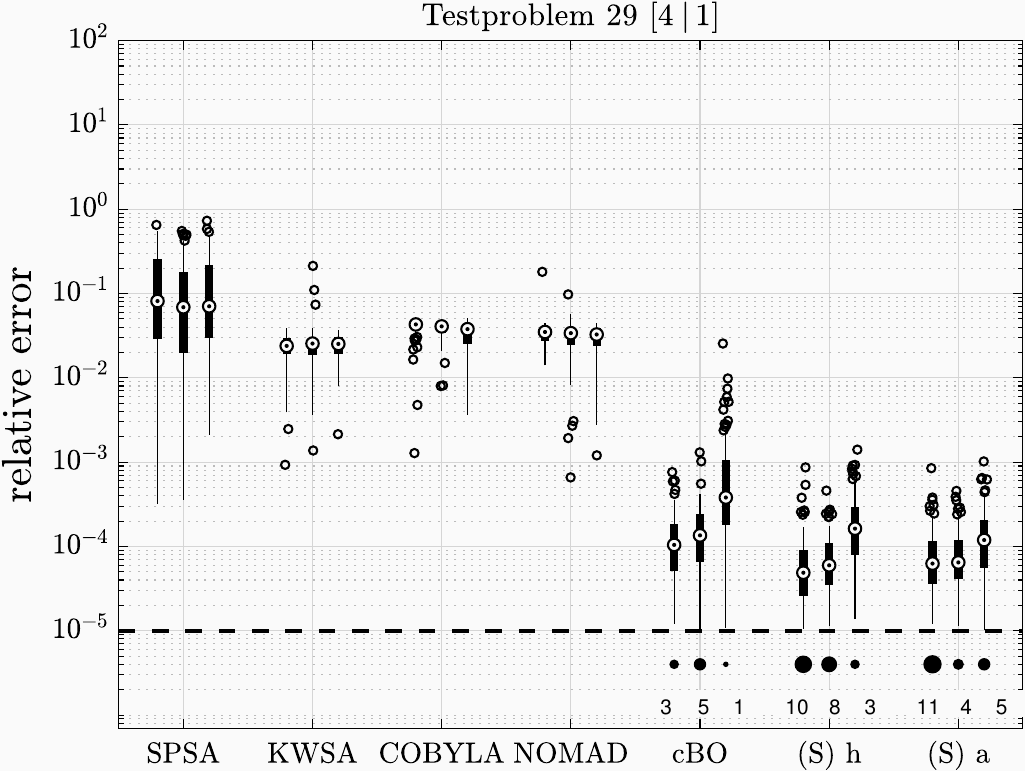}\hspace{2mm}
\includegraphics[width=0.3\textwidth]{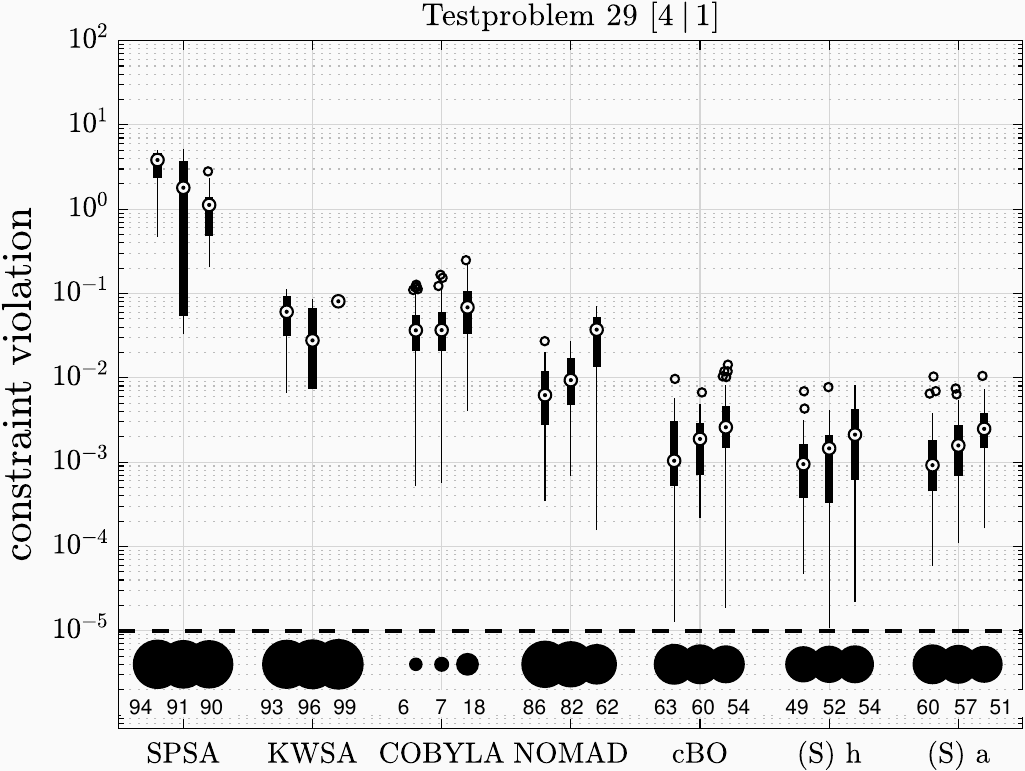}\hspace{2mm}
\includegraphics[width=0.3\textwidth]{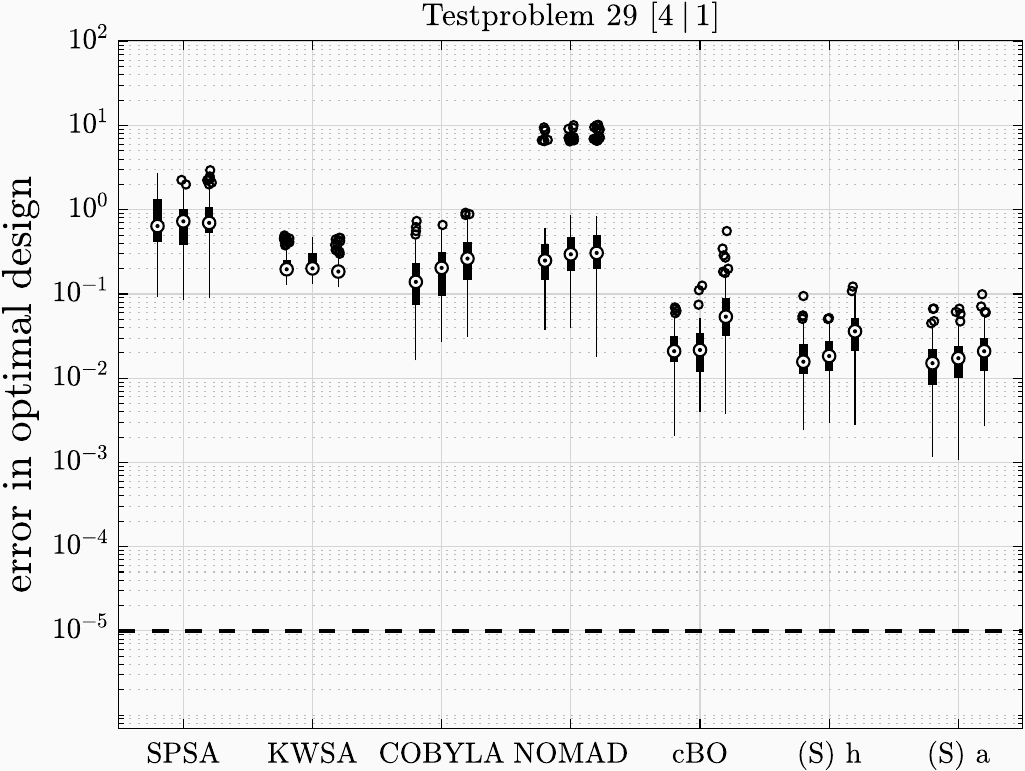}\\[3mm]
\includegraphics[width=0.3\textwidth]{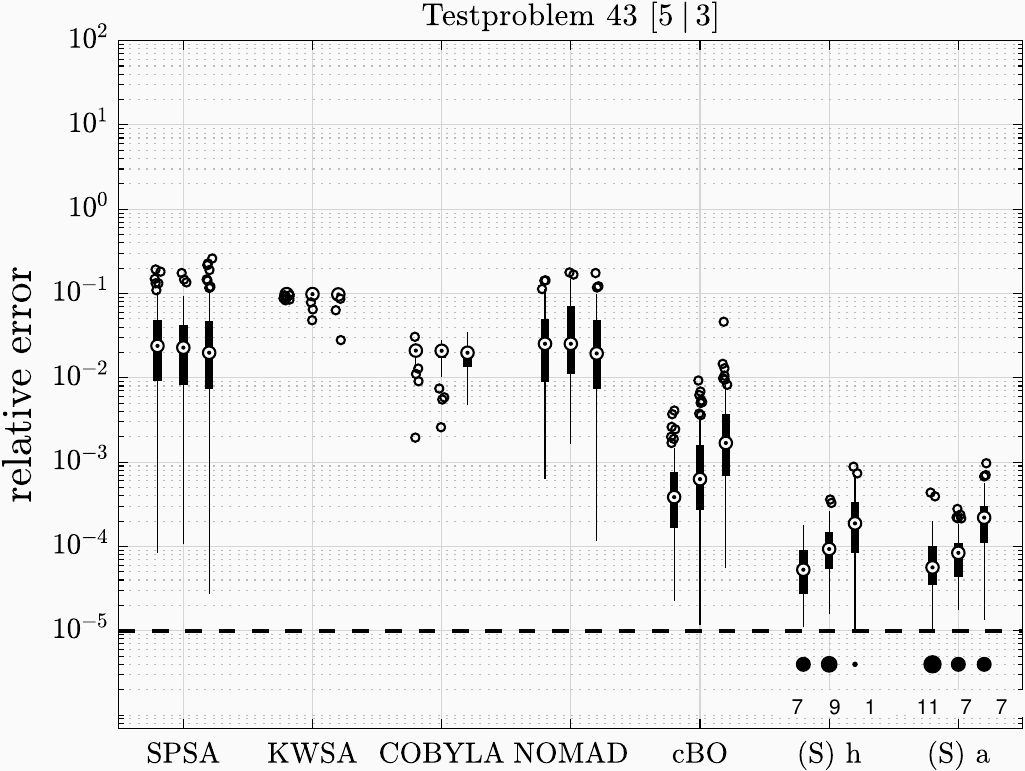}\hspace{2mm}
\includegraphics[width=0.3\textwidth]{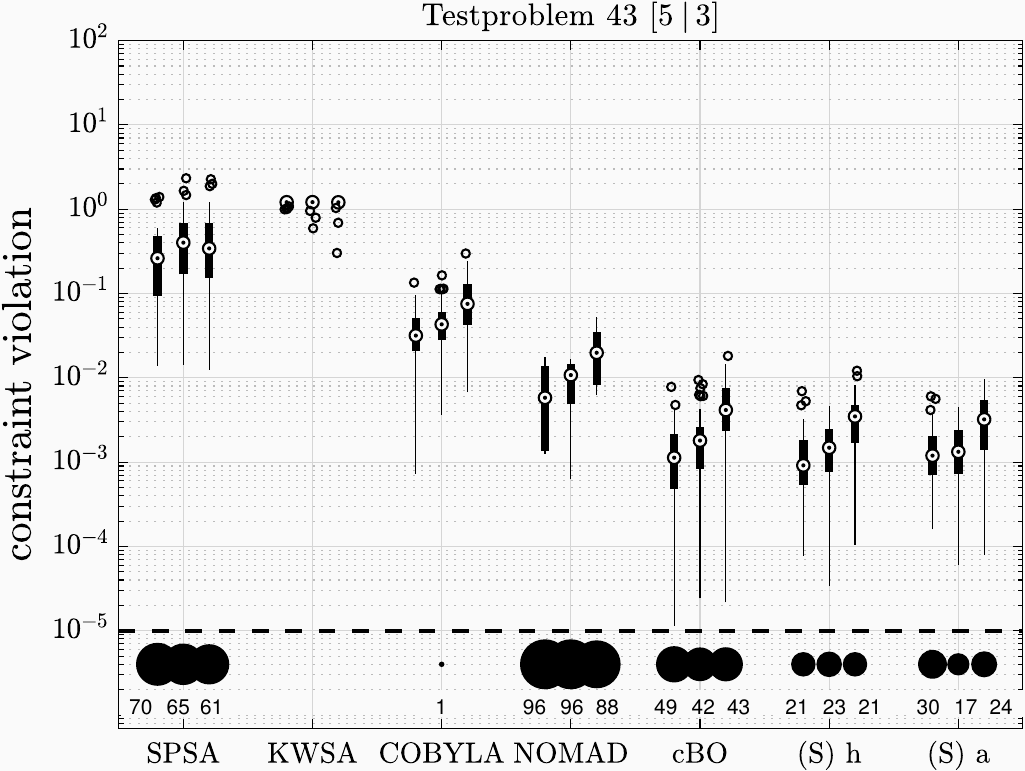}\hspace{2mm}
\includegraphics[width=0.3\textwidth]{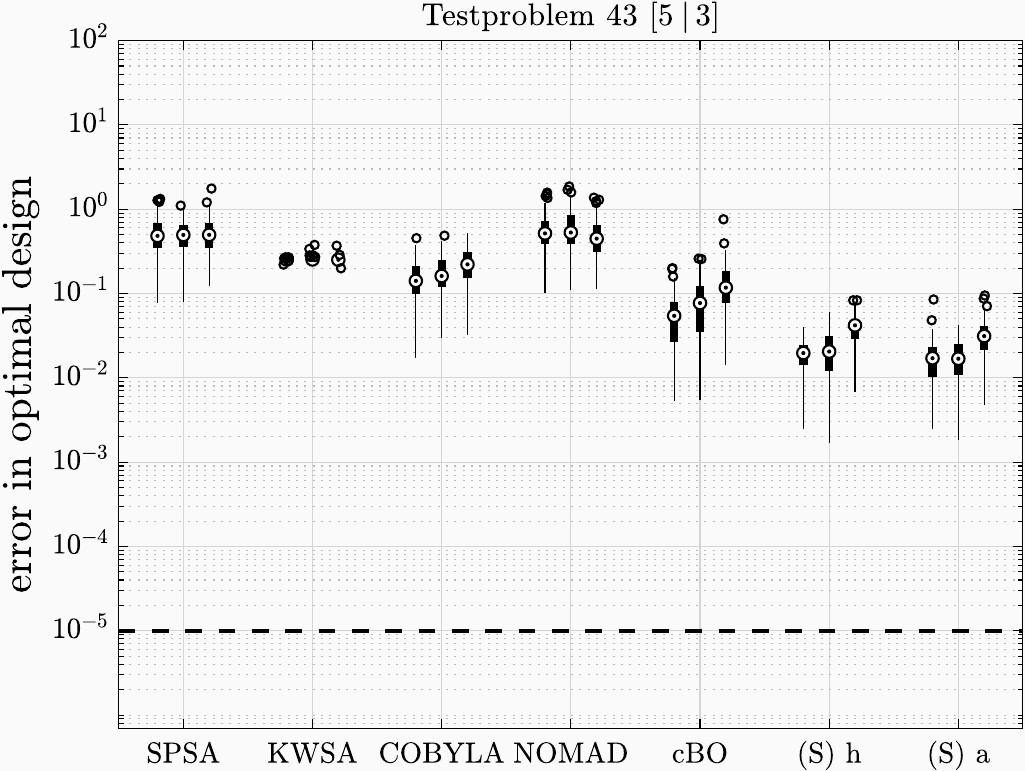}\\[3mm]
\includegraphics[width=0.3\textwidth]{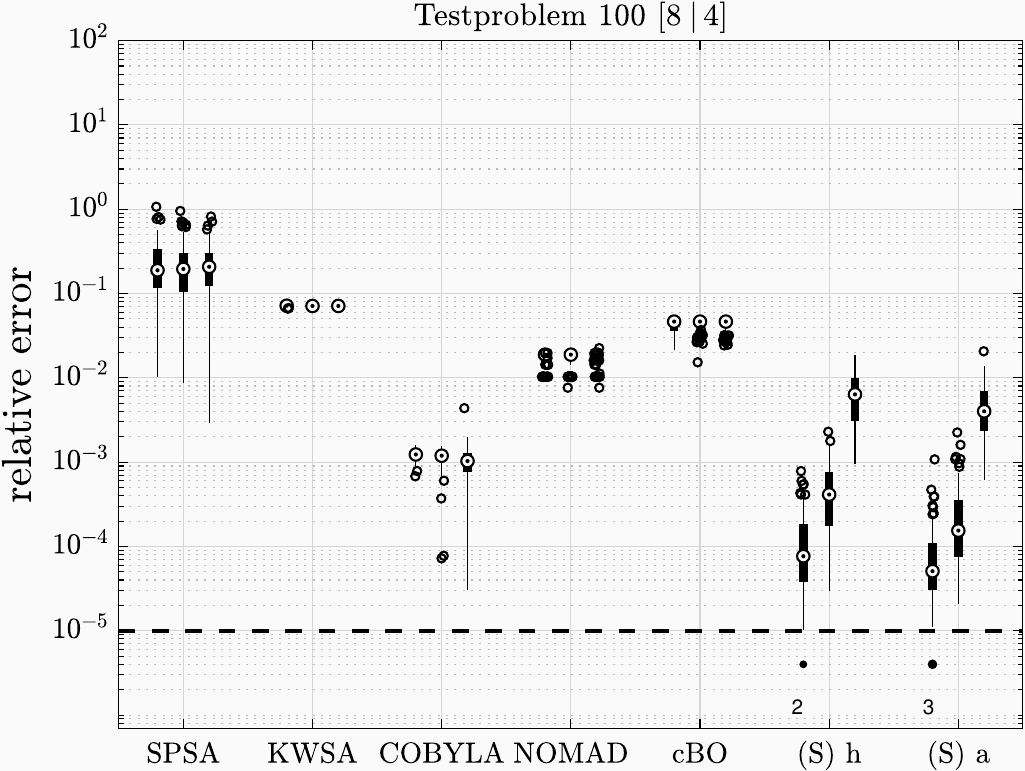}\hspace{2mm}
\includegraphics[width=0.3\textwidth]{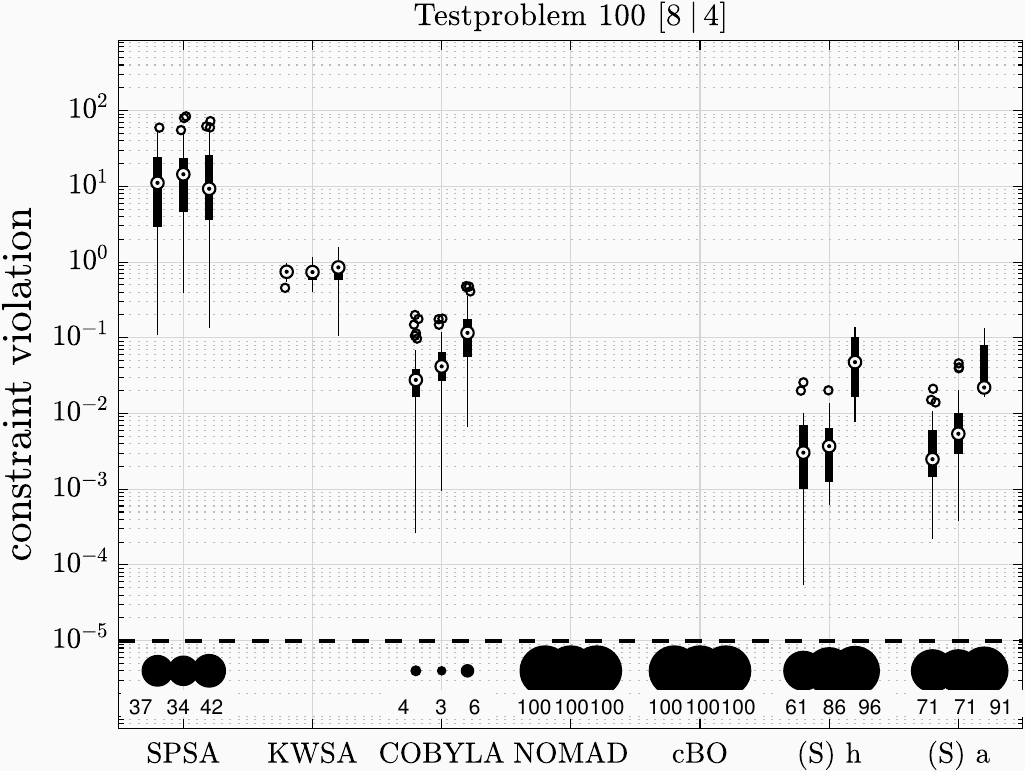}\hspace{2mm}
\includegraphics[width=0.3\textwidth]{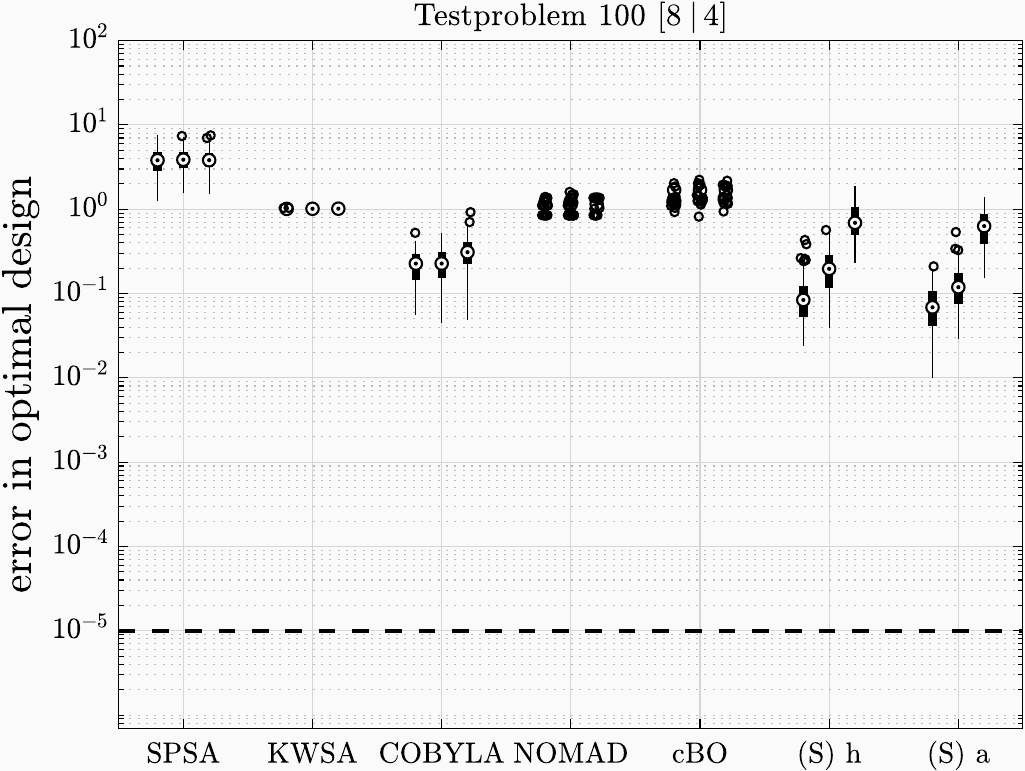}\\[3mm]
\includegraphics[width=0.3\textwidth]{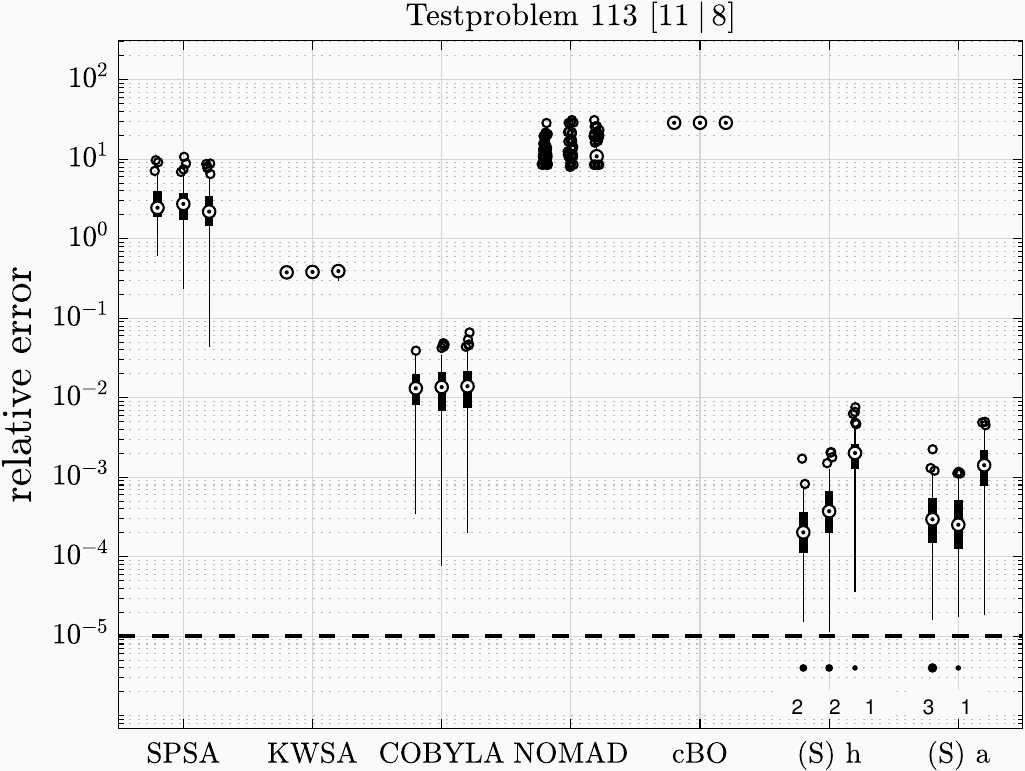}\hspace{2mm}
\includegraphics[width=0.3\textwidth]{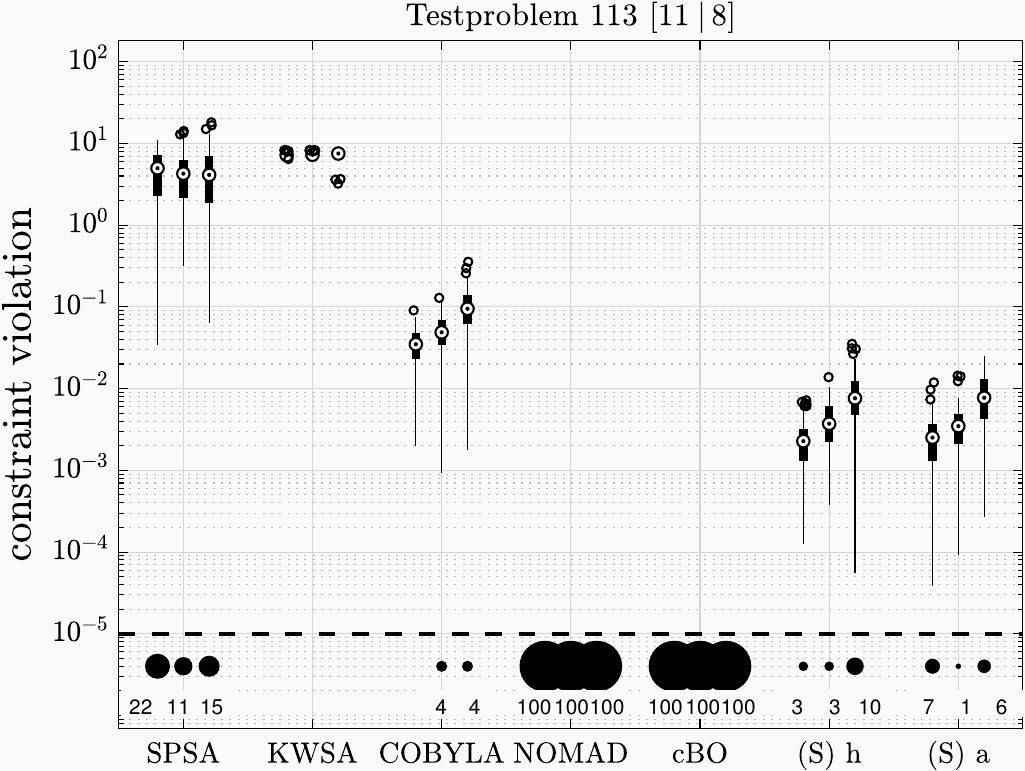}\hspace{2mm}
\includegraphics[width=0.3\textwidth]{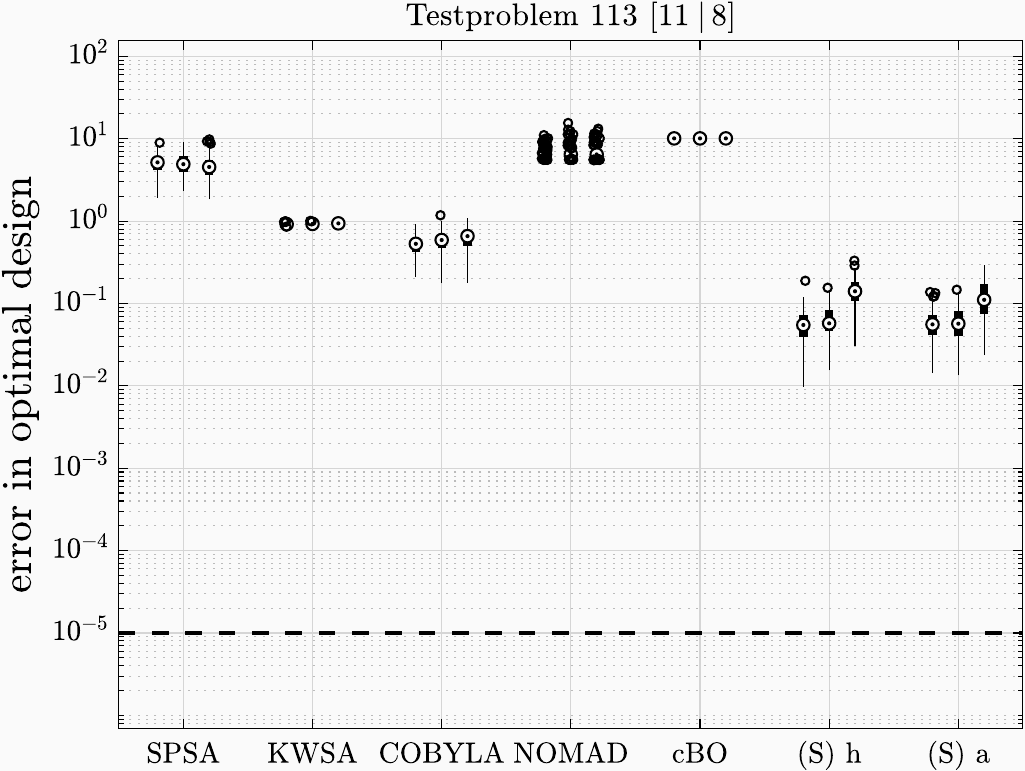}
\caption{Box plots of the errors in the approximated optimal objective values (left plots), the constraint violations (middle plots) and the $l_2$ distance to the exact optimal solution (right plots) of $100$ repeated optimization runs for the Schittkowski test problems number $29$, $43$, $100$, and $113$ for (\ref{eq:robust_cvar_value_Schittkowski_problems}). The plots show results of the exact objective function and constraints evaluated at the approximated optimal design computed by (S)NOWPAC, cBO, COBYLA, NOMAD, SPSA and KWSA. Thereby all errors or constraint violations below $10^{-5}$ are stated separately below the $10^{-5}$ threshold and the box plots only contain data above this threshold.}\label{fig:schittkowski_testset3a}
\end{center}
\end{figure}
\begin{figure}[!htb]
\begin{center}
\includegraphics[width=0.3\textwidth]{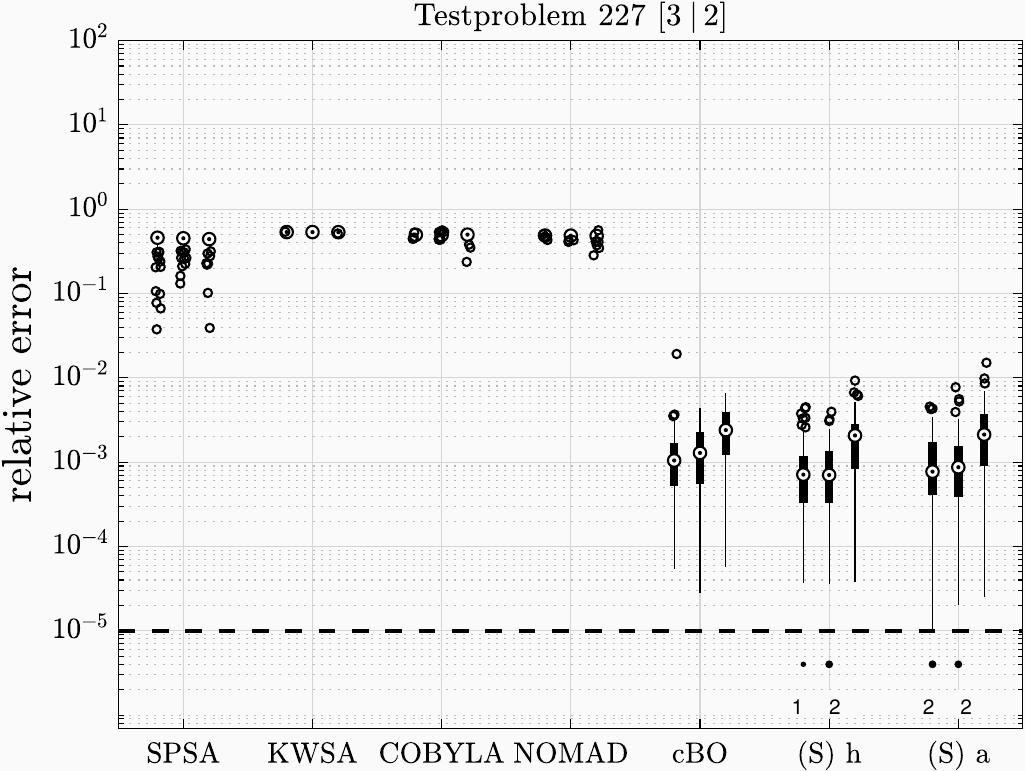}\hspace{2mm}
\includegraphics[width=0.3\textwidth]{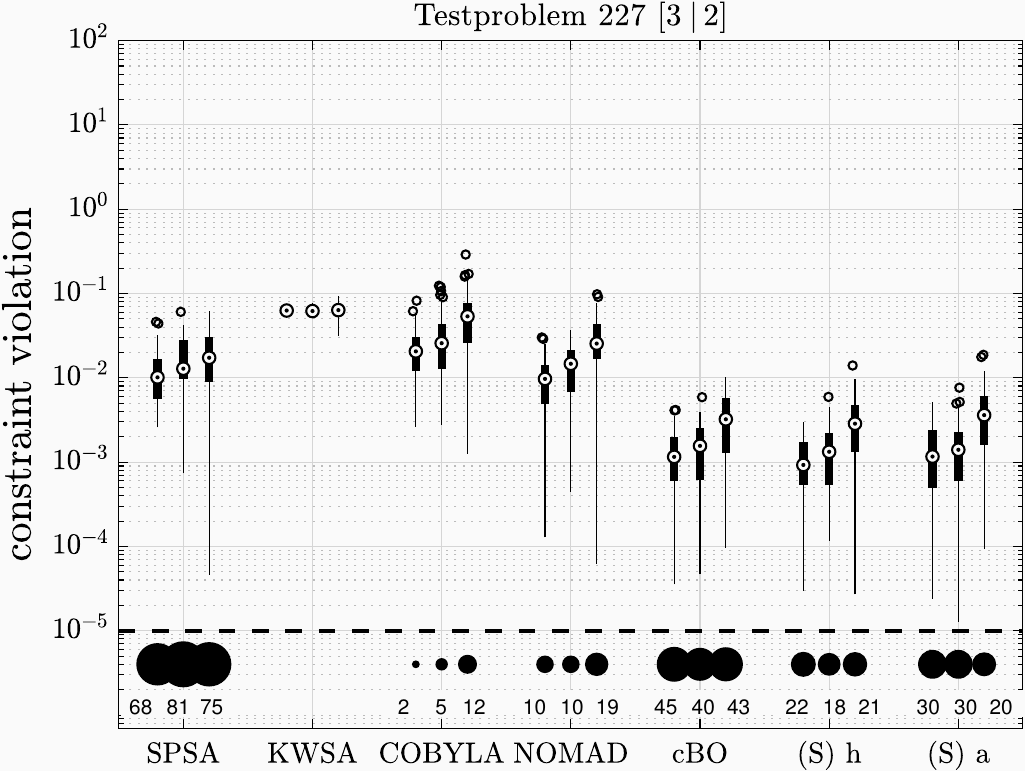}\hspace{2mm}
\includegraphics[width=0.3\textwidth]{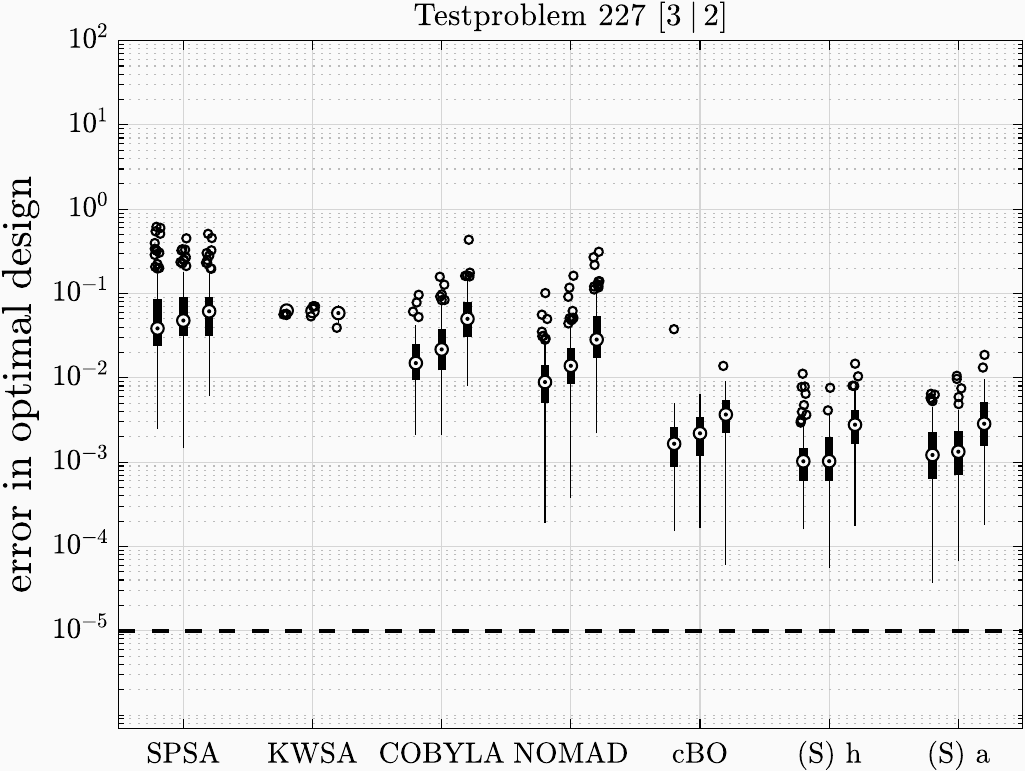}\\[3mm]
\includegraphics[width=0.3\textwidth]{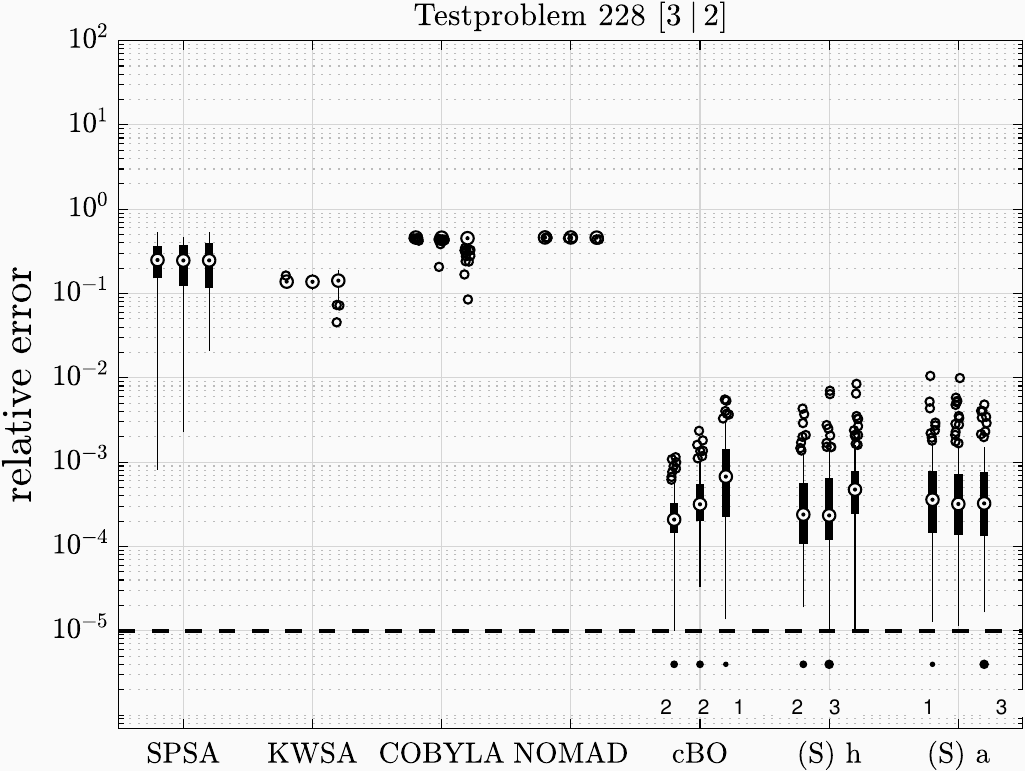}\hspace{2mm}
\includegraphics[width=0.3\textwidth]{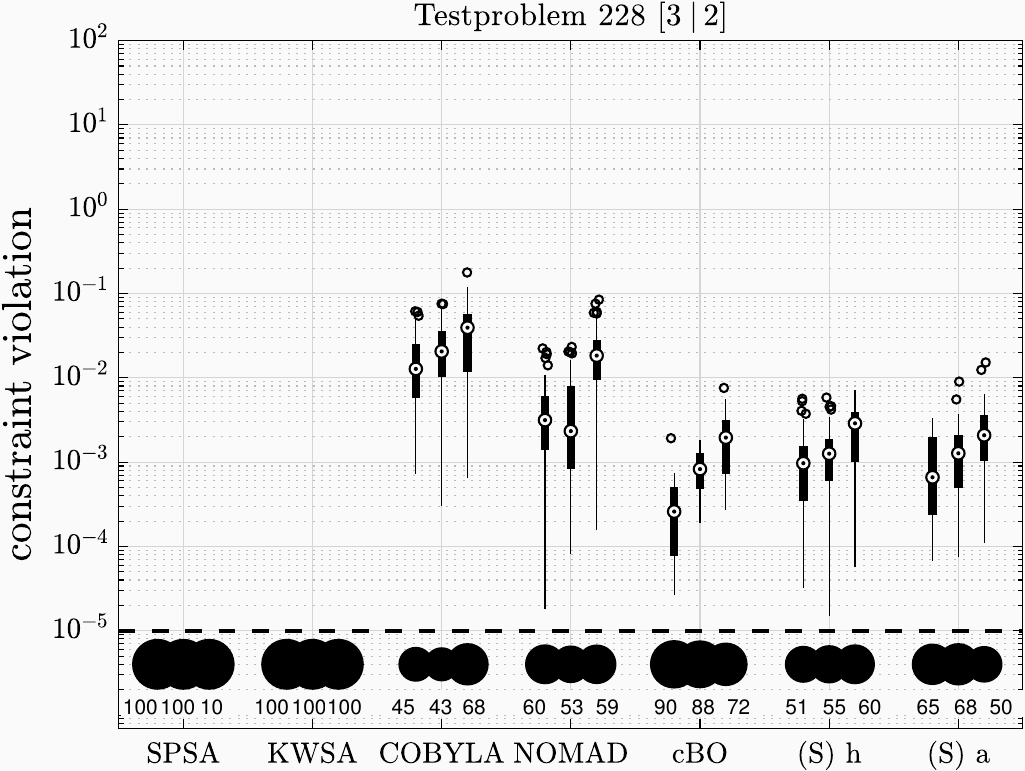}\hspace{2mm}
\includegraphics[width=0.3\textwidth]{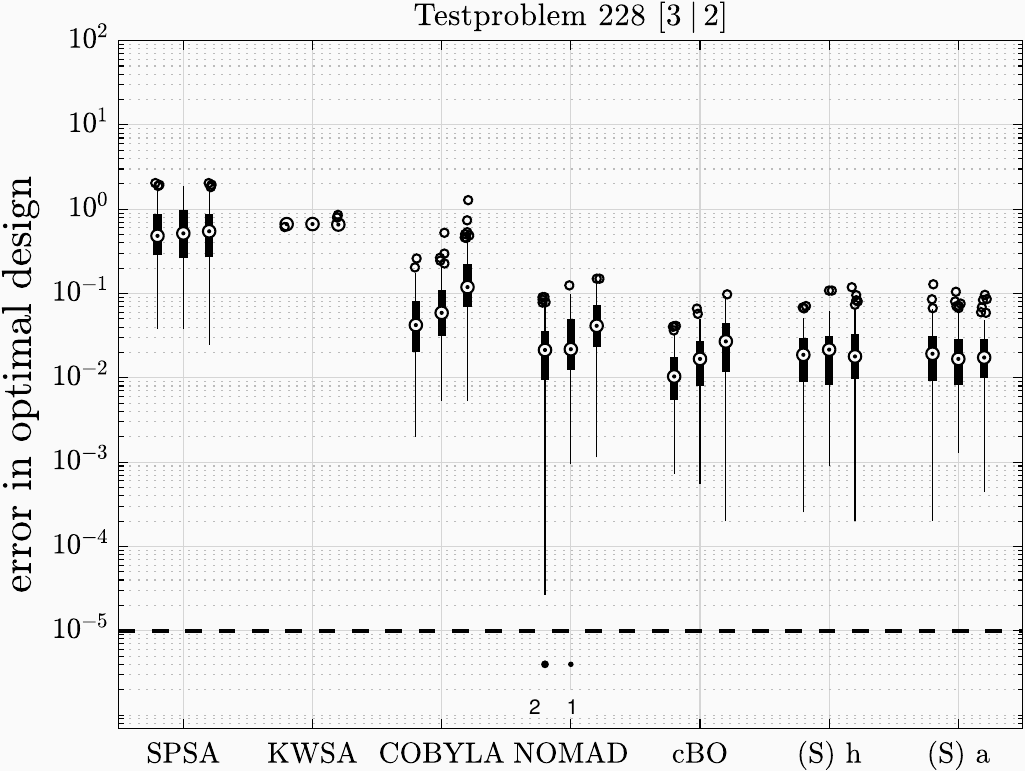}\\[3mm]
\includegraphics[width=0.3\textwidth]{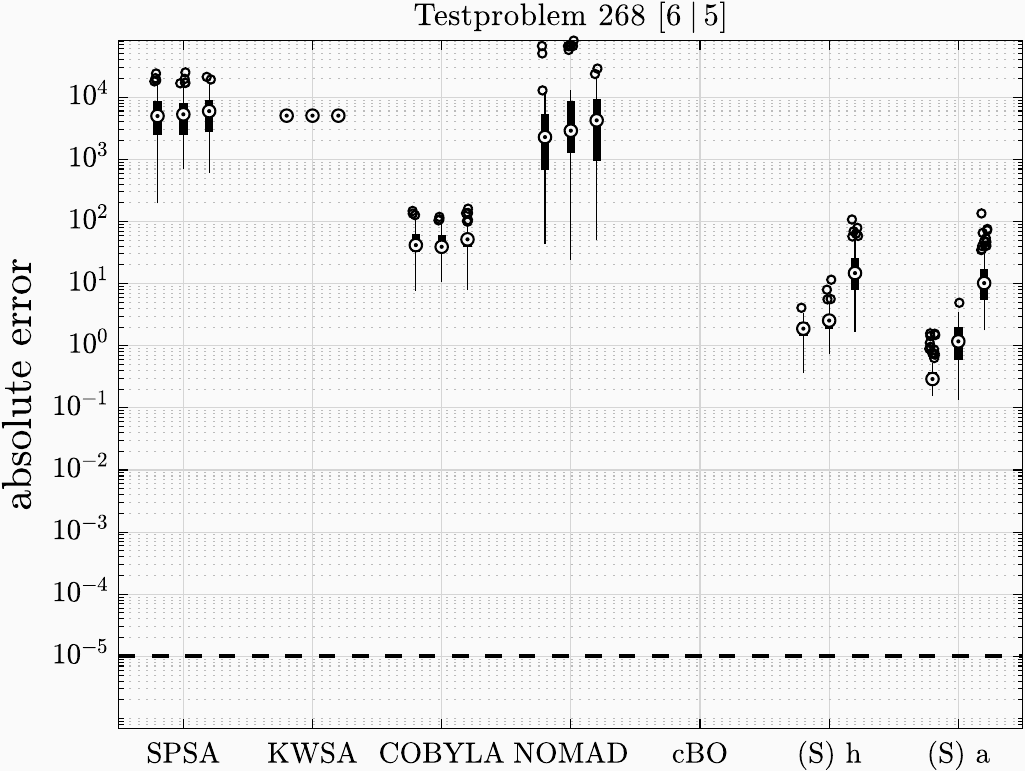}\hspace{2mm}
\includegraphics[width=0.3\textwidth]{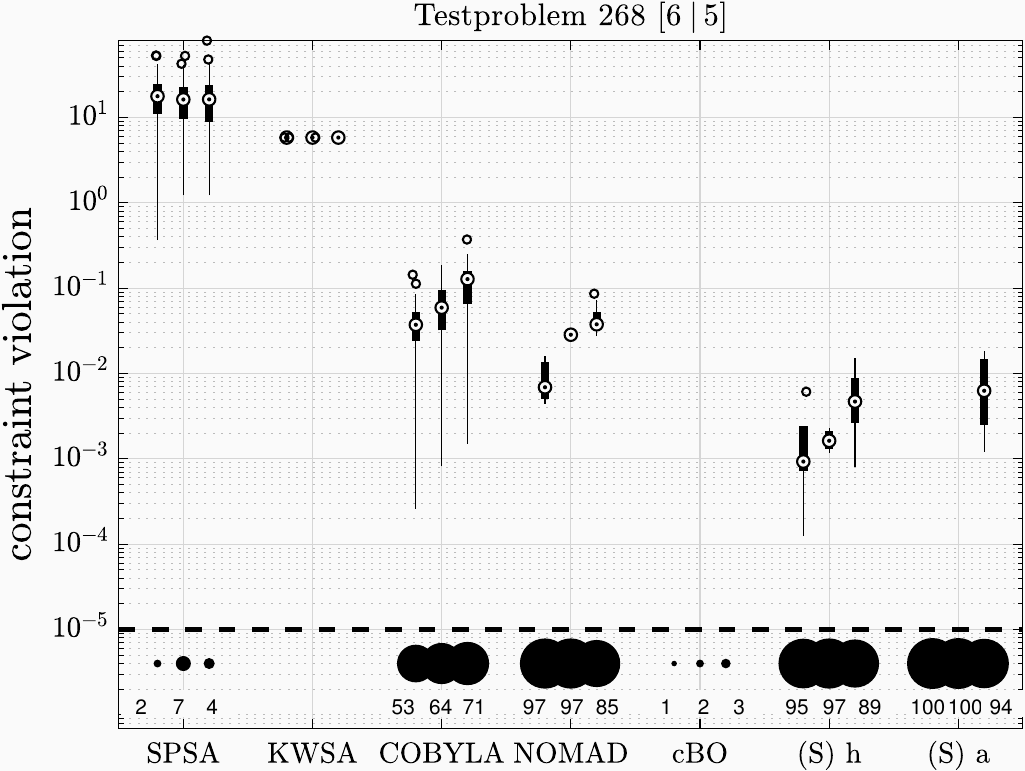}\hspace{2mm}
\includegraphics[width=0.3\textwidth]{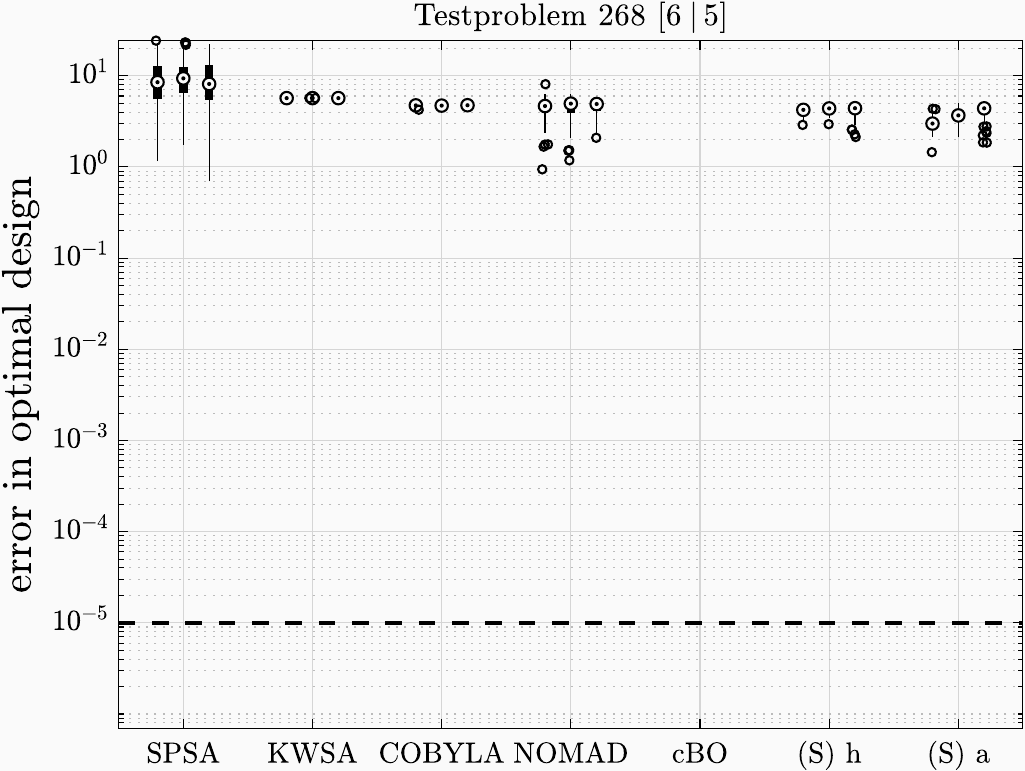}\\[3mm]
\includegraphics[width=0.3\textwidth]{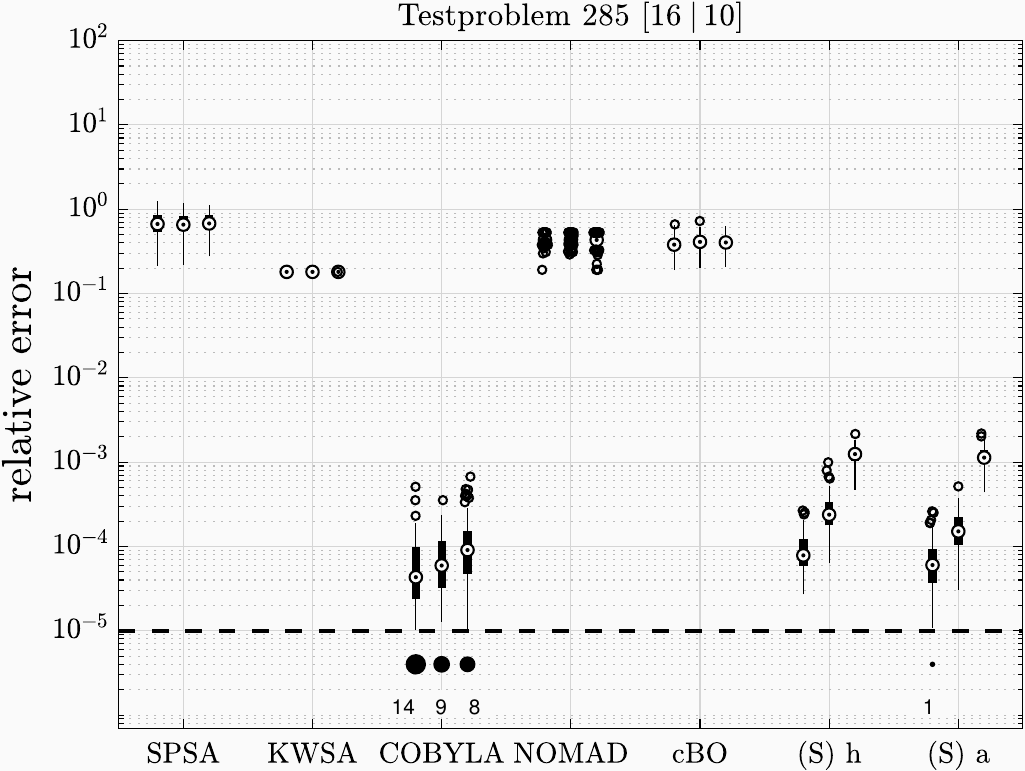}\hspace{2mm}
\includegraphics[width=0.3\textwidth]{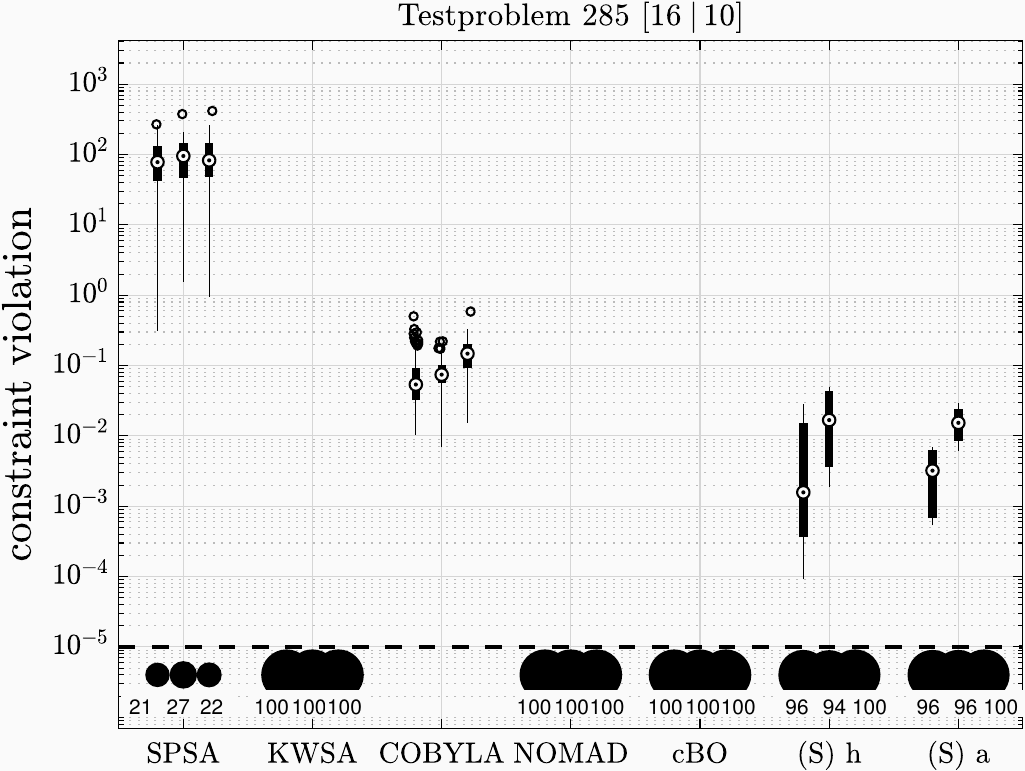}\hspace{2mm}
\includegraphics[width=0.3\textwidth]{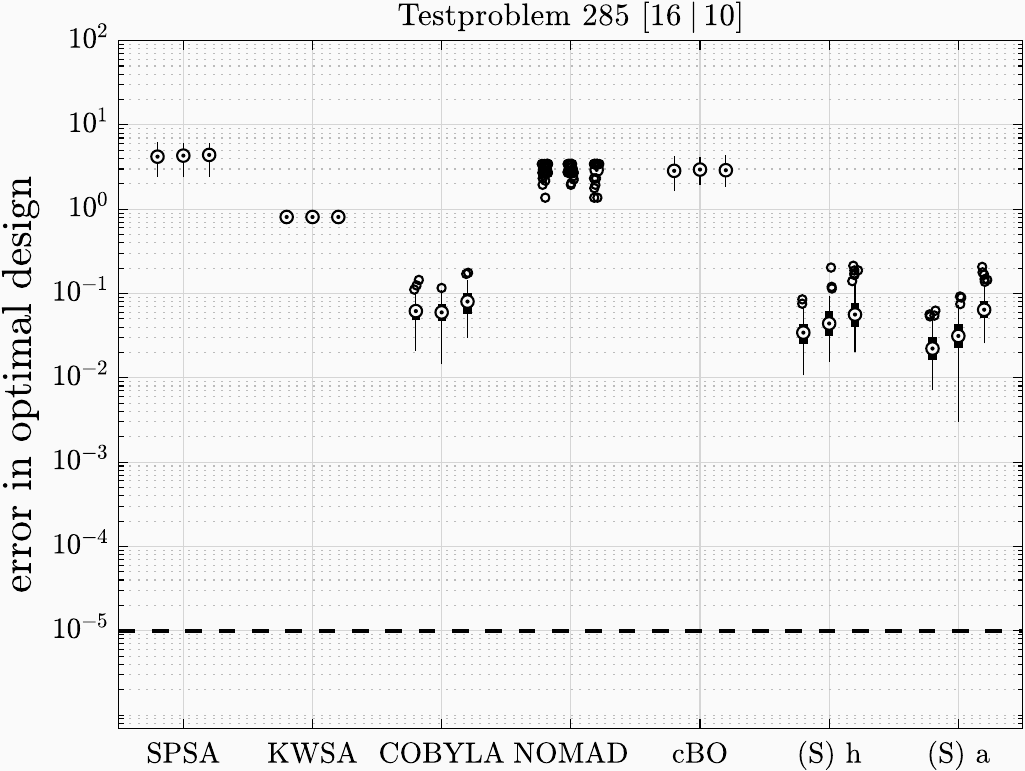}
\caption{Box plots of the errors in the approximated optimal objective values (left plots), the constraint violations (middle plots) and the $l_2$ distance to the exact optimal solution (right plots) of $100$ repeated optimization runs for the Schittkowski test problems number $227$, $228$, $268$, and $285$ for (\ref{eq:robust_cvar_value_Schittkowski_problems}). The plots show results of the exact objective function and constraints evaluated at the approximated optimal design computed by (S)NOWPAC, cBO, COBYLA, NOMAD, SPSA and KWSA. Thereby all errors or constraint violations below $10^{-5}$ are stated separately below the $10^{-5}$ threshold and the box plots only contain data above this threshold.}\label{fig:schittkowski_testset3b}
\end{center}
\end{figure}

\clearpage
\newpage
\listoftables
\newpage
\listoffigures
\end{document}